\documentclass[reqno, 11pt]{amsart}
\usepackage{amsmath}
\usepackage{amsthm}
\usepackage{amsfonts}
\usepackage{amssymb}
\usepackage{mathrsfs}
\usepackage{epsfig}
\usepackage{slashed}
\usepackage{mathtools}
\usepackage{enumerate}
\usepackage{fullpage}

\newcommand{\mb}{\mathbf}
\newcommand{\mc}{\mathcal}

\renewcommand{\Re}{\mathrm{Re}\,}
\renewcommand{\Im}{\mathrm{Im}\,}
\newcommand{\rg}{\mathrm{rg}\,}
\newcommand{\N}{\mathbb{N}}
\newcommand{\R}{\mathbb{R}}
\newcommand{\C}{\mathbb{C}}
\newcommand{\Z}{\mathbb{Z}}
\newcommand{\B}{\mathbb{B}}
\newcommand{\Oc}{\mathcal O}

\newcommand{\I}{\,\mathrm{i}\,}

\newcommand{\sgn}{\mathrm{sgn}\,}

\newcommand{\Lap}{\Delta}
\newcommand{\rad}{L^2(\R^7)}

\newcommand{\rst}[1]{\ensuremath{{\mathbin |}%
\raise-.5ex\hbox{$#1$}}} 

\newcommand{\norm}[1]{{\left\vert\kern-0.25ex\left\vert\kern-0.25ex\left\vert #1 
    \right\vert\kern-0.25ex\right\vert\kern-0.25ex\right\vert}}

\newtheorem{lemma}{Lemma}[section]
\newtheorem{theorem}[lemma]{Theorem}
\newtheorem{corollary}[lemma]{Corollary}
\newtheorem{proposition}[lemma]{Proposition}
\theoremstyle{remark}
\newtheorem{remark}[lemma]{Remark}

\theoremstyle{definition}
\newtheorem{definition}[lemma]{Definition}

\numberwithin{equation}{section}

\title[Stable blowup for the Yang-Mills heat flow]{Stable blowup for the supercritical Yang-Mills heat flow}

\author{Roland Donninger}
\address{Rheinische Friedrich-Wilhelms-Universit\"at Bonn, 
Mathematisches Institut,  Endenicher Allee 60, D-53115 Bonn, Germany}
\email{donninge@math.uni-bonn.de}

\author{Birgit Sch\"orkhuber}
\address{Universit\"at Wien, Fakult\"at f\"ur Mathematik,
Oskar-Morgenstern-Platz 1, A-1090 Vienna, Austria}
\email{birgit.schoerkhuber@univie.ac.at}

\begin{document}
\begin{abstract}
In this paper, we consider the heat flow for Yang-Mills connections on $\R^5 \times SO(5)$.
In the $SO(5)-$equivariant setting, the Yang-Mills heat equation reduces
to a single semilinear reaction-diffusion equation for which an explicit self-similar blowup solution 
was found by Weinkove \cite{Wei04}. We prove the nonlinear asymptotic stability of this solution under small perturbations.
In particular, we show that there exists an open set of initial conditions 
in a suitable topology such that the corresponding solutions blow up in finite time and converge to 
a non-trivial self-similar blowup profile on an unbounded domain. Convergence is obtained in suitable Sobolev norms
and in $L^{\infty}$. 
\end{abstract}

\maketitle

\section{Introduction}

\subsection{Equivariant Yang-Mills connections on $\R^d \times SO(d)$}
For $\mu = 1,\dots,d$, we consider mappings $A_{\mu}: \R^d \to  \mathfrak{so}(d)$, where 
$ \mathfrak{so}(d)$ denotes the Lie algebra of the Lie group $SO(d)$, i.e., $ \mathfrak{so}(d)$ can be considered as the set
of skew-symmetric $(d\times d)$-matrices endowed with
the commutator bracket. In the following, Einstein's summation convention is in force. For  
\[ F_{\mu \nu} := \partial_{\mu} A_{\nu} - \partial_{\nu} A_{\mu} + [ A_{\mu}, A_{\nu} ], \]
the Yang-Mills functional is defined as 
\begin{equation}\label{Eq:YangMillsFunc}
\mc F_A = \int_{\R^d} \mathrm{tr} (F_{\mu \nu} F^{\mu \nu}).
\end{equation}  

The Euler-Lagrange equations associated to this functional are given by
\begin{align}\label{Eq:YangMills}
\partial_{\mu}  F^{\mu \nu}(x) + [A_{\mu}(x),F^{\mu \nu}(x)] = 0
\end{align}
and solutions are referred to as Yang-Mills connections. We note that there is still a gauge freedom in this equation
which, however, is of no relevance here since we will only consider the equivariant case. 
A standard approach to find solutions to Eq.~\eqref{Eq:YangMills} is to add an artificial time dependence to the model which implies that 
Yang-Mills connections are static solutions to the corresponding heat flow equation
\begin{equation}\label{Eq:YMGradientflow}
\partial_t A_{\mu}(t,x) + \partial^{\nu}  F_{\mu \nu}(t,x)+ [A^{\nu}(t,x),F_{\mu \nu}(t,x)] = 0, \quad t > 0,
\end{equation}
for some initial condition $A_{\mu}(0,x) = A_{0\mu}(x)$. This model is usually referred to as the \textit{Yang-Mills heat flow}
for connections on the trivial bundle $\R^{d} \times SO(d)$. 
The equation enjoys scale invariance in the sense that
if $A_{\mu}$ is a solution, then, for any $\lambda > 0$,
\[A^{\lambda}_{\mu}(t,x) := \lambda A_{\mu}(\lambda^2 t, \lambda x),\]
solves Eq.~\eqref{Eq:YMGradientflow} with initial condition $A^{\lambda}_{\mu}(0,x) = \lambda A_{0 \mu}(\lambda x)$. 
Since the Yang-Mills functional in four space dimensions is invariant under this scaling transformation, the model is referred 
to as  \textit{critical} for $d= 4$.  Consequently, the Yang-Mills equation is
\textit{supercritical} for $d \geq 5$.
To simplify matters, we consider $SO(d)-$equivariant connections of the form
\begin{equation}\label{Def:Equiv}
A^{ij}_{\mu}(t,x) = u(t,|x|) ~ \sigma^{ij}_{\mu}(x) ,
\end{equation}
where $\sigma^{ij}_{\mu}(x) =  \delta^i_{\mu} x^j - \delta^j_{\mu} x^i$, 
for  $i,j = 1, \dots, d$, cf.~for example \cite{Grotowski01}, \cite{Gas02} and the references therein.
This ansatz reduces Eq.~\eqref{Eq:YMGradientflow} to a single equation for a radial function $u: [0,\infty) \times [0,\infty) \to \R$,
given by 
\begin{align} \label{Eq:EquivarEq}
\begin{split}
\partial_t u(t,r) -   \partial^2_r  u(t,r) - \frac{d+1}{r} \partial_r u(t,r)   + 3 (d-2)  u^2(t,r)  +  (d-2)  r^2 u^3(t,r) = 0, \quad t > 0,
\end{split}
\end{align} 
see \cite{Gas02} for a detailed derivation. It was shown by Grotowski \cite{Grotowski01} that this symmetry
is preserved by the flow. Hence, for equivariant initial data $A_{0 \mu}(x) =  u_0(|x|) ~ \sigma^{ij}_{\mu}(x)$
it suffices to consider Eq.~\eqref{Eq:EquivarEq} with $u(0,r) = u_0(r)$. 
The scale invariance of Eq.~\eqref{Eq:YMGradientflow} implies that
Eq.~\eqref{Eq:EquivarEq} is invariant under the transformation $u \mapsto u_{\lambda}$, where
\[ u_{\lambda}(t,r)= \lambda^2 u(\lambda^2 t, \lambda r) \]
and $\lambda > 0$. It was shown by Gastel  \cite{Gas02} that self-similar blowup solutions to Eq.~\eqref{Eq:EquivarEq} 
exist in dimensions $5 \leq d \leq 9$.
Explicit examples were found by Weinkove \cite{Wei04} and they are of the form 
\begin{align}\label{Eq:WeinkoveSol}
\mb {w}_T(t,r)  = \frac{1}{T-t} \mb W\left (\frac{r}{\sqrt{T-t}} \right), \quad \mb W(\rho) = - \frac{1}{a_1(d) \rho^2 + a_2(d)},
\end{align}
for some $T >0$, with constants
\begin{equation}\label{Def:ConstantsWeinkove}
a_1(d) =\frac {\sqrt{d-2}} {2 \sqrt{2}} , \quad a_2(d)=  \frac{1}{2} (6d-12-(d+2) \sqrt{2d-4} ),
\end{equation}
for $5 \leq d \leq 9$. By setting
\begin{equation}
\mb{A}_{\mu,T}(t,x) =  \mb{w}_T(t,|x|) ~ \sigma_{\mu}(x),
\end{equation} 
a one-parameter family of blowup solutions for Eq. \eqref{Eq:YMGradientflow} is obtained. 
These solutions are invariant under the natural scaling of the equation (up to a change of the blowup time) and obviously, 
they blow up in $L^{\infty}$ as $t \to T^{-}$. 
In this paper, we address the stability of these solutions under equivariant perturbations in five space dimensions. 
The main results are summarized below.

\subsection{Stable self-similar blowup for $d=5$} 
We consider the equation
\begin{align} \label{Eq:EquivarEq5}
\begin{split}
\partial_t u(t,r) =  \partial^2_r & u(t,r)  + \frac{6}{r} \partial_r u(t,r)  -  9 u^2(t,r)  -  3 r^2 u^3(t,r),
\end{split}
\end{align} 
for $r \in [0,\infty)$ and $t > 0$, subject to the initial condition $u(0,\cdot) = u_0$. We set $\mc E=C^{\infty}_{\mathrm{e},0}(\R)$, where 
\begin{align}\label{Def:EvenTestFunc}
C^{\infty}_{\mathrm{e},0}(\R) := \{ u \in C_0^{\infty}(\R):  u(x) =  u(-x) \}.
\end{align}
On $\mc E$ a norm is defined by $\| u \|_{\mc E}^2 : = \| u \|_{1}^2  +  \| u \|_{2}^2$, where 
\[ \|  u \|_{1}^2: = \int_0^{\infty} | r^3  u''(r) + 6 r^2  u'(r) |^2 dr, \]
\[ \|  u \|_{2}^2 := \int_0^{\infty} |r^3  u^{(4)}(r) + 12 r^2  u^{(3)}(r) + 24 r u''(r) - 24  u'(r)|^2 dr. \]

The following result shows that the blowup described by $\mb{A}_{\mu,T}$ is stable under equivariant
perturbations. 

\begin{theorem}\label{Th:Main}
Fix $T_0 > 0$. Then there are constants $\delta, K > 0$ such that for all real-valued functions $u_0 \in \mc E$ satisfying
\begin{align}\label{Eq:CondDataThMain}
\| u_0 - \mb{w}_{T_0}(0,\cdot) \|_{\mc E}   \leq \frac{\delta}{K}, 
\end{align}
the following holds. 
\begin{enumerate}
\item There is a $T \in [T_0- \delta, T_0 + \delta]$ such that a unique classical solution $u(t,\cdot)$ to
Eq.~\eqref{Eq:EquivarEq5} exists for all $t \in (0,T)$ with $u(0,\cdot) = u_0$. 
\item At $t=T$, the solution blows up at the origin and converges to $\mb w_{T}$ according to 
\begin{align}\label{Eq:ConvSolHnorm}
 \frac{\|u(t,\cdot) - \mb{w}_{T}(t,\cdot)\|_1 }{ \|\mb{w}_{T}(t,\cdot)\|_1 } \lesssim (T-t)^{\frac{1}{150}}, 
\quad    \frac{\|u(t,\cdot) - \mb{w}_{T}(t,\cdot)\|_2 }{\|\mb{w}_{T}(t,\cdot)\|_2} \lesssim (T-t)^{\frac{1}{150}}.
\end{align}
\item Furthermore, we have convergence in $L^{\infty}(\R^+)$, i.e.,
\begin{align}\label{Eq:ConvSolLinfty}
\frac{ \| u(t,\cdot) - \mb{w}_{T}(t,\cdot) \|_{L^{\infty}(\R^+)}}{ \|\mb{w}_{T}(t,\cdot) \|_{L^{\infty}(\R^+)}} 
 \lesssim (T-t)^{\frac{1}{150}}.
\end{align}
\end{enumerate}
\end{theorem}

Some remarks are in order. 
\begin{remark}
\begin{enumerate}[(i)]
\item Although the problem is posed on $\R^{d}$, $d =5$, the effective dimension of the Laplacian in Eq.~\eqref{Eq:EquivarEq}
is $d+2$. From a mathematical point of view it is therefore reasonable to consider $u$ as a radial function
on $\R^7$. Obviously, 
\[ \|u\|_{1} \simeq \|\Lap u(|\cdot|) \|_{L^2(\R^7)}, \quad \|u\|_{2} \simeq \|\Lap^2 u(|\cdot|) \|_{L^2(\R^7)}.\]
\item In all of the above bounds the left hand side is normalized to
the behavior of $\mb{w}_{T}$ in the respective norm, i.e., we obtain convergence of the solution relative 
to $\mb{w}_{T}$. However, the given rate of convergence is not sharp.
\item The smoothness assumptions on the initial data seem to be quite restrictive. In fact, the result holds
for a much larger class of initial conditions, see the discussion in Section \ref{Sec:OutlineProof}. 
A more general version for solutions that satisfy the equation in a suitable weak sense is given in Theorem \ref{Th:MainSim} below.
\item The restriction to $d=5$ is by no means crucial and our techniques easily extend to the cases $d=7$ and $d=9$. For the sake of simplicity, however, we only consider the case $d=5$.
\end{enumerate}
\end{remark}

We note that our approach is very robust since we do not make use of any Lyapunov functionals or monotonicity formulas. Instead, we rewrite Eq.~\eqref{Eq:EquivarEq}
in similarity coordinates and study perturbations around $\mb{w}_T$ by means of strongly continuous semigroups, 
operator theory and spectral analysis. This is described in more detail in Section \ref{Section:Formulation_SimVar}.

\subsection{Acknowledgments}
Roland Donninger is supported by the Alexander von Humboldt Foundation via
a Sofja Kovalevskaja Award endowed by the German Federal Ministry of Education
and Research. Birgit Sch\"orkhuber is supported by the Austrian Science Fund
(FWF) via the Hertha Firnberg Program, Project Nr. T 739-N25. Partial support by the Deutsche Forschungsgemeinschaft 
(DFG), CRC 1060 'The Mathematics of Emergent Effects', is also gratefully acknowledged.

\medskip
The second author would like to thank Herbert Koch and the group Analysis and Partial Differential Equations for the warm hospitality 
and the inspiring atmosphere
during a three-months stay at the Mathematical Institute of the University of Bonn, where this work has been initiated.

\subsection{A brief history of the problem and known results} 
In this paper, we are studying Yang-Mills connections on the trivial bundle $\R^d \times SO(d)$. 
However, the theory is usually formulated in the much more general
language of principal $G$-bundles over manifolds, where oftentimes $G$ is a compact Lie group
such as $SO(d)$ or $SU(d)$. 
The study of Yang-Mills functionals in mathematics
was initiated in the 1980's and triggered profound developments
in differential geometry, see for example \cite{Donaldson_MathUse} for a review and 
a list of references. The corresponding heat flow equation has been studied extensively 
in various geometric settings. 
Over closed Riemannian manifolds of dimension $d=2,3$, R{\aa}de \cite{Rade09} showed
that solutions exist globally in time and converge to a Yang-Mills connection as $t \to \infty$.

In the critical case $d=4$, the existence of 
global weak solutions has been investigated for example in \cite{StruweYM1994}, \cite{Schlatter} and \cite{Kozono1995}.
For $SU(2)-$bundles over $\R^4$, global existence of smooth solutions was established by Schlatter, 
Struwe and Tahvildar-Zadeh  \cite{SchlStruwZah98} in the equivariant setting.
In fact, they studied a slightly different 
formulation of Eq.~\eqref{Eq:EquivarEq} that is obtained by setting 
$h(t,r) := - r^2 u(t,r)$, which yields
\begin{align}\label{Eq:YMVer_3dLap}
h_t -  h_{rr} - \frac{d-3}{r} h_r   + \frac{(d-2)h(h-1)(h-2)}{r^2} =0 .
\end{align}

In supercritical dimensions, blowup in finite time for the YM heat
flow over $\mathbb{S}^d$, $d \geq 5$, is due to 
Naito \cite{Naito94}. Grotowski \cite{Grotowski01}
considered the problem in the setting of a trivial $SO(d)$--bundle
over $\R^d$ for $d \geq 5$.
By constructing suitable subsolutions to Eq.~\eqref{Eq:YMVer_3dLap}, he proved singularity formation
from smooth initial data. As mentioned above, the existence of self-similar blowup solutions to Eq.~\eqref{Eq:EquivarEq} was first established by Gastel \cite{Gas02} for $5 \leq d \leq 9$.
Weinkove \cite{Wei04} investigated the nature of singularities over compact Riemannian manifolds under the assumption 
of $\textit{type I}$ blowup
(which implies a certain upper bound on the blowup rate of the curvature), and it was found that
locally around the blowup point solutions converge in a suitable sense to so-called \textit{homothetically shrinking solitons},
which are also referred to as \textit{YM-solitons}. These objects correspond to solutions of the YM heat flow on the trivial bundle over $\R^d$, which is our main motivation to study the problem in this geometrical setting.
Moreover, Weinkove gave explicit examples of such solitons on $\R^d \times SO(d)$, $5 \leq d \leq 9$, see Eq.~\eqref{Eq:WeinkoveSol}.
Very recently, a description of general blowup solutions for the YM heat flow over closed Riemannian manifolds
was obtained by Kelleher and Streets \cite{KelStreets2016} for $d \geq 4$ and it was shown that singularities can be described
either by suitably rescaled Yang-Mills connections, i.e., by static solutions, or by YM-solitons.

The results of Weinkove have raised interest in the stability of YM-solitons in recent years and
notions of variational stability have been introduced by Kelleher and Streets \cite{KelStreets2014} as well
as by Chen and Zhang \cite{ChenZhang2015}. However, to the best of our knowledge no rigorous proof on the stability
of the Weinkove solution given in Eq.~\eqref{Eq:WeinkoveSol} has been obtained so far and Theorem \ref{Th:Main} is the first
result in this direction. 

We note that in higher dimensions, $d \geq 10$, the existence of self-similar solutions to Eq.~\eqref{Eq:EquivarEq} 
was excluded by Bizo{\'n} and Wasserman \cite{BizonWasserman2015}
and it is expected that in this case the generic blowup is of 
\textit{type II}, see also the discussion on the harmonic map heat flow below. 

\subsection{Related problems} 

\subsubsection{Supercritical harmonic maps heat flow}
The above model bears many similarities with the heat flow
for harmonic maps from $\R^d \to \mathbb{S}^d$, which in co-rotational symmetry 
reduces to 
\begin{align}\label{Eq:HarmMaps}
h_t -  h_{rr} - \frac{d-1}{r} h_r   + \frac{(d-1)\sin(2h)}{ 2 r^2} = 0.
\end{align}
The model is supercritical for $d \geq 3$. For $3 \leq d \leq 6$, 
Fan \cite{Fan99} has constructed a family of self-similar blowup solutions to
Eq.~\eqref{Eq:HarmMaps}. According to the numerical work of 
Biernat and Bizo{\'n} \cite{Bizon_Biernat2011},  the ground state of this family describes the behavior 
of generic blowup solutions. In this scenario, the gradient of the solution diverges, while the solution itself remains bounded. 
In particular, self-similar solutions blow up in a type I fashion, which is defined by the bound
$\sup_{t < T} | \sqrt{T-t}  ~ \partial_r h(0,t)| < \infty$.
In contrast to the YM heat flow, none of the self-similar blowup solutions 
is known in closed form and their stability is largely open.
We are convinced that the techniques developed in this paper can be used to study the corresponding problem
for Eq.~\eqref{Eq:HarmMaps}, provided one has a sufficiently good approximation to
a self-similar blowup profile. 

In analogy to the YM heat flow, self-simlar blowup is excluded
for $d \geq 7$, see \cite{BizonWasserman2015}. In this regime, Biernat and Seki \cite{BiernatSeki16} have recently 
constructed explicit examples for type II blowup solutions. It is most likely that a similar result for the YM heat
flow can be obtained in dimensions $d \geq 10$.

We note that the problem of non-uniqueness of weak solutions 
to Eq.~\eqref{Eq:HarmMaps} in the supercritical case (which is naturally tied 
to the question of continuation after the blowup) has recently been addressed by
Germain, Ghoul and Miura \cite{GermainGhoulMiura16}.

Furthermore, it is worth mentioning that in the critical case $d=2$, finite-time blowup for solutions of Eq.~\eqref{Eq:HarmMaps}
has been proved by Chang and Ding \cite{ChangDing92} which 
contrasts the result of \cite{SchlStruwZah98} for the YM heat flow, see also \cite{GrotShat07} for a comparison of 
the two models. Stable type II blowup for $d=2$ is due to Rapha{\"e}l and Schweyer \cite{SchweyerRaph_2013} and involves dynamically
rescaled static solutions of Eq.~\eqref{Eq:HarmMaps} as asymptotic blowup profiles, see also \cite{SchweyerRaph_2014}.

\subsubsection{The supercritical heat equation}
The formation of singularities for the nonlinear heat equation 
\begin{align}\label{Eq:NLH}
 \partial_t u(t,x) - \Delta u(t,x) = |u(t,x)|^{p-1} u(t,x)
 \end{align}
on $\R^d$ (or on some bounded domain $\Omega \subset \R^d$)
has been investigated extensively in the past decades and reviewing the vast body of literature on the subject is clearly beyond the scope of this introduction. At the very least, however, we would like to briefly mention some results in the supercritical case $p > \frac{d+2}{d-2} =: p_c$, where apart from the usual ODE blowup solution $u(t,x)= \pm \kappa (T-t)^{-\frac{1}{p-1}}$, $\kappa = (p-1)^{-\frac{1}{p-1}}$, infinitely many (presumably unstable) nontrivial self-similar blowup solutions exist, giving rise to type I blowup behavior. In a series of papers,  Matano and Merle \cite{MatanoMerle2004}, \cite{MatanoMerle2009}, \cite{MatanoMerle2011} investigated various aspects of singularity formation for radial solutions. Among many other  results, they excluded the existence of type II blowup for a certain range of nonlinearities $p < p^*$. For $p > p^*$ such solutions were in fact constructed in \cite{HerreroVelazquez94}, \cite{HerreroVelazquezUnpub} and \cite{Mizoguchi2004}, see also \cite{Collot2016NLH} for a remarkable result in the non-radial case. The role of the ODE blowup solution as a generic local blowup profile for $p > p_c$ has been established in \cite{MatanoMerle2011}.
In the non-radial case, another line of investigation was pursued by Blatt and Struwe \cite{StruweBlatt2015a}, \cite{StruweBlatt2015b}, \cite{StruweBlatt2015c}, who studied global existence and blowup of solutions to Eq.~\eqref{Eq:NLH} for $p > p_c$ by using Morrey spaces and monotonicity formulas.

Finally, we would like to mention a recent paper by Tayachi and Zaag \cite{TayachiZaag2015}, who constructed a stable
(non-selfsimilar) blowup solution to a version of Eq.~\eqref{Eq:NLH} with an additional nonlinear gradient term.

In the context of Eq.~\eqref{Eq:NLH} it is also worth mentioning that one of the key steps in our approach is the reformulation of Eq.~\eqref{Eq:EquivarEq5} in similarity coordinates. This goes back to the work of Giga and Kohn \cite{GigaKohn1985} on blowup for Eq.~\eqref{Eq:NLH} in the subcritical case and also plays an important role in many of the above cited works. However, as already mentioned above, our methods do not rely on parabolic techniques and are thus very robust concerning their application to other types of nonlinear evolution equations, see below.

\subsubsection{Yang-Mills wave equation on $\R^{5+1}$} 
Yang-Mills theory has its true origins in particle physics, where it provides a standard framework 
to formulate non-abelian gauge theories.
There, the Yang-Mills functional is usually considered 
over $\R^{d+1}$-Minkowski space and Eq.~\eqref{Eq:YangMills} corresponds to a nonlinear time dependent PDE. 
%Hier noch Referenzen einfuegen?? T
In a similar geometric setting as described above, there is a 'wave counterpart' of Eq.~\eqref{Eq:EquivarEq} 
(with $\partial_t u$  replaced by $\partial_t^2 u$), cf. for example \cite{B_YM2002}. 
For $d=5$, singularity formation for the resulting  YM wave equation 
is known and the existence of infinitely many self-similar blowup solutions has been proved by  Bizo{\'n} \cite{B_YM2002},
who also found a closed form expression for the ground state solution. The stability of this object has been established
by the first author in \cite{Donninger14_YM} under a spectral assumption which, however, was recently removed in a joint 
paper of the first author with Costin, Glogi\'c and Huang \cite{CostinDonn}. The proof in \cite{Donninger14_YM} is 
based on a canonical method that was developed by the authors to investigate stable self-similar blowup for wave equations, 
cf. similar results for supercritical wave maps  \cite{DonSchAic12} and wave equations with focusing power-nonlinearities,
\cite{DonSch13}, \cite{DonSch14b}, \cite{DonSch15a}.

In this paper, we generalize our method
to parabolic problems. In stark contrast to our results on wave equations, which are restricted to a backward lightcone, 
we are here able to treat the problem on the whole space. In particular, we construct open sets of initial data such that the 
corresponding solutions converge to a nontrivial self-similar blowup profile on an unbounded domain. 
%To the best of our knowledge, this is the first result in this direction for an energy supercritical model. 

%In five dimensions case an 'energy' functional for Eq.~\eqref{Eq:EquivarEq} is given by
%\begin{align*}
%\mc E_u(t) = \int_0^{\infty} r^{2} \left ( [r^2 u(t,r) + 2 r u(t,r)]^2 + \frac{3 r^4  u(t,r)^2 [r^2 u(t,r) +2]^2}{2 r^2}  \right) 
%\end{align*}
%and one can easily check that $\mc E'_u(t) \leq 0$ for classical solution that decay fast enough at infinity.
%Under scaling this quantity behaves like 
%\[\mc E_{u_{\lambda}}(t) = \lambda^{-1} \mc E_{u}(t) \]
%which classifies the problem as 'energy supercritical'.

\subsection{Formulation of the problem in similarity coordinates}\label{Section:Formulation_SimVar}

We fix $T_0 > 0$ and write the initial condition as
\begin{align}\label{Eq:OrigData}
u(0,r) =  \mb{w}_{T_0}(0,r) + v_0(r),
\end{align}
for $r \in [0,\infty)$, where $\mb {w}_{T_0}$ denotes the self-similar blowup solution with fixed blowup time $t = T_0$, see 
Eq.~\eqref{Eq:WeinkoveSol}. We emphasize that this is no restriction on the 
data since $v_0$ is a free function. 
The aim is to show that if $v_0$ is small in a suitable sense, then the corresponding solution to 
Eq.~\eqref{Eq:EquivarEq5} blows up and  converges
to $\mb {w}_T$, for some suitable $T > 0$, that will in general depend on $v_0$.
Therefore, the blowup time $T$ enters the analysis as a free parameter that will be fixed only
at the very end of the argument. Having this in mind, we introduce similarity coordinates
 $(\tau, \rho) \in [0,\infty) \times [0,\infty)$ defined by 
\[ \tau = - \log(T-t) + \log T, \quad \rho = \frac{r}{\sqrt{T-t}}, \]
for $T > 0$. By setting 
\[  \psi(- \log(T-t) + \log T , \tfrac{r}{\sqrt{T-t}})  :=  (T-t)  u(t, r), \]
the initial value problem given by Eq.~\eqref{Eq:EquivarEq5} and Eq.~\eqref{Eq:OrigData} transforms into
\begin{align}\label{Eq:SelfSim}
\begin{split}
\partial_{\tau}  \psi(\tau, \rho) =   \partial^2_{\rho} &  \psi(\tau, \rho)  + \frac{6}{\rho} \partial_{\rho} \psi(\tau, \rho)  
- \frac{1}{2} \rho \partial_{\rho}   \psi(\tau, \rho)    -  \psi(\tau, \rho)   - 9   \psi^2(\tau, \rho)  - 3 \rho^2   \psi^3(\tau, \rho) 
\end{split}
\end{align}
for $\tau > 0$ with initial condition
\[\psi(0,\rho) = T u(0,\sqrt{T} \rho) = \tfrac{T}{T_0} \mb W \left (\tfrac{ \sqrt{T} }{\sqrt{T_0}} \rho \right )   + T v_0 ( \sqrt{T} \rho ). \]
We note that the parameter $T$ does not appear in the equation itself but only in the initial condition.
The Weinkove solution
\begin{align}\label{Eq:GroundState}
\mb W(\rho) = - \frac{1}{a_1 \rho^2 + a_2},
\end{align}
with 
\begin{align}\label{Eq:ConstantsWein5d}
a_1 = \frac12 \sqrt{\frac{3}{2}}, \quad a_2 = \frac12 (18 - 7 \sqrt{6}),
\end{align}
is a static solution to Eq.~\eqref{Eq:SelfSim}.

The differential operator on the right hand side of Eq.~\eqref{Eq:SelfSim} has a natural extension to 
$\R^7$. In fact,  Eq.~\eqref{Eq:SelfSim} can be written as 
\begin{align}\label{Eq:SelfSimNonSym}
\begin{split}
\partial_{\tau}  \Psi(\tau, \xi) =  \Lap \Psi (\tau, \xi) - \tfrac{1}{2} \xi \cdot \nabla \Psi (\tau, \xi) - \Psi (\tau, \xi)   -
9   \Psi^2 (\tau, \xi)  -  3 |\xi|^2 \Psi^3 (\tau, \xi), \quad \tau > 0,
\end{split}
\end{align}
for $\xi \in \R^7$ and a radial function $\Psi(\tau, \xi) = \psi(\tau, |\xi|)$. Hence, we study
Eq.~\eqref{Eq:SelfSimNonSym} with initial data of the form  
\[ \Psi(0, \xi)  = \tfrac{T}{T_0} \mb W \left (\tfrac{ \sqrt{T} }{\sqrt{T_0}} |\xi| \right )   + T v_0 ( \sqrt{T} |\xi| ).\]
In a more abstract way, the problem can be formulated as 
\begin{align}\label{Eq:YM_Abstract}
\begin{split}
\frac{d}{d \tau} \Psi(\tau)  & = L_0 \Psi(\tau)  + F(\Psi(\tau)) ,\quad \tau > 0, \\
\Psi(0) & = \Psi^{T}_0,
\end{split}
\end{align}
where $L_0$ represents the linear differential operator on the right hand side of Eq.~\eqref{Eq:SelfSimNonSym} and $F$ denotes the nonlinearity. 
In the following, we study Eq.~\eqref{Eq:YM_Abstract} on a Hilbert space $(\mc H, \| \cdot \|)$ defined as the
completion of the set of radial, compactly supported functions on $\R^7$ with respect to the norm
\begin{align}\label{Def:Hnorm}
\| u \|^2 = \| \Lap u \|_{L^2(\R^7)}^2 + \| \Lap^2 u \|_{L^2(\R^7)}^2,
\end{align}
see Section \ref{Sec:FuncSpace} for details.
In this setting, $L_0$ has a realization as an unbounded, closed operator 
on a suitable domain $\mc D(L_0) \subset \mc H$.
In order to study small perturbations of $\mb W$, we insert the ansatz
\[\Psi(\tau) = \mb W + \Phi(\tau) ,\]
which yields 
\begin{align}\label{Eq:YM_AbstractPertubation}
\begin{split}
\frac{d}{d \tau} \Phi(\tau)  & = (L_0 + L') \Phi(\tau)  + N(\Phi(\tau)) ,\quad \tau > 0 \\
\Phi(0) & =   U(v_0, T).
\end{split}
\end{align}
Here, $L'u = V(|\cdot|)u$ is a linear perturbation with
\begin{align}\label{Eq:Potential5dim} 
  V(\rho) = - 18 \mb W(\rho) - 9\rho^2 \mb W(\rho)^2= \frac{72(36-14 \sqrt{6} +  (\sqrt{6} -2) \rho^2 )}{(36-14 \sqrt{6} +  \sqrt{6} \rho^2)^2}
\end{align}
for $\rho \in [0,\infty)$ and $N$ denotes the remaining nonlinearity.
A short calculation shows that 
\[ N(\Phi(\tau)) = - 9 \big [1+ |\cdot|^2   \mb W(|\cdot|)\big ] \Phi(\tau)^2 -  3 |\cdot|^2  \Phi(\tau)^3. \]
Furthermore, 
\[  U(v_0, T)  :=  T v_0 ( \sqrt{T} |\cdot| ) + \tfrac{T}{T_0} \mb W \left (\tfrac{ \sqrt{T} }{\sqrt{T_0}} |\cdot| \right ) - \mb W(|\cdot|) \]
denotes the transformed initial condition. 
In Section $3$ we show that the operator $L := L_0 + L'$, equipped with a suitable domain, generates 
a strongly continuous one-parameter semigroup $\{S(\tau) : \tau \geq 0\}$ 
of bounded linear operators on $\mc H$. 
For general initial conditions in $\mc H$ we do not expect to obtain classical solutions
to Eq.~\eqref{Eq:YM_AbstractPertubation}. Thus, 
we look for mild solutions $\Phi \in C([0,\infty), \mc H)$ that satisfy the integral equation
\begin{align}\label{Eq:Duhamel}
\Phi(\tau) = S(\tau)  U(v_0, T)   + \int_0^{\tau} S(\tau - \tau') N ( \Phi(\tau')) d \tau'
\end{align}
for $\tau \geq 0$. With these preliminaries we can state the following theorem.

\begin{theorem}\label{Th:MainSim}
Fix $T_0 > 0$. Let $M > 0$ be sufficiently large and $\delta > 0$ sufficiently small. 
For every $v_0 \in \mc H$ with 
\[ \| v_0 \| \leq \frac{\delta}{M^2} \] 
there exists  a $T = T_{v_0} \in [T_0 - \frac{\delta}{M}, T_0 + \frac{\delta}{M}]$ and a function $\Phi  \in C([0,\infty), \mc H)$ that satisfies
Eq.~\eqref{Eq:Duhamel} for all $\tau \geq 0$. Furthermore, 
\[ \| \Phi (\tau) \| \lesssim e^{-\frac{1}{150} \tau} \]
for all $\tau \geq 0$. 
\end{theorem}

Sections \ref{Sec:FuncSpace} - \ref{Sec:NonlinearStability} are mainly  devoted to the proof of Theorem \ref{Th:MainSim}.
In Section \ref{Subsec:FinalArg} we show that Theorem \ref{Th:MainSim} implies Theorem \ref{Th:Main}.

\subsection{Outline of the proof}\label{Sec:OutlineProof}

\subsubsection{Functional analytic setup}
The choice of the topology is basically determined by two requirements: First, the function space must allow
for a rigorous analysis of the linearized problem in terms of strongly continuous semigroups, 
including suitable growth bounds for the linearized time evolution. Furthermore, we need local Lipschitz estimates 
for the nonlinearity in order to run a fixed-point argument. 
We note that the operator $L_0$ corresponds to an Ornstein-Uhlenbeck operator on $\R^7$, cf. \cite{Met01} .
The natural topology for the 'free' problem would therefore be the weighted space $L_{\mu}^2(\R^7)$ with
$\mu(\xi) = e^{-|\xi|^2/4}$ since in this setting, $L_0$ has a self-adjoint realization.
However, the fact that the weight function decays at infinity renders a nonlinear perturbation theory 
hopeless. To circumvent this problem, we omit the weight function altogether and study the problem on
$L^2$-based Sobolev spaces. On $L^2(\R^7)$, the operator $L_0$ is truly non-selfadjoint, however, 
it still has a realization as the generator of a $C_0$-semigroup, cf.~\cite{Met01}.
As described in \cite{Met01}, the spectrum of the generator on $L^2$ is given by a half-plane of
spectral points and the spectral bound is positive. This causes the corresponding semigroup
to grow exponentially. In this paper, we work instead with homogeneous Sobolev norms 
and introduce a Hilbert space $(\mc H, \| \cdot\|)$ consisting of radial functions, with norm given by 
\[ \| u \|^2 = \| \Lap u \|_{L^2(\R^7)}^2 + \| \Lap^2 u \|_{L^2(\R^7)}^2, \]
see Section \ref{Sec:FuncSpace}. We show that if $u \in \mc H$, then 
$u \in C_{\mathrm{rad}}(\R^7) \cap  C_{\mathrm{rad}}^3(\R^7\setminus\{0\})$.
%and $u \in H_{\mathrm{rad}}^4(\B^7_R)$ for any $R > 0$.  
Furthermore, elements of $\mc H$ decay at infinity according to 
\[ \lim_{|\xi| \to \infty} |\xi|^{\frac32} |u(\xi)| = 0. \]

\subsubsection{Semigroup theory and spectral analysis}
We prove that $L_0$ equipped with
a suitable domain $\mc D(L_0)$ generates a $C_0$-semigroup $\{S_0(\tau) : \tau \geq 0 \}$ on $\mc H$ that satisfies
 \[ \|S_0(\tau) u \| \leq e^{- \frac{1}{4} \tau} \|u \|, \quad \forall \tau \geq 0.\]
Furthermore, we show that $ L_0 u(\xi) =  \ell_0 u(|\xi|)$, where 
\[ \ell_0 u(\rho) = \Lap_{\mathrm{rad}} u(\rho) - \tfrac{1}{2} \rho u'(\rho) - u(\rho), \quad \rho > 0, \]
in a classical sense. 
By exploiting the decay of the potential $V$ at infinity, we show that the perturbation $L'$ is bounded and
compact relative to $(L_0,\mc D(L_0))$. Standard semigroup theory implies that 
$L = L_0 + L'$, $\mc D(L) = \mc D(L_0)$, generates a $C_0$-semigroup $\{S(\tau) : \tau \geq 0 \}$. However, the abstract
theory only yields
\[ \|S(\tau) u\| \leq e^{( \| L' \| - \frac{1}{4} ) \tau} \|u  \| \]
for $\tau\geq 0$.
To improve this bound, we analyze the spectrum of the generator. Standard results from semigroup theory 
and the compactness of the perturbation relative to $L_0$ imply that spectral points of $L$ with $\mathrm{Re} \lambda > -\frac14$
are eigenvalues. Hence, for such spectral points it suffices to study the 
eigenvalue problem which reduces to a second order ODE with a regular singular point 
at zero and an essential singularity at infinity. For $\lambda = 1$, the eigenvalue equation can be solved 
explicitly which reveals the existence of an unstable eigenvalue. However, this instability is caused by the
time translation symmetry of the problem (therefore, the corresponding eigenfunction $\mb g$ is referred to as the \textit{symmetry mode})
and can be controlled in the further analysis. To exclude other unstable eigenvalues, we apply a well-known trick from 
supersymmetric quantum mechanics, cf.~ for example \cite{CostinDonn}, which allows us to prove that 
\[ \sigma(L) \subseteq \{ \lambda \in \C: \mathrm{Re} \lambda  \leq - \tfrac{1}{75} \} \cup \{ 1 \}, \]
see Section \ref{Sec:Superpartner} and Section \ref{Sec:SpectrumL}. Since the symmetry eigenvalue is isolated,
we can define a spectral projection $P$ that commutes with the semigroup such that $\sigma(L|_{\mathrm{rg P}}) = \{1\}$. 
Furthermore, we show that $\rg P = \ker (1- L) = \mathrm{span}(\mb g)$.

\subsubsection{Resolvent estimates}
The aim is to translate the spectral bound for $L$ on the stable subspace $\ker P \subset \mc H$ into a growth bound for
the corresponding restricted semigroup. It is well-known that the growth bound $\omega_S$ of the semigroup relates to the spectral radius $r(S(\tau))$
of the operator $S(\tau)$ via $r(S(\tau)) = e^{\omega_S \tau}$. Moreover,  
all spectral points of the generator are mapped to spectral points of $S(\tau)$, i.e.,
\[ \{ e^{\lambda \tau }: \lambda \in \sigma(L) \} \subset \sigma(S(\tau))\setminus \{0\} . \]
If the converse inclusion is also true, one speaks of \textit{spectral mapping} and the latter is known to hold for large classes of semigroups with certain 'nice'
properties (such as eventual norm-continuity).

\noindent However, even if spectral mapping seems to be the generic case \cite{Ren96}, 
many counterexamples are available in the
literature \cite{engel}.  A particular striking one comes from PDE theory \cite{Ren94} where spectral mapping is destroyed 
by a lower-order (and relatively compact) perturbation of the wave equation. 
Already on $L^2(\R^d)$, the properties of semigroups generated by Ornstein-Uhlenbeck operators with potential 
do not allow for the application of simple abstract results that would imply the desired relation between spectral bound and growth bound,
see \cite{DonSch14a}. It cannot be expected that the situation improves much for the semigroup generated by the perturbed OU-operator $(L,\mc D(L))$ on $\mc H$. 

\noindent Fortunately, in Hilbert spaces, one can resort to the classical theorem
by Gearhardt, Pr\"uss, Huang and Greiner which guarantees spectral mapping provided the resolvent of the generator
is uniformly bounded. By proving such estimates, the authors have established
 spectral mapping for a large class of (not necessarily radial) OU-operators with potential on $L^2(\R^d)$, see 
\cite{DonSch14a}.
In Section \ref{Sec:Resolventest}, we apply this machinery to prove the required bounds for the resolvent $R_L(\lambda)$ on $\mc H$. 
In fact, the construction of the resolvent reduces to solving the ordinary differential equation 
\begin{equation*}
\lambda u(\rho) - u''(\rho) - \tfrac{6}{\rho} u'(\rho) + \tfrac{1}{2} \rho u'(\rho) + u(\rho) - V(\rho) u(\rho)  = f(\rho),
\end{equation*}
for $f(\cdot) = f(|\cdot|) \in \mc H$ and $V$ given by Eq.~\eqref{Eq:Potential5dim}. For $\lambda = \alpha + \I \omega$ 
and $\alpha > - \frac{1}{75}$ fixed, we consider the 
problem in normal form and perturbatively construct fundamental systems for the homogeneous equation that are valid in different 
regimes (for small/large/intermediate values of $\rho$), provided $\omega \gg 1$. Furthermore, we show that the corresponding
functions behave well under differentiation and keep track of the dependence
of all quantities on $\omega$, which is crucial to prove uniform bounds later on. A global fundamental system is obtained
by 'gluing' together 
fundamental systems in  different regimes. With the variation of constants
formula, we get an explicit representation of the resolvent as an integral operator on $\mc H$ that is valid for
large imaginary parts. By a detailed analysis of weighted derivatives of the resolvent
up to fourth order we are able to prove that 
\[ \| R_L(\alpha + \I \omega) \| \leq C_{\alpha} \]
for all $|\omega| \gg 1$. An application of the Pr\"uss theorem \cite{Pru84} implies that
\[ \|P S(\tau) u \| = e^{\tau} \|P u \|, \quad \|(1-P)S(\tau) u \| \leq C_{\varepsilon} e^{-(\frac{1}{75} - \varepsilon) \tau} \|(1-P)u \|.  \]
For convenience we set $\varepsilon = \frac{1}{150}$. We note that in the corresponding analysis for wave equations, the
problem of spectral mapping is trivial since there, $L'$ defines a compact perturbation and standard results apply.

\subsubsection{Nonlinear perturbation analysis}
For the nonlinear problem, we first show that the nonlinearity $N$ is locally Lipschitz 
and continuously Fr\'echet-differentiable.
Furthermore, we define the initial data operator $U( v, T)$ on $\mc H$ and show some continuity
properties. In particular, we convince ourselves that $U( v , T)$ is small provided $v$ is small and $T$ is close to $T_0$. 
For the investigation of Eq.~\eqref{Eq:Duhamel} we then proceed as in our previous works. We study the problem
on a function space
\[ \mc X := \{ \Phi \in C([0,\infty), \mc H) : \| \Phi \|_{\mc X} := \sup_{\tau \geq 0}  e^{\frac{1}{150} \tau } \|\Phi(\tau) \|   < \infty \}. \]
To suppress the exponential growth of the semigroup on the unstable subspace, we define a correction term
\[C(\Phi, U(v, T)) :=  P U(v, T) + \int_0^{\infty} e^{-\tau'} P N(\Phi(\tau')) d\tau' \]
which is added to Eq.~\eqref{Eq:Duhamel} (this is reminiscent of center manifold theory).
The modified problem reads 
\begin{align}\label{Eq:DuhamelModified_Outline}
\Phi(\tau) = S(\tau)  U(v, T)   + \int_0^{\tau} S(\tau - \tau') N ( \Phi(\tau')) d \tau' - e^{\tau}  C(\Phi, U(v, T) ).
\end{align}
By applying the Banach fixed point theorem, we show that for every small $v$ and every $T$ close to $T_0$,
a unique solution $\Phi \in C([0,\infty), \mc H)$ to Eq.~\eqref{Eq:DuhamelModified_Outline} exists and the latter satisfies
$\| \Phi(\tau) \| \lesssim e^{-\frac{1}{150} \tau}$. In a final step, we show that for every small $v \in \mc H$
there exists a particular $T = T_{v}$ such that $C(\Phi, U(v, T_{v}) ) = 0$. This crucially relies on the fact that 
the correction $C$ has values in $\mathrm{span(\mb g)}$. Consequently, the equation $C(\Phi, U(v, T) ) = 0$ can be reformulated  
as a one-dimensional fixed point problem. This proves Theorem \ref{Th:MainSim}.

Under the assumptions of Theorem \ref{Th:Main} we can apply Theorem \ref{Th:MainSim} to
obtain a unique solution $\Phi \in  C([0,\infty), \mc H)$
to the equation 
\begin{align*}
\Phi(\tau) = S(\tau)\Phi^T_0  + \int_0^{\tau} S(\tau - \tau') N ( \Phi(\tau')) d \tau', \quad \tau \geq 0,
\end{align*}
for $\Phi^T_0 = T u_0(T^{\frac12} |\cdot|) - \mb W(|\cdot|)$. The conditions on the 
initial data in Theorem \ref{Th:MainSim} imply that 
 $\Phi^T_0 \in \mc D(L)$ (obviously, this is true for a much larger class
of initial conditions).
Since the nonlinearity is continuously differentiable, 
standard semigroup arguments imply that $\Phi \in  C([0,\infty), \mc H) \cap C^1((0,\infty),\mc H)$ and
$\Phi(\tau) \in \mc D(L)$ for all $\tau > 0$. In particular, $\Phi(\tau)(\cdot) = \varphi(\tau, |\cdot|)$
and $\varphi$ satisfies
\begin{align}\label{Eq:SelfSimPert}
\begin{split}
 \partial_{\tau} \varphi  =  &   \partial^2_{\rho} \varphi  + \tfrac{6}{\rho}  \partial_{\rho} \varphi
- \tfrac{1}{2} \rho  \partial_{\rho} \varphi   -  \varphi  + V \varphi + N(\varphi)
\end{split}
\end{align}
in a classical sense for all $\tau > 0$ with $\varphi(0,\rho) = T u_0(0,\sqrt{T} \rho) - \mb W(\rho)$. 
By setting 
\[ u(t, r) := (T-t)^{-1} \varphi(- \log(T-t) + \log T , \tfrac{r}{\sqrt{T-t}}) +   (T-t)^{-1} \mb W(\tfrac{r}{\sqrt{T-t}}) \]
we obtain a classical solution to Eq.~\eqref{Eq:EquivarEq5} for all $t > 0$ with $u(0,\cdot) = u_0$. The exponential decay rate of $\Phi$
translates into the convergence rates for $u$ given in Theorem $\ref{Th:Main}$.

\subsection{Notation and Conventions}
We write $\N$ for the natural numbers $\{1,2,3, \dots\}$ and set $\N_0 := \{0\} \cup \N$. Furthermore, $\R^+ := \{x \in \R: x >0\}$.
The notation $a\lesssim b$ means $a\leq Cb$ for an absolute constant $C>0$ and we also write $a\simeq b$ if $a\lesssim b$ and $b \lesssim a$.  	
If $a \leq C_{\varepsilon} b$ for a constant $C_{\varepsilon}>0$ depending on some parameter $\varepsilon$, we write $a \lesssim_{\varepsilon} b$. 
We use the common notation $\langle x \rangle := \sqrt{1+|x|^2}$ also known as the \textit{Japanese bracket}.

For a function $x \mapsto g(x)$, we denote by $g^{(n)}(x) = \frac{d^n g(x)}{dx^n} $ the derivatives of order $n \in \N$. 
For $n=1,2$, we also write $g'(x)$ and $g''(x)$, respectively. For a function $(x,y) \mapsto f(x,y)$, partial derivatives of order $n$ will be denoted by 
$\partial^{n}_x f (x,y) = \frac{\partial^n}{\partial x^n} f(x,y)$. Throughout the paper,
$W(f,g)$ denotes the Wronskian of two functions $x \mapsto f(x)$ and $x \mapsto g(x)$, where we use the convention
$W(f,g)=fg'-f'g$.

By $C^{\infty}_0(\R^d)$ we denote the set of compactly supported functions on $\R^d$, $d \geq 1$. 
The spaces $L^2(\R^d)$ and $H^{k}(\R^d)$, $k \in \N_0$, denote the standard Lebesgue  
and Sobolev spaces with the usual norm
\[ \|  u \|^2_{H^k(\R^d)} := \sum_{\alpha: |\alpha| \leq k} \|\partial^{\alpha}  u \|^2_{L^2(\R^d)}.\]

The set of bounded linear operators on a Hilbert space $\mc H$ is denoted by $\mc B(\mc H)$. 
For a closed linear operator $(L, \mc D(L))$, we write $\sigma(L)$ for the spectrum.
The resolvent set is defined as $\rho(L) := \C \setminus \sigma(L)$ and we write 
$R_{L}(\lambda):=(\lambda- L)^{-1}$ for $\lambda \in \rho(L)$.

\section{Function spaces}\label{Sec:FuncSpace}

\subsection{Preliminaries} 
On $\R^7$, we introduce the set of radial test functions,
\begin{align*}
C^{\infty}_{\mathrm{rad},0}(\R^7) := \{ \tilde u \in C_{0}^{\infty}(\R^7): \tilde u \text{ is radial} \},
\end{align*}
and define $C^{\infty}_{\mathrm{e},0}(\R)$, the set of even test functions on $\R$, as in Eq.~\eqref{Def:EvenTestFunc}. 
Note that if $u \in C^{\infty}_{\mathrm{e},0}(\R)$, then $u^{(2k+1)}(0)=0$ for all $k \in \N_0$. 
Furthermore, every $u \in C^{\infty}_{\mathrm{e},0}(\R)$ defines a function $\tilde u \in  C^{\infty}_{\mathrm{rad},0}(\R^7)$ 
by $\tilde u(\xi) := u(|\xi|)$. Conversely, if $\tilde u \in C^{\infty}_{\mathrm{rad},0}(\R^7)$, 
then $\tilde u(\xi) = u(|\xi|)$ for some $u \in C^{\infty}_{\mathrm{e},0}(\R)$.

\medskip
In the following, we set 
\[ L^2_{\mathrm{rad}}(\R^7):= \{ \tilde u \in L^2(\R^7):  \tilde u \text{ is radial} \}.\]
For $\tilde u \in L^2_{\mathrm{rad}}(\R^7)$, we have $\tilde u(\xi) = u(|\xi|)$ a.e. on $\R^7$, for some $u \in L^2(\R)$,
satisfying $u(x) = u(-x)$ a.e on $\R$. By using polar coordinates $(\rho, \omega)$, $\rho = |\xi|, \omega = \frac{\xi}{|\xi|}$, the
inner product on $ L^2_{\mathrm{rad}}(\R^7)$ can be written as  
\[ (\tilde u|\tilde v)_{L^2(\R^7)} = \int_{\R^7} u(|\xi|) \overline{v(|\xi|)} d\xi = C \int_0^{\infty} u(\rho) \overline{v(\rho)}\rho^6 d\rho, \]
where the constant comes from the integration over $\mathbb{S}^6$. 
Furthermore, we consider Sobolev spaces 
\[ H^k_{\mathrm{rad}}(\R^7):= \{ \tilde u \in H^k(\R^7):  \tilde u \text{ is radial} \}, \]
for $k \in \N$, where the norm on $H^k$ is defined in the usual manner. For $\tilde u \in H_{\mathrm{rad}}^{2j}(\R^7)$, 
$j=1,2$, the operators $\Lap^j$ exist in the weak sense, i.e., there are functions  $\tilde f_j \in L^2_{\mathrm{rad}}(\R^7)$ such that 
\begin{align}\label{Eq:TestFunc}
 \int_{\R^7} \tilde u(\xi) \Lap^j \tilde \phi(\xi) d\xi =   \int_{\R^7} \tilde f_j(\xi) \tilde \phi(\xi) d\xi 
\end{align}
for all $\tilde \phi \in C^{\infty}_{\mathrm{rad},0}(\R^7)$ (this is not a restriction since for non-radial testfunctions, 
Eq.~\eqref{Eq:TestFunc} can be recast by introducing polar coordinates and writing $\tilde \phi$ as a spherical mean).
%By change of coordinates this is equivalent to the existence of even functions $f_j \in L^2(\R)$ such that 
%\begin{align*}
% \int_{0}^{\infty} u(\rho) \Lap_{\mathrm{rad}}^j \phi(\rho) \rho^6 d\rho =    \int_{0}^{\infty} f_j(\rho) \phi(\rho) \rho^6 d\rho,
%$\end{align*}
%for all $\phi \in C^{\infty}_{\mathrm{e},0}(\R)$, where 

In the following, we use the notation
\[ \Lap_{\mathrm{rad}} u(\rho) = u''(\rho) + \tfrac{6}{\rho} u'(\rho),\]
for the radial Laplace operator. If $\tilde u = u(|\cdot|) \in C^{\infty}_{\mathrm{rad},0}(\R^7)$, then $\nabla \Lap \tilde u \in C^{\infty}_{0}(\R^7)$ and 
\[ |\nabla \Lap \tilde u(\xi)|^2 =  |(\Lap_{\mathrm{rad}} u)'(\rho)|^2. \]

\subsection{Definition of the Hilbert space $\mc H$}

On $C^{\infty}_{\mathrm{rad},0}(\R^7)$ we introduce the inner product
\begin{align*}
(\tilde u|\tilde v) := (\Lap \tilde u|\Lap \tilde v)_{L^2(\R^7)} + (\Lap^2 \tilde u|\Lap^2 \tilde v)_{L^2(\R^7)},
\end{align*}
%where 
%\[ (\Lap \tilde u|\Lap \tilde v)_{L^2(\R^7)} = C_{\Omega} \int_0^{\infty} \Lap_{\mathrm{rad}} u(\rho) \overline{\Lap_{\mathrm{rad}} v(\rho)} \rho^6 d\rho^, \]
%and 
%\[ (\Lap^2 \tilde u|\Lap^2 \tilde v)_{L^2(\R^7)} = C_{\Omega}  \int_0^{\infty} \Lap^2_{\mathrm{rad}} u(\rho) \overline{\Lap^2_{\mathrm{rad}} v(\rho)} \rho^6 d\rho. \]
and a norm $\|\tilde u \| := \sqrt{(\tilde u|\tilde u)}$, i.e.,   
\[ \| \tilde u \|^2 =  \| \Lap \tilde u \|_{L^2(\R^7)}^2 + \| \Lap^2 \tilde u\|_{L^2(\R^7)}^2. \]
Next, we state some bounds that will be very useful in the following. 

\begin{lemma}\label{Le:AllBounds}
Let $\tilde u = u(|\cdot|) \in C^{4}_{\mathrm{rad}}(\R^7)$, such that $\|\tilde u \| < \infty$ and 
\begin{align}\label{Eq:CondLeAllB1}
 \lim_{\rho \to \infty} \rho^{3} |\Lap_{\mathrm{rad}}u(\rho)| = 0, \quad  \lim_{\rho \to \infty} \rho^{3} 
|(\Lap_{\mathrm{rad}}u)'(\rho)| = 0. 
\end{align}
Then we have the bound
\[ \|\nabla \Lap  \tilde u\|_{L^2(\R^7)}  \lesssim \|\tilde u \|.\]

If in addition
\[ \lim_{\rho \to \infty} \rho^{\frac32} |u(\rho)| = 0, \quad  \lim_{\rho \to \infty} \rho^{\frac52} |u'(\rho)| = 0, \]
then 
\begin{align}\label{Eq:L2Bounds}
\| (\cdot)^{\alpha} u \|_{L^2(\R^+)}   & \lesssim \|\tilde u \| \quad \text{ for }  \alpha \in [0,1], & 
\| (\cdot)^{\alpha} u' \|_{L^2(\R^+)}    & \lesssim \|\tilde u \| \quad  \text{ for }  \alpha \in [0,2],  \nonumber  \\
\| (\cdot)^{\alpha} u'' \|_{L^2(\R^+)}   & \lesssim \|\tilde u \|\quad  \text{ for }  \alpha \in [ 1,3],  
&  \| (\cdot)^{\alpha} u^{(3)} \|_{L^2(\R^+)}   & \lesssim \|\tilde u \|  \quad \text{ for }  \alpha \in [2,3],  \\
\| (\cdot)^{3} u^{(4)} \|_{L^2(\R^+)}     & \lesssim \| \tilde u \|,  \nonumber 
\end{align}
and
\begin{align}\label{Eq:LinftyBounds}
\| (\cdot)^{\alpha}  u \|_{L^{\infty}(\R^+)}  & \lesssim \|\tilde u \| \quad \text{ for } \alpha \in [0, \tfrac32], &
\| (\cdot)^{\alpha}  u' \|_{L^{\infty}(\R^+)}    & \lesssim \|\tilde u \| \quad  \text{ for } \alpha \in [1, \tfrac52],  \nonumber  \\
\| (\cdot)^{\alpha}  \Lap_{\mathrm{rad}}  u \|_{L^{\infty}(\R^+)}   & \lesssim \|\tilde u \|\quad  \text{ for }  \alpha \in [2, 3] ,  & 
\| (\cdot)^3  (\Lap_{\mathrm{rad}}  u)' \|_{L^{\infty}(\R^+)}   & \lesssim \|\tilde u \|.
\end{align}
\end{lemma}

\begin{proof}
By scaling, it is natural to expect the bounds
\begin{align*}
 \sum_{j=1}^4 \|(\cdot)^{j-1}u^{(j)}\|_{L^2(\R^+)}\lesssim \|\tilde u\|_{\dot H^4(\R^7)},\qquad 
 \sum_{j=0}^2 \|(\cdot)^{j+1}u^{(j)}\|_{L^2(\R^+)}\lesssim \|\tilde u\|_{\dot H^2(\R^7)} 
\end{align*}
and Hardy's inequality shows that these estimates are indeed correct.
Based on this, the stated assertions follow easily by interpolation and Sobolev embedding.
A self-contained and elementary proof including all details is given in Appendix \ref{Sec:ProofAllBounds}.
\end{proof}

\begin{lemma}\label{Le:PropH}
Let $\mc H$ denote the completion of $(C^{\infty}_{\mathrm{rad},0}(\R^7), \|\cdot\|)$. 
Then $\mc H$ is a Hilbert space and its elements can be identified with 
functions $\tilde u = u(|\cdot|)  \in C_{\mathrm{rad}}(\R^7) \cap C_{\mathrm{rad}}^3(\R^7\setminus\{0\})$, that satisfy 
\[ \lim_{\rho \to \infty} \rho^{\frac32} |u(\rho)| = 0, \quad  \lim_{\rho \to \infty} \rho^{\frac52} |u'(\rho)| = 0. \]
%Furthermore, $\tilde u \in H^4(\B^7_R)$ for any $R > 0$. 
The norm induced by the inner product on $\mc H$ is given by 
\[ \|\tilde u\|^2 = \| \Lap \tilde u\|^2_{L^2(\R^7)} +  \| \Lap^2 \tilde u\|^2_{L^2(\R^7)}, \]
where $\Lap \tilde u$ can be interpreted as a classical differential operator and $\Lap^2 \tilde u$ has 
to be understood in a weak sense, cf. Eq.~\eqref{Eq:TestFunc}.
Finally, we have 
\begin{align}\label{Eq:NormInterpolation}
\| \nabla \Lap \tilde u\|_{L^2(\R^7)}  \lesssim \|\tilde u \|,
\end{align}
as well as 
\begin{align}\label{Eq:LinftyH}
\| u \|_{L^{\infty}(\R^+)} \lesssim \|\tilde u \|,
\end{align}
for all $\tilde u \in \mc H$
\end{lemma}

\begin{proof}
Let $(\mc H, \| \cdot \|_{\mc H})$ denote the completion of $(C^{\infty}_{\mathrm{rad},0}(\R^7), \|\cdot\|)$ and let $U \in \mc H$. 
By construction, $U = [(\tilde u_n)_{n \in \N}]$ is an equivalence class of Cauchy sequences in $C^{\infty}_{\mathrm{rad},0}(\R^7)$
and 
\[\| U \|^2_{\mc H} = \lim_{n \to \infty} \| \tilde u_n \|^2 
= \lim_{n \to \infty}  \|\Lap \tilde u_n \|^2_{L^2(\R^7)} +  \lim_{n \to \infty}  \| \Lap^2 \tilde u_n \|^2_{L^2(\R^7)},\]
for some representative $(\tilde u_n)_{n \in \N}$.
By completeness of $L^2(\R^7)$, $\Lap \tilde u_n \to f_1$, $ \Lap^2 \tilde u_n \to f_2$ in $L_{\mathrm{rad}}^2(\R^7)$. In particular,
we have that 
\[  (\Lap^j \tilde u_n|  \tilde \varphi)_{L^2(\R^7)}  \to (f_j| \tilde \varphi)_{L^2(\R^7)},\]
for $j=1,2$, and all $\tilde \varphi \in C^{\infty}_{\mathrm{rad},0}(\R^7)$.
Our aim is to show that  
the functions $f_j$ can be identified with weak derivatives in the sense that there is a unique 
$\tilde u \in C_{\mathrm{rad}}(\R^7)$ such that for all 
$\tilde \varphi \in C^{\infty}_{\mathrm{rad},0}(\R^7)$,
\begin{align}\label{Eq:WeakLap}
(\tilde u| \Lap^j \tilde \varphi)_{L^2(\R^7)} = (f_j| \tilde \varphi)_{L^2(\R^7)}.
\end{align}

The sequence $(\tilde u_n)_{n \in \N}$,  $\tilde u_n = u_n(|\cdot|)$, is a Cauchy sequence in $\mc H$ and by  
Lemma \ref{Le:AllBounds}, $(u_n)_{n \in \N}$ converges in $H^1(\R)$. We denote the corresponding limit function by $u$.
Convergence of $(u_n)_{n \in \N}$ in $L^{\infty}(\R)$ implies that $u \in C(\R)$ and by pointwise convergence, $u$ is even.
The $L^2$-bounds show that $u \in H^4(\delta,\infty)$ for any $\delta > 0$ and by Sobolev embedding we get that $u \in C^3(0,\infty)$.
Since $u$ is approximated by functions with compact support, the $L^{\infty}$-bounds yield
\[ \lim_{\rho \to \infty} \rho^{\frac32} |u (\rho)| = 0, \quad  \lim_{\rho \to \infty} \rho^{\frac52} |u' (\rho)| = 0. \]

We define $\tilde u := u(|\cdot|)$, such that $\tilde u \in C_{\mathrm{rad}}(\R^7) \cap C_{\mathrm{rad}}^3(\R^7\setminus\{0\})$.
In general $\tilde u \notin  L^2(\R^7)$, however,
$(\cdot)^3 u_n \to (\cdot)^3 u$ in $L^2(0,R)$, which implies that $\tilde u_n \to \tilde u \in L^2(\B^7_R)$ for any $R > 0$.
In particular, 
\[(\tilde u_n| \Lap^j \tilde \varphi)_{L^2(\R^7)} \to (\tilde u| \Lap^j \tilde \varphi)_{L^2(\R^7)},\]
for $j = 1,2$ and all $\tilde \varphi \in C^{\infty}_{\mathrm{rad},0}(\R^7)$.
Since $(\tilde u_n| \Lap^j \tilde \varphi)_{L^2(\R^7)} = (\Lap^j \tilde u_n|  \tilde \varphi)_{L^2(\R^7)}$ by integration by parts, we 
get that 
\[  (\Lap^j \tilde u_n|  \tilde \varphi)_{L^2(\R^7)} \to (\tilde u| \Lap^j \tilde \varphi)_{L^2(\R^7)},\]
and as a consequence Eq.~\eqref{Eq:WeakLap}, i.e., $\Lap^j u $ exists in a weak sense and we can write 
\[\| U \|^2_{\mc H} = \|\Lap \tilde u\|^2_{L^2(\R^7)} + \| \Lap^2 \tilde u \|^2_{L^2(\R^7)}.\]
The above arguments do not depend on the particular choice of the representative  and
the constructed limit function $\tilde u = u(|\cdot|)$ is unique. Hence, we can identify $U$ with $\tilde u$ and use the
notation $\| U \|_{\mc H} = \| \tilde u\| $.

Note that by Lemma \ref{Le:AllBounds}, $(\cdot) u'_n \to  (\cdot) u' \in L^{\infty}(\R^+)$. In particular, $\lim_{\rho \to 0} \rho u'(\rho) = 0$.
Hence, we can integrate by parts to obtain 
\[(\tilde u| \Lap \tilde \varphi)_{L^2(\R^7)} = C \int_0^{\infty} u(\rho) \Lap_{\mathrm{rad}} \varphi(\rho) \rho^6 d\rho 
= C \int_0^{\infty} \Lap_{\mathrm{rad}} u(\rho) \varphi(\rho) \rho^6 dr = (\Lap \tilde u|  \tilde \varphi)_{L^2(\R^7)}. \]
This shows that  $\Lap \tilde u$ can be interpreted as a classical differential operator in the above norm. 
%Corollary \ref{Cor:H4Bound} and fact that $\tilde u_n \to \tilde u$ in $L^2(\B^7_R)$ imply that 
%$\tilde u \in H^4(\B^7_R)$. 
Finally, Eq.~\eqref{Eq:NormInterpolation} and Eq.~\eqref{Eq:LinftyH} follow from Lemma \ref{Le:AllBounds} by approximation.
\end{proof}

%\begin{remark}\label{Re:FinitenormimpliesH}
%Assume that $\tilde u \in C_{\mathrm{rad}}^4(\R^7)$, such that $\|\tilde u \| < \infty$. 
%By using standard techniques, one can always construct a sequence $(\tilde u_n)_{n \in %\N} \in C^{\infty}_{\mathrm{rad},0}(\R^7)$ 
%such that $\lim_{n \to \infty} \| \tilde u_n - \tilde u \| =0,$
%which implies that $\tilde u \in \mc H$. 
%\end{remark}

After these preliminaries we can now turn to the investigation of Eq.~\eqref{Eq:YM_AbstractPertubation}.

\section{Semigroup theory and spectral analysis}\label{Sec:SemigroupTh}

\subsection{The Ornstein-Uhlenbeck operator on $\mc H$} 

We set 
\[ \Lambda \tilde  u(\xi) := \frac{1}{2} \xi \cdot \nabla \tilde u(\xi),\]
for $\xi \in \R^7$ and define the formal differential expression
\[ \mc L_0 \tilde u(\xi) := \Lap \tilde  u(\xi) -  \Lambda \tilde  u(\xi) - \tilde  u(\xi). \]
In polar coordinates, $\mc L_0$ decouples into a radial and an angular part. In particular, for $\tilde u = u(|\cdot|) \in  C_{\mathrm{rad}}^{\infty}(\R^7)$,  
$\mc L_0 \tilde u(\cdot) = \tilde \ell_0 u(|\cdot|)$, where 
\[\tilde \ell_0 u(\rho)  =    \Lap_{\mathrm{rad}} u(\rho) - \tfrac{1}{2} \rho u'(\rho)  - u(\rho).\]
We define $\tilde L_0 \tilde u := \mc L_0 \tilde u$ on the domain  
\begin{align*}
\mc D(\tilde L_0) :=  
\left \{ \tilde u = u(|\cdot|) \in \mc H \cap C_{\mathrm{rad}}^{6}(\R^7): \Lap^3 \tilde u \in L^2(\R^7),
 \mc L_0 \tilde u \in \mc H,  \right. \\
\left.  \lim_{\rho \to \infty} \rho^{3} |\Lap^{j}_{\mathrm{rad}}u(\rho)| = 0,  \lim_{\rho \to \infty} \rho^{3} 
|(\Lap^{j}_{\mathrm{rad}}u)'(\rho)| = 0, \text{ for } j =1,2 \right \}.
\end{align*}
The operator $\tilde L_0$ is densely defined since $C^{\infty}_{\mathrm{rad},0}(\R^7) \subset \mc D(\tilde L_0)$.
Note that functions in $\mc D(\tilde L_0)$ satisfy the assumptions of Lemma \ref{Le:AllBounds} by Lemma \ref{Le:PropH}. 
We need the following result, which is based on \cite{Met01}. 

\begin{lemma}\label{Le:Metafune}
Consider the operator  $(\Lambda, \mc D(\Lambda))$, with 
\[ \mc D( \Lambda) = \{ \tilde u \in L^2(\R^7) :  \Lambda \tilde u \in L^2(\R^7) \}, \]
where $\Lambda \tilde u$ is understood in the sense of distributions. Then,
\begin{align}\label{Eq:EstLambda}
\mathrm{Re}(-  \Lambda \tilde u | \tilde u)_{L^2(\R^7)} \leq \frac{7}{4} \| \tilde u \|^2_{L^2(\R^7)},
\end{align}
for all $\tilde u \in \mc D( \Lambda)$
\end{lemma}

\begin{proof}
As defined above, the operator is closed by Lemma $2.1$ in \cite{Met01}. Furthermore,
$C_0^{\infty}(\R^7)$ is a core by \cite{Met01}, Proposition $2.2$. Integration by parts shows
that Eq.~\eqref{Eq:EstLambda} holds  for all $\tilde u \in C_0^{\infty}(\R^7)$. This implies the claim. 
\end{proof}

\begin{lemma}\label{Le:LuPhi}
For all $\tilde u \in \mc D(\tilde L_0)$, we have  
\begin{align}\label{EQ:LumerPh}
\mathrm {Re} (\tilde L_0 \tilde u | \tilde u ) \leq - \frac{1}{4} \|\tilde u \|^2. 
\end{align}
\end{lemma}

\begin{proof}
First, one can easily check that for all $\tilde u \in C^{6}(\R^7)$,
\begin{align} \label{Eq:Commutator}
\Delta  \Lambda \tilde u  =  \Lambda \Delta \tilde  u + \Delta \tilde u, 
\quad \Delta^2   \Lambda \tilde u = \Lambda \Delta^2 \tilde u + 2 \Delta^2 \tilde u.
\end{align}
 Hence,
\begin{align}\label{Eq:CommL0}
\Lap \tilde {L}_0 \tilde u = \Lap^2 \tilde u -  \Lambda \Lap  \tilde u - 2\Lap \tilde u, 
\quad  \Lap^2 \tilde {L}_0 \tilde u = \Lap^3 \tilde u - \Lambda \Lap^2  \tilde u - 3 \Lap^2 \tilde u. 
\end{align}

Let $\tilde u \in \mc D(\tilde L_0)$, then $\Lap^j \tilde u \in L_{\mathrm{rad}}^2(\R^7)$, $j=1,2,3$. 
In view of Eq.~\eqref{Eq:CommL0} and the fact that $\tilde L_0 u \in \mc H$,
we get that  
$\Lambda \Lap \tilde u,  \Lambda \Lap^2 \tilde u \in L_{\mathrm{rad}}^2(\R^7)$. 
In particular, $\Lap \tilde u, \Lap^2 \tilde u \in  \mc D(\Lambda)$ and 
\[ \mathrm{Re}(- \Lambda \Lap \tilde u |\Lap \tilde u)_{L^2(\R^7)} \leq \tfrac{7}{4} \| \Lap \tilde u \|^2_{L^2(\R^7)}, \quad 
 \mathrm{Re}(- \Lambda \Lap^2 \tilde u |\Lap^2 \tilde u)_{L^2(\R^7)} \leq \tfrac{7}{4} \| \Lap^2 \tilde u \|^2_{L^2(\R^7)},
\]
by Lemma \ref{Le:Metafune}. 
Integration by parts implies that 
\[ \mathrm{Re} (\Lap^2 \tilde  u | \Lap \tilde  u)_{L^2(\R^7)} \leq 0, \quad   \mathrm{Re} (\Lap^3 \tilde  u | \Lap^2 \tilde  u)_{L^2(\R^7)} \leq 0. \] 
In view of Eq.~\eqref{Eq:CommL0} we infer that Eq.~\eqref{EQ:LumerPh} holds.
\end{proof}

\begin{lemma}\label{Le:DenseRange}
The range of $(\mu - \tilde L_0)$ is dense in $\mc H$ for $\mu = \frac{5}{2}$.
\end{lemma}

\begin{proof}
Let $\tilde f \in C^{\infty}_{\mathrm{rad},0}(\R^7)$ and let $\hat f =\mc F(\tilde f)$ denote its Fourier transform.
The properties of $\tilde f$ imply that $\hat f$ is a radial Schwartz function, i.e., $\hat f(\zeta) = \hat g(|\zeta|)$ for
some even function $g \in \mc S(\R)$. For $\eta = |\zeta|$, we set
\begin{align*}
\hat h(\eta) =  \int_{\eta}^{\infty} \frac{2 e^{-(s^2 - \eta^2)}}{s} \hat g(s) ds
\end{align*}
and define $\hat u(\cdot) := \hat h(|\cdot|)$. It is easy to check that $\hat u$ satisfies the equation 
\begin{align}\label{Eq:DenseRange}
|\zeta|^2 \hat u(\zeta)  - \tfrac{1}{2} \zeta \cdot \nabla \hat u(\zeta)  = \hat f(\zeta),
\end{align}
which reduces to 
\[ \eta^2 \hat h(\eta) -  \tfrac{1}{2} \eta \hat h'(\eta) = \hat g(\eta) \]
in polar coordinates. We consider the integral
\begin{align*}
\eta^{3} \hat h(\eta) =  \int_{0}^{\infty} K(\eta,s) \hat g(s) s^3 ds
\end{align*}
with $ K(\eta,s) :=   2  s^{-4}  \eta^{3}  e^{- (s^2  - \eta^2)} 1_{[0,\infty)}(s-\eta) $. It is easy to check that   
\[|K(\eta,s) | \lesssim \mathrm{min} \{ \eta^{-1}, s^{-1} \},\]
for all $\eta, s \in [0,\infty)$. By \cite{DonKri13}, Lemma 5.5, 
the kernel induces a bounded integral operator on $L^2(\R^+)$. Hence, 
\[\| \hat u \|_{L^2(\R^7)} \simeq \|(\cdot)^3 \hat h \|_{L^2(\R^+)}  \lesssim \|(\cdot)^3 \hat g \|_{L^2(\R^+)} \simeq \| \hat f \|_{L^2(\R^7)}.\]
In a similar manner one can show that
 \[ \| \langle \cdot \rangle ^k \hat u  \|_{L^2(\R^7)} \lesssim_k \|  \langle \cdot \rangle ^k \hat f  \|_{L^2(\R^7)}\]
for all $k \in \N_0$. The right hand side is finite since $\hat f \in \mc S(\R^7)$.
We define $\tilde u := \mc F^{-1}(\hat u)$.
Then  $\tilde u \in  H^k(\R^7)$ for all $k \in \N_0$ by Plancherel's theorem.
In particular, $\tilde u$ can be approximated with respect to $\|\cdot\|$ by functions in $C^\infty_{0,\mathrm{rad}}(\R^7)$ and thus, $\tilde u \in C^{6}_{\mathrm{rad}}(\R^7) \cap \mc H$.
Applying $\mc F^{-1}$ to Eq.~\eqref{Eq:DenseRange} shows that $\tilde u$ is a solution to 
the equation
\[\tfrac{7}{2} \tilde u(\xi) -  \Lap \tilde u(\xi)  + \tfrac{1}{2} \xi \cdot \nabla \tilde u(\xi) = \tilde f(\xi), \]
and we infer that $ (\tfrac{5}{2} - \tilde L_0) \tilde u  =   \tilde f$.
Note that for all $\tilde v = v(|\cdot|) \in C^{\infty}_{\mathrm{rad},0}(\R^7)$, all $\rho \geq 1$ and $j =1,2$,
\begin{align*}
\rho^6 |(\Lap^j_{\mathrm{rad}} v)'(\rho)|^2  & \lesssim \int_{\rho}^{\infty} s^6 
|(\Lap^j_{\mathrm{rad}} v)'(s)|^2 ds + \int_{\rho}^{\infty} s^6 
|(\Lap^j_{\mathrm{rad}} v)''(s)|^2 ds  \\
& \lesssim \int_{\rho}^{\infty} s^6 
|(\Lap^j_{\mathrm{rad}} v)'(s)|^2 ds + \int_{\rho}^{\infty} s^6 
|\Lap^{j+1}_{\mathrm{rad}} v(s)|^2 ds  \\
& \lesssim 
\| \nabla \Lap^j \tilde v \|^2_{L^{\infty}(\R^7)} +  \| \Lap^{j+1} \tilde v \|^2_{L^2(\R^7)}.
\end{align*}
This implies that for $j=1,2$,
\[  \| (\cdot)^3 (\Lap^j_{\mathrm{rad}} v)' \|_{L^{\infty}(1,\infty)}  \lesssim \| \tilde v \|_{H^6(\R^7)}.\]
Similarly, we obtain 
\[  \| (\cdot)^3 \Lap^j_{\mathrm{rad}} v\|_{L^{\infty}(1,\infty)}  \lesssim \| \tilde v \|_{H^6(\R^7)}.\]
By density of $C^{\infty}_{0}(\R^7)$ in $H^6(\R^7)$, we can always find a radial sequence 
$(\tilde u_n)_{n \in \N}$, $\tilde u_n = u_n(|\cdot|) \in C^{\infty}_{\mathrm{rad},0}(\R^7)$ 
that approximates $\tilde u$ in $ H^6(\R^7)$. 
The above calculations show  that 
\[  (\cdot)^3 (\Lap^j_{\mathrm{rad}} u_n)'  \to  (\cdot)^3 (\Lap^j_{\mathrm{rad}} u)' \]
in $L^{\infty}(1,\infty)$, which implies that 
\[ \lim_{\rho \to \infty}  \rho^3 |(\Lap^j_{\mathrm{rad}} u)'(\rho)| = 0 \]
for $j=1,2$.  Similarly, 
\[ \lim_{\rho \to \infty} \rho^3 |\Lap^j_{\mathrm{rad}} u(\rho)| = 0. \]
Hence, $\tilde u \in \mc D(\tilde L_0)$.   
The density of $C^{\infty}_{\mathrm{rad},0}(\R^7)$ in $\mc H$ finally implies the claim. 
\end{proof}

\begin{lemma}\label{Cor:SemigoupS0}
The operator $(\tilde L_0, \mc D( \tilde L_0))$ is closable and the closure $( L_0, \mc D(L_0))$ generates
a strongly continuous one-parameter semigroup $\{ S_0(\tau) : \tau \geq 0\}$ of bounded operators on $\mc H$. The semigroup satisfies 
\[ \|S_0(\tau) \tilde u \| \leq e^{- \frac{1}{4} \tau} \|\tilde u \|\]
for all $\tilde u \in \mc H$ and $\tau \geq 0$.  Furthermore,  $L_0 \tilde u(\xi) =  \ell_0 u(|\xi|)$, where 
\[ \ell_0 u(\rho) = \Lap_{\mathrm{rad}} u(\rho) - \tfrac{1}{2} \rho u'(\rho) - u(\rho), \quad \rho > 0, \]
in a classical sense and $\lim_{\rho \to 0} \ell_0 u(\rho)$ exists.
\end{lemma}

\begin{proof}
The first part of the statement follows from Lemma \ref{Le:LuPhi}, Lemma \ref{Le:DenseRange} and an application of the 
Lumer-Phillips Theorem \cite{engel}, p.~83. The closure is constructed in the usual way: Let 
$(\tilde u_n)_{n \in \N} \subset \mc D(\tilde L_0)$, $\tilde u_n(\cdot) = u_n(|\cdot|)$, be such that $\tilde u_n \to \tilde u$ and 
$\tilde L_0 \tilde u_n \to \tilde f$ in $\mc H$. Then we say $\tilde u \in \mc D(L_0)$ and define $L_0 \tilde u := \tilde f$. Next, we describe $L_0 \tilde u$ in more detail. Convergence in $\mc H$ implies that
\begin{align}\label{Eq:ClosureEq1}
\Lap \tilde u_n \to \Lap \tilde u \in L^2_{\mathrm{rad}}(\R^7)
\end{align}
and $\tilde u  \in C_{\mathrm{rad}}(\R^7) \cap C^3_{\mathrm{rad}}(\R^7\setminus\{0\})$. 
Furthermore, $\tilde u(\cdot) = u(|\cdot|)$ and we have $u_n \to  u$ in $L^{\infty}(\R^+)$ (Lemma \ref{Le:AllBounds}).
Analogously, $L_0 \tilde u \in  C_{\mathrm{rad}}(\R^7)$, $L_0 u(\cdot) = \ell_0 u(|\cdot|)$  and 
\[ \Lap_{\mathrm{rad}} u_n - \tfrac{1}{2} (\cdot) u'_n -   u_n \to \ell_0 u \text{ in } L^{\infty}(\R^+).\]
By Lemma \ref{Le:AllBounds}, $(\cdot) u'_n \to (\cdot) u' $ in $L^{\infty}(\R^+)$.
This and Eq.~\eqref{Eq:ClosureEq1}
imply that $\Lap_{\mathrm{rad}} u_n \to \Lap_{\mathrm{rad}} u$ in $L^{\infty}(\R^+)$. By uniqueness of the limit function 
we get that 
\[ \ell_0 u(\rho) = \Lap_{\mathrm{rad}} u(\rho) - \tfrac{1}{2} \rho u'(\rho) -   u(\rho), \quad \rho > 0. \]
Finally, $\lim_{\rho \to 0} \ell_0 u(\rho)$ exists by continuity of $L_0 \tilde u$ at the origin.
\end{proof}

\subsection{The perturbed problem} 

Next, we convince ourselves that the perturbation $L'$ in Eq.~\eqref{Eq:YM_AbstractPertubation} defines a bounded operator on $\mc H$. 
We actually prove a much stronger statement that will be crucial later on.

\begin{lemma}\label{Le:Compact1}
Let $V$ be defined as in Eq.~\eqref{Eq:Potential5dim}. Then 
\[L'\tilde u :=  V(|\cdot|) \tilde u \]
defines a bounded operator on $\mc H$. Moreover, $L'$ is compact relative to $(L_0, \mc D( L_0))$.
\end{lemma}

\begin{proof}
First, we observe that the potential satisfies 
\begin{align}\label{Eq:PotBounds1}
\begin{split}
|V^{(2k)}(\rho)| & \lesssim_{k}  \langle \rho \rangle^{-2-2k}, \quad |V^{(2k+1)}(\rho)| \lesssim_k  \rho \langle \rho \rangle^{-4 - 2k},
\end{split}
\end{align}
for all $k \in \N_0$ and $\rho \in [0, \infty)$. We show below that for all $\tilde u \in \mc D(\tilde L_0)$ and $R \geq 1$,
\begin{align}\label{Eq:Pot1}
\| \Lap L' \tilde u \|_{H^1(\B^7)} \lesssim \|\tilde u \|, \quad \| \Lap L' \tilde u \|_{H^1(\R^7\setminus \B^7_R)} \lesssim R^{-1} \|\tilde u \|,
\end{align}
as well as 
\begin{align}\label{Eq:Pot2}
\| \Lap^2 L' \tilde u \|_{L^2(\B^7)} \lesssim \|\tilde u \|, \quad \| \Lap^2 L' \tilde u \|_{L^2(\R^7\setminus \B^7_R)} \lesssim R^{-1} \|\tilde u \|.
\end{align}
These bounds imply that 
\[ \| L' \tilde u \| \lesssim \| \tilde u  \| \]
for all $ \tilde u  \in \mc D(\tilde L_0)$. By density of $\mc D(\tilde L_0)$ in $\mc H$, $L'$ extends to a bounded operator on $\mc H$. 
Next, we set 
$\mc G := (\mc D(L_0), \| \cdot \|_{\mc G})$, where
\[ \|\tilde u \|_{\mc G} := \|\tilde u \| + \|L_0 \tilde u \|, \]
denotes the graph norm. We show below that for all $\tilde u \in \mc D(\tilde L_0)$ and $R \geq 1$, 
\begin{align}\label{Eq:Pot3}
\| \Lap^2 L' \tilde u \|_{\dot H^1(\B^7)} \lesssim  \|\tilde u \|_{\mc G} ,
\quad \| \Lap^2 L'\tilde u \|_{\dot H^1(\R^7 \setminus \B^7_R)} \lesssim  R^{-1}  \|\tilde u \|_{\mc G}.
\end{align}
By definition of the closure, these bounds extend to all of $\mc D(L_0)$. 
In view of Eq.~\eqref{Eq:Pot1} - \eqref{Eq:Pot3} we get that for $j= 1,2$, 
and all $\tilde u \in \mc D(L_0)$, 
\begin{align}\label{Eq:Pot3a}
 \| \Lap^j L' \tilde u \|_{H^1(\R^7)} \lesssim \| \tilde u \|_{\mc G}
\end{align}
and 
\begin{align}\label{Eq:Pot3b}
\| \Lap^j L'\tilde u \|_{H^1(\R^7 \setminus \B^7_R)} \lesssim R^{-1}  \| \tilde u \|_{\mc G}.
\end{align} 

Let  $\mc B_{\mc G} := \{\tilde  u \in \mc D(L_0): \|\tilde  u \|_{\mc G} \leq 1 \}$. 
To see that the perturbation is compact as an operator $L': \mc G \to \mc H$,
we convince ourselves that the sets 
\begin{align*}
 K_1 := \Lap L'(\mc B_{\mc G} ) , \quad  K_2 := \Lap^2 L'(\mc B_{\mc G}) 
\end{align*}
are totally bounded in $L^2(\R^7)$. By Eq.~\eqref{Eq:Pot3a}, $K_1, K_2 \subset H^1(\R^7)$. 
By equation \eqref{Eq:Pot3b}, there exists a constant $C >0$ such that for all $\tilde u \in  \mc B_{\mc G}$ and $R \geq 1$,
\[\| \Lap^j L'\tilde u \|_{H^1(\R^7 \setminus \B^7_R)} \leq C R^{-1}. \]
The right hand side becomes arbitrarily small (uniformly in $\tilde u$) by choosing $R$ large enough. 
Hence, we can apply the result of \cite{HanchHold2010}, Theorem 10,  which implies that $K_1, K_2 \subset L^2(\R^7)$ are totally bounded. 
Now let $(\tilde u_n)_{n \in \N} \subset \mc B_{\mc G}$ and consider the sequence
$(L'\tilde u_n)_{n \in \N}$. Then $(\Lap L' \tilde u_n)_{n \in \N} \subset K_1$. Since $K_1$ is totally bounded in $L^2(\R^7)$, 
there is a subsequence, 
still denoted by $(\Lap L' \tilde u_n)_{n \in \N} $, that is a Cauchy sequence in $L^2(\R^7)$.
By applying the Laplace operator we obtain a sequence $(\Lap^2 L' \tilde u_n)_{n \in \N}  \subset K_2$. Again, we find a subsequence
that converges in $L^2(\R^7)$ and thus, we have identified a subsequence of $(L'\tilde u_n)_{n \in \N}$ that converges in $\mc H$. This implies
the claim. 

\medskip
It is left to prove Eq.~\eqref{Eq:Pot1} - \eqref{Eq:Pot3}. By using the bounds in Lemma \ref{Le:AllBounds}, we get 
\begin{align*}
\| \Lap [V(|\cdot|) \tilde u ] \|_{L^2(\B^7)}  & \lesssim  \| V \|_{L^{\infty}(\R^+)} \| \Lap \tilde u   \|_{L^2(\R^7)}
+  \|(\cdot)^3 \Lap_{\mathrm{rad}} V \|_{L^{\infty}(\R^+)} \| u \|_{L^{2}(\R^+)}
 \\ & +   \|(\cdot) V'\|_{L^{\infty}(\R^+)} \|(\cdot)^2 u' \|_{L^{2}(\R^+)} \lesssim \| \tilde u \| 
\end{align*}
for all $\tilde u \in \mc D(\tilde L_0)$. For $R \geq 1$, we can exploit the decay of the potential to infer that 
\begin{align*}
\| \Lap  [V(|\cdot|) \tilde u ]\|_{L^2(\R^7 \setminus \B^7_R)} & \lesssim  R^{-2} \|\Lap \tilde u \|_{L^2(\R^7)}    + R^{-2}  \|(\cdot) u \|_{L^2(\R^+)}    + R^{-2}  \|(\cdot)^2 u'\|_{L^2(\R^+)} 
\lesssim   R^{-2} \| \tilde u \|.
\end{align*}
An explicit calculation shows that 
\begin{align*}
 \| \nabla \Lap [V(|\cdot|) \tilde u ]\|_{L^2(\R^7)} &  \lesssim 
 \|(\cdot)^2 \langle \cdot \rangle^{-4}  u\|_{L^2(\R^+)} + \|(\cdot) \langle \cdot \rangle^{-2}  u'\|_{L^2(\R^+)} \\
 & +
  \|(\cdot)^2 \langle \cdot \rangle^{-2}  u''\|_{L^2(\R^+)} +  \|(\cdot)^3 \langle \cdot \rangle^{-2}  u^{(3)}\|_{L^2(\R^+)} \lesssim \| \tilde   u \|,
\end{align*}
and thus,
\begin{align*}
 \| \nabla \Lap [V(|\cdot|) \tilde u ]\|_{L^2(\R^7 \setminus \B^7_R)}  \lesssim R^{-2} \| \tilde   u \|.
\end{align*}
This proves Eq.~\eqref{Eq:Pot1}. To get Eq.~\eqref{Eq:Pot2}, we estimate 
\begin{align*}
 \| \Lap^2 [V(|\cdot|) \tilde u ]\|_{L^2(\R^7)} &  \lesssim 
 \|(\cdot) \langle \cdot \rangle^{-4}  u\|_{L^2(\R^+)} + \|\langle \cdot \rangle^{-2}  u'\|_{L^2(\R^+)} \\
 & +
  \|(\cdot) \langle \cdot \rangle^{-2}  u''\|_{L^2(\R^+)} +  \|(\cdot)^2 \langle \cdot \rangle^{-2}  u^{(3)}\|_{L^2(\R^+)} 
 + \|(\cdot)^3 \langle \cdot \rangle^{-2}  u^{(4)}\|_{L^2(\R^+)}  \lesssim \| \tilde   u \|,
\end{align*}
From this calculation it follows easily that 
\begin{align*}
 \| \Lap^2 [V(|\cdot|) \tilde u ]\|_{L^2(\R^7 \setminus \B^7_R)}  \lesssim R^{-2} \| \tilde   u \|.
\end{align*}
It is left to prove Eq.~\eqref{Eq:Pot3}. First, we observe that 
\begin{align*}
 \| \nabla \Lap^2 [V(|\cdot|) \tilde u ]   \|_{L^2(\R^7)} &  \lesssim 
\|V(|\cdot|) \nabla  \Lap^2 \tilde  u \|_{L^2(\R^7)} + \| \langle \cdot \rangle^{-4} u \|_{L^2(\R^+)}   
+  | (\cdot) \langle \cdot \rangle^{-4} u' \|_{L^2(\R^+)}  \\
& + \| (\cdot)^2 \langle \cdot \rangle^{-4} u'' \|_{L^2(\R^+)}   +  \| (\cdot)^{3}  \langle \cdot \rangle^{-4}  u''' \|_{L^2(\R^+)} 
+ \|(\cdot)^{4} \langle \cdot \rangle^{-4} u^{(4)} \|_{L^2(\R^+)} \\
& \lesssim \|V(|\cdot|) \nabla  \Lap^2 \tilde  u \|_{L^2(\R^7)} + \| \tilde   u \|,
\end{align*}
and 
\begin{align*}
 \| \nabla \Lap^2 [V(|\cdot|) \tilde u ]   \|_{L^2(\R^7 \setminus \B^7_R)} 
 \lesssim \|V(|\cdot|) \nabla  \Lap^2 \tilde  u \|_{L^2(\R^7 \setminus \B^7_R)} + R^{-2} \| \tilde   u \|.
\end{align*}
For $\tilde u \in \mc D(\tilde L_0)$, we have
\[ \| \nabla \Lap \tilde L_0 \tilde  u \|_{L^2(\R^7)} \lesssim \| \tilde L_0 \tilde u \|,\]
by Lemma \ref{Le:PropH}. Hence,
\begin{align*}
\|V(|\cdot|)  & \nabla  \Lap^2 \tilde  u  \|_{L^2(\R^7)}   \lesssim \| V(|\cdot|)\nabla \Lap \tilde L_0 \tilde u \|_{L^2(\R^7)}
+ \|  V(|\cdot|) \nabla \Lap \Lambda \tilde u  \|_{L^2(\R^7)}   + \|  V(|\cdot|)  \nabla \Lap \tilde u \|_{L^2(\R^7)}. 
\end{align*} 
With Lemma \ref{Le:AllBounds} one easily gets that 
\begin{align}\label{Eq:NablaLapLamba}
\|\langle \cdot \rangle^{-1}  \nabla \Lap \Lambda  \tilde u \|_{L^2(\R^7)} \lesssim  \|\tilde u\| ,
\end{align} 
for all $\tilde u \in \mc D(\tilde L_0)$.
This and the decay of the potential at infinity yield Eq.~\eqref{Eq:Pot3}.
\end{proof}

The next statement is a consequence of Lemma \ref{Le:Compact1} and the Bounded Perturbation Theorem, \cite{engel}, p.~158.

\begin{corollary}\label{Cor:SemigroupS}
The operator $L := L_0 + L'$, $\mc D(L) = \mc D(L_0)$ generates a strongly continuous one-parameter 
semigroup $\{ S(\tau) : \tau \geq 0\}$ of bounded operators on $\mc H$ satisfying
\[ \|S(\tau) \tilde  u \| \leq e^{(\| L' \| - \frac{1}{4}) \tau} \| \tilde u \|\]
for all  $ \tilde u \in \mc H$ and all $\tau \geq 0$. 
\end{corollary}

Unfortunately, the growth bound in Corollary \ref{Cor:SemigroupS} is too weak for our purpose and we have to analyze 
the spectrum of the generator $(L, \mc D(L))$ in order to improve it. For this, we exploit 
the compactness of the perturbation relative to $L_0$ to infer that spectral points with $\mathrm{Re} \lambda > - \frac{1}{4}$
must be eigenvalues. The time-translation invariance of the problem
induces the unstable eigenvalue $\lambda = 1$. To prove that this is the only eigenvalue with 
non-negative real part, we use information on the spectrum of the supersymmetric partner of the 
perturbed Ornstein-Uhlenbeck operator, cf.~\cite{CostinDonn}, in the self-adjoint setting.

\subsection{Spectral analysis for a self-adjoint operator}\label{Sec:Superpartner}

We introduce the formal differential expression 
\begin{align}\label{Eq:DefAlpha}
\mc A v(\rho) = - v''(\rho) + \left (\frac{\rho^2}{16}  + \frac{12}{\rho^2} +\frac34  + Q(\rho) \right ) v(\rho), 
\end{align}
with
\[ Q(\rho) = \frac{384 \sqrt{6} - \rho^2 \left(\rho^2+24 \sqrt{6}-44\right) -956}{\left(\rho^2+6
 \sqrt{6}-14\right)^2} . \]

\begin{lemma}\label{Le:Supersymm_Op}
Let
\[ \mc D(A) = \{v \in L^2(\R^+): v,v' \in AC_{\mathrm{loc}}(\R^+), \mc A v \in L^2(\R^+)  \}\]
and set $Av = \mc A v$ for $v \in \mc D(A)$. 
Then, $(A,\mc D(A))$ is a self-adjoint operator on $L^2(\R^+)$ and $\sigma(A) \subseteq [\omega_A,\infty)$ for $\omega_A = \frac{1}{75}$.
\end{lemma}

\begin{proof}
First, we convince ourselves that $\mc A$ is limit-point at both endpoints of the  interval $(0,\infty)$.
For $\rho = 0$ we consider the equation 
\begin{align}\label{Eq:EigenvalueA}
\lambda v(\rho) + v''(\rho) - q(\rho) v(\rho) = 0,
\end{align}
for 
\[q(\rho) := \frac{\rho^2}{16}  + \frac{12}{\rho^2} +\frac34  + Q(\rho)  \]
and some fixed $\lambda \in \C$. The origin is a regular singular point and we can use 
Frobenius' method to obtain the existence of a fundamental system around zero given by $\{ v_0 , v_1 \}$,
\[ v_0(\rho) = \rho^{4} h_0(\rho), \quad  v_1(\rho) = c h_0(\rho) \rho^4  \log \rho + \rho^{-3} h_1(\rho), \]
where $c \in \R$ is some constant and $h_0, h_1$ are analytic functions with $h_0(0) = h_1(0) = 1$. 
Clearly, $v_1$ does not belongs to $L^2(0,\delta)$ for any $\delta > 0$ and the Weyl alternative implies that $\mc A$ is
limit point at $\rho = 0$.  For the other endpoint, we use the fact that $\lim_{\rho \to \infty} q(\rho)= \infty$, cf.~ \cite{Weidmann}, Theorem 6.6, p.~96.
We infer that the (maximal) operator defined as above is self-adjoint. Furthermore, a core is given by
 $(\tilde A,\mc D(\tilde A))$, where $\tilde A v := \mc A u$ and 
\[ \mc D(\tilde A )=\{v \in \mc D(A): v \text{ has compact support}\}.\]

We show that 
\begin{align}\label{Eq:BoundA}
(\tilde A v,v)_{L^2(\R^+)}  \geq \omega_A \|v\|^2_{L^2(\R^+)} 
\end{align}
for $\omega_A = \frac{1}{75}$ and all $v \in \mc D(\tilde A )$. 
Since $\tilde A$ is a core of $A$, this implies that 
\begin{align*}
( A v,v)_{L^2(\R^+)}  \geq \omega_A \|v\|^2_{L^2(\R^+)} 
\end{align*}
for all $v \in \mc D( A )$.
From this, we get that the spectrum of $A$ is bounded below and that $\sigma(A) \subset  [\omega_A,\infty)$, see for example \cite{kato}, p.~278.
The proof of Eq.~\eqref{Eq:BoundA} depends crucially on the properties of $q$. First, we 
use the fact that $v \in \mc D(\tilde A )$ vanishes at the origin and apply the Cauchy-Schwarz inequality to get that
\begin{align}\label{Eq:Proof_SuSy}
\int_0^{\gamma} |v(\rho)|^2 d\rho \leq \gamma^2 \int_0^{\gamma} |v'(\rho)|^2 d\rho
\end{align}
for all $\gamma > 0$. One can easily check that the function $q$ attains
its global minimum $q(\rho_{min})=q_{min}$ at some $\rho_{min} \in (0,\gamma)$, $\gamma := \frac52$, and that
\[  \gamma^{-2} + q_{min} > q(\gamma) > \omega_A.\]
On $[\gamma, \infty)$, $q$ is strictly positive and monotonically increasing. 
Using integration by parts we estimate
\begin{align*}
(\tilde A & v| v)_{L^2(\R^+)}   =   - \int_0^{\infty} v''(\rho) \overline{v(\rho)} d\rho + \int_0^{\infty} q(\rho) |v(\rho)|^2 d\rho \\
		   & =   \int_0^\gamma |v'(\rho)|^2 d\rho  + \int_\gamma ^{\infty} |v'(\rho)|^2 d\rho   + \int_0^{\gamma} q(\rho) |v(\rho)|^2 d\rho   +\int_\gamma^{\infty} q(\rho) |v(\rho)|^2 d\rho \\
		   & \geq \left(\gamma^{-2} + q_{min} \right) \int_0^{\gamma} |v(\rho)|^2 d\rho + q(\gamma) \int_\gamma^{\infty} |v(\rho)|^2 d\rho  
		     \geq   q(\gamma) \int_0^{\infty}|v(\rho)|^2 d\rho.
\end{align*}
This implies the claim. 
\end{proof}

\subsection{Characterization of the spectrum of $L$}\label{Sec:SpectrumL}

\begin{lemma}\label{Le:Spec}
Let $\lambda \in \C$ be a spectral point of the operator $(L, \mc D(L))$. 
Then either $\mathrm{Re} \lambda \leq - \frac{1}{75}$ or $\lambda = 1$ is an eigenvalue.
\end{lemma}

\begin{proof}
First, we note that in polar coordinates the equation $(\lambda - L)\tilde u =0$
reduces to 
\begin{align}\label{Eq:EigenvalueEq}
\lambda u(\rho) - u''(\rho) - \frac{6}{\rho} u'(\rho) + \frac{1}{2} \rho u'(\rho) + (1 - V(\rho))u(\rho) = 0.
\end{align}
Smoothness of the coefficients on $(0,\infty)$ implies that eigenfunctions are smooth
away from zero. By Frobenius' method there exists a  
fundamental system $\{ u^0, u^1 \}$ around zero, where $u^0$ is analytic and satisfies $u^0(0) =1$. 
The second solution is given by
\[  u^1(\rho) = c u^0(\rho)  \log  \rho+ \rho^{-5} h(\rho) \]
for some $c \in \R$ and some analytic function $h$ with $h(0)=1$. For an eigenfunction $u$, we have that
$u \in C[0,\infty)$, and thus, $u = \alpha u^0$ for some $\alpha \in \C$. This shows that 
$u$ is smooth on $[0,\infty)$.

For $\lambda = 1$, a direct calculation shows that  
\[\mb g(\rho) = (a_1 \rho^2 + a_2)^{-2}\]
solves Eq.~\eqref{Eq:EigenvalueEq}, where the constants $a_1,a_2$ are given in Eq.~\eqref{Eq:ConstantsWein5d}.
By exploiting the decay of $\mb g$ it follows that $\tilde{\mb g}$ can be approximated by functions in $C^\infty_{0,\mathrm{rad}}(\R^7)$, with respect to the norm $\|\cdot\|$.
Consequently, $\mb{\tilde g}:= \mb g(|\cdot|) \in \mc H$ and  
it is easy to check that $\mb{\tilde g} \in \mc D(\tilde L_0) \subset \mc D(L_0)$. Hence,
$\mb{\tilde g}$ is an eigenfunction.

To prove the rest of Lemma \ref{Le:Spec}, we assume that $\lambda \in \sigma(L)$. 
If $\mathrm{Re} \lambda \leq - \frac{1}{75}$, then the statement is true. 
So let $\mathrm{Re} \lambda > - \frac{1}{75}$. By standard results from semigroup theory and Lemma \ref{Cor:SemigoupS0},
$\lambda \not \in \sigma(L_0)$ and the resolvent $R_{L_0}(\lambda) \in \mc B(\mc H)$ exists.  
The identity $\lambda - L = (1 - L'R_{L_0}(\lambda))(\lambda - L_0)$ shows that
$1 - L'R_{L_0}(\lambda)$ is not bounded invertible. The operator $L'R_{L_0}(\lambda)$ is compact
by the compactness of $L'$ relative to $L_0$, hence, $1 \in \sigma(L'R_{L_0}(\lambda))$ is an eigenvalue
and there exists an eigenfunction $\tilde f \in \mc H$ such that $(1 - L'R_{L_0}(\lambda)) \tilde f = 0$.
Set $\tilde u = R_{L_0}(\lambda) \tilde f$. Then, $\tilde u \in \mc D(L)$ and 
$(\lambda - L) \tilde u = (1 - L'R_{L_0}(\lambda))(\lambda - L_0) \tilde u = 0$. 
Thus, $\lambda$ must be an eigenvalue.
The case $\lambda = 1$ has already been discussed, so assume that $\lambda \neq 1$
and let $\tilde u_{\lambda} =u_{\lambda}(|\cdot|)$ denote the corresponding eigenfunction. 
We set 
\[ v_{\lambda}(\rho) = \rho^{3} e^{-\frac{\rho^2}{8}} u_{\lambda}(\rho). \]
Since $u_{\lambda} \in C^{\infty}[0,\infty)$ is a solution to Eq.~\eqref{Eq:EigenvalueEq}, $v_{\lambda}$ is smooth and solves the equation
\begin{align}\label{Eq:EigenvalueNormal}
\lambda v_{\lambda}(\rho) - v_{\lambda}''(\rho) +  \left( \frac{\rho^2}{16}  + \frac{6}{\rho^2} - \frac34  - V(\rho) \right) v_{\lambda}(\rho) = 0.
\end{align}

We define the differential expressions
\begin{align}\label{Def:Factorization}
 Bv(\rho) := -v'(\rho) + \beta (\rho) v(\rho), \quad B^+ v(\rho) := v'(\rho) + \beta (\rho) v(\rho), 
\end{align} 
for 
\begin{align*}
\beta(\rho) := \frac{3}{\rho} - \frac{\rho}{4} - \frac{4 \rho}{6 \sqrt{6} -14 + \rho^2}.
\end{align*}
Using this, Eq.~\eqref{Eq:EigenvalueNormal} can be written as 
\begin{align}\label{Eq:NormFormFac}
(\lambda + B^+ B - 1)v_{\lambda}(\rho) = 0. 
\end{align}
It is easy to see that the kernel of $B$ is spanned by the transformed symmetry mode $\rho \mapsto \rho^{3} e^{-\frac{\rho^2}{8}} \mb g(\rho)$.
We set $w_{\lambda} :=  B v_{\lambda}$. 
By applying $- B$ to Eq.~\eqref{Eq:NormFormFac} we infer that $w_{\lambda}$ satisfies the equation
\begin{align}\label{Eq:Supersym}
(-\lambda - B B^+ + 1)w_{\lambda}(\rho) = 0.
\end{align}
A straightforward calculation shows that Eq.~\eqref{Eq:Supersym} can be written as
\[ -\lambda w_{\lambda}(\rho)  - \mc A(\rho)w_{\lambda}(\rho) = 0, \]
where $\mc A$ is given by Eq.~\eqref{Eq:DefAlpha}. Recall that $\mc A$ is the defining differential expression
for the self-adjoint operator $A$ described in Lemma \ref{Le:Supersymm_Op}.
By inserting the definition we get that 
\[ w_{\lambda}(\rho) = - \rho^{3} e^{-\frac{\rho^2}{8}} u_{\lambda}'(\rho) + \rho^2 e^{-\frac{\rho^2}{8}} h(\rho) u_{\lambda}(\rho), \]
where $h \in C^{\infty}[0,\infty)$, $\lim_{\rho \to \infty} h(\rho) = c$. Obviously, $w_{\lambda} \in C^{\infty}[0,\infty)$
and with Lemma \ref{Le:PropH} it is easy to see that $w_{\lambda}  \in L^2(\R^+)$ and thus 
$w_{\lambda} \in \mc D(A)$. Hence, $\tilde \lambda := - \lambda$ is an eigenvalue of the operator $(A, \mc D(A))$ with eigenfunction 
$w_{\lambda}$. Our assumption on $\lambda$ implies that $\mathrm{Re} \tilde \lambda < \frac{1}{75}$. However, this contradicts Lemma \ref{Le:Supersymm_Op}. 
We conclude that there are no spectral points of the operator $L$ with $\mathrm{Re} \lambda > - \frac{1}{75}$, except for $\lambda = 1$. 
\end{proof}

\begin{lemma}\label{Le:Ker1mL}
We have that $\mathrm{ker} ( 1 - L) = \mathrm{span}(\mb{\tilde g}) $, where 
$\mb{ \tilde g} = \mb g(|\cdot|)$ and
\[\mb g(\rho) = \frac{1}{(a_1 \rho^2 + a_2)^{2}},\]
for $a_1 = \frac12 \sqrt{\frac{3}{2}}$, $a_2 = \frac12 (18 - 7 \sqrt{6})$. 
\end{lemma}

\begin{proof}
It was already shown above that $\mb{\tilde g}$ is an eigenfunction
corresponding to the eigenvalue $\lambda = 1$. Assume that there is another eigenfunction $\tilde u = u(|\cdot|) \in
\mathrm{ker} ( 1 - L)$. Then $u$ is a solution to the ODE 
\begin{align}\label{Eq:SymmEigenvalueEq}
u''(\rho) + \tfrac{6}{\rho} u'(\rho) - \tfrac12 \rho u'(\rho) + V(\rho) u(\rho) - 2 u(\rho) = 0.
\end{align}
For this equation a fundamental system is given by $\{ \mb g, \mb h \}$, where 
$\mb g$ is as above and
\begin{align}\label{Eq:SecSolSymmEigenvalEq}
 \mb h(\rho) = e^{\frac{\rho^2}{4}} \left( h_1(\rho) + h_2(\rho) e^{-\frac{\rho^2}{4}} \int_0^{\rho} e^{\frac{s^2}{4}} ds \right),
\end{align}
\[ h_1(\rho) = \frac{\sum_{j=0}^3 \alpha_j \rho^{2j}}{20\rho^5(6 \sqrt{6} - 14 + \rho^2)^2}, \quad h_2(\rho) = \frac{2(61-24 \sqrt{6})}{5 (6 \sqrt{6} - 14 + \rho^2)^2}, \]
for constants $\alpha_0 = 24 ( 8652 \sqrt{6} - 21193)$, $\alpha_1  = 4 ( 8347-3408 \sqrt{6})$, $\alpha_2 = 2 (372 \sqrt{6} - 923)$, $\alpha_3 = 15$. 
Their Wronskian is given by 
\[W(\mb g,\mb h)(\rho) = \rho^{-6} e^{\frac{\rho^2}{4}}.\]
The function $\mb h$ can also be written as 
\[\mb h(\rho) = \rho^{-5} \langle \rho \rangle^2 e^{\frac{\rho^2}{4}} H(\rho),\]
where $H$ is regular around zero, $H(0) \neq 0$ and 
$\lim_{\rho \to \infty} H(\rho) = c$ for some $c \in \R\setminus\{0\}$. Hence,
\[ u(\rho) = c_1 \mb g(\rho) + c_2 \mb h(\rho), \]
for some constants $c_1,c_2 \in \C$. Since $\mb h$ diverges at the origin as well as for $\rho \to \infty$, 
we must have $c_2 = 0$ to guarantee that $u(|\cdot|) \in \mc D(L) \subset \mc H$. 
Thus, $u$ is a multiple of $\mb g$, which implies the claim.
\end{proof}

\section{Resolvent estimates}\label{Sec:Resolventest}

We need suitable bounds for the resolvent $R_L(\lambda)$ in order to translate the spectral information into growth estimates for the semigroup generated by $L$.
First, note that
\[ \mathrm{Re} (L \tilde u| \tilde u) \leq (\| L'\|  -\tfrac14 ) \|\tilde u \|^2 ,\]
for all $\tilde u \in \mc D(L)$. By standard semigroup theory,
\[ \|R_{L}(\lambda)\| \leq \frac{1}{\mathrm{Re} \lambda - \|L'\| +\frac14}, \]
for all $\lambda \in \C$ with $\mathrm{Re} \lambda > \| L'\| -\frac14 $.
This shows uniform boundedness of the resolvent in a right half plane for $\mathrm{Re} \lambda$ 
sufficiently large. The aim of this section is to improve this considerably by establishing the following result.  

\begin{proposition}\label{Prop:ResolventBounds}
Fix $\alpha > - \frac{1}{75}$ and choose $M_{\alpha} > 0$ sufficiently 
large. Define 
\[ \Omega_{\alpha} := \{\lambda \in \C: \lambda = \alpha + \I \omega, \text{ where } \omega \in \R, |\omega| \geq M_{\alpha} \}. \]
There exists a constant  $C_{\alpha} >0$ such that 
\begin{align}\label{Eq:ResolventBounds}
\| R_L(\lambda ) \tilde f \|  \leq C_{\alpha}  \|\tilde f \|
\end{align}
for all $\lambda \in \Omega_{\alpha}$ and $\tilde f \in \mc H$.
\end{proposition}

To prove Proposition \ref{Prop:ResolventBounds}, we use the machinery developed by the authors in \cite{DonSch14a}. 
There, an analogous statement was established for perturbed (non-radial) Ornstein-Uhlenbeck operators on $L^2(\R^d)$. 
The situation here is to some extent simpler, since we only have to deal with a radial problem in $\R^7$.
However, in contrast to \cite{DonSch14a} we need uniform bounds on weighted derivatives of the resolvent up to the fourth order.
This requires additional work. 

To keep formulas within margin, we introduce the following useful notation.

\begin{definition}
For a function $f: I \subset \R \to \C$ and $\gamma \in \R$, we write $f(x) = \mathcal O(x^{\gamma})$ if
\[ | f^{(k)}(x)| \leq C_k |x|^{\gamma - k}\]
for all $x \in I$ and $k \in \N_0$.  Similarly, 
$f(x) = \mathcal O(\langle x \rangle^{\gamma})$, if 
\[|f^{(k)}(x)| \leq C_k \langle x \rangle^{\gamma - k},\]
for all $x \in I$ and $k \in \N_0$. 
Functions with this property are said to be of \textit{symbol type} or to have \textit{symbol behavior}. We note that symbol behavior is stable
under algebraic operations, e.g., $\Oc(x^{\beta})\Oc (x^{\gamma}) = \Oc(x^{\beta + \gamma})$ 
for $\beta, \gamma \in \R$. An analogous definition holds for functions depending on more than one variable. 
\end{definition} 

\subsection{Reduction to an ODE problem}\label{SubSecRes:ODEform}
If $\lambda \in \C$, $\lambda \neq 1$ and $\mathrm{Re} \lambda > - \frac{1}{75}$, then we know from Lemma \ref{Le:Spec} that 
$\lambda \in \rho(L)$ and that the resolvent $R_L(\lambda): \mc H \to \mc D(L)$ exists as a bounded operator. For 
$\tilde f = f(|\cdot|) \in \mc H$ we set $\tilde w:= R_L(\lambda) \tilde  f$, $\tilde w = w(|\cdot|)$. By definition,
$(\lambda - L) \tilde w = \tilde f$, which reduces to 
\begin{equation}\label{Eq:ODEResolvent}
\lambda w(\rho) - w''(\rho) - \frac{6}{\rho} w'(\rho) + \frac{1}{2} \rho w'(\rho) + w(\rho) - V(\rho) w(\rho)  = f(\rho),
\end{equation}
in polar coordinates. Constructing the resolvent explicitly amounts to solving Eq.~\eqref{Eq:ODEResolvent}
for given right hand side. 
We are left with a pure ODE problem, which is however challenging due to the nontrivial form of the potential.
In order to apply some of the outcomes of \cite{DonSch14a}, we slightly
reformulate our problem. By setting $u(r) := w(2 \rho)$ and changing the sign, the above equation transforms into
\begin{equation}\label{Eq:ODEResolvent_Rescaled}
u''(r) + \frac{6}{r} u'(r) - 2 r u'(r) +  \tilde V(r)  u(r) - \tilde \lambda  u(r)  = - 4 f(2 r) 
\end{equation}
for $\tilde \lambda = 4(\lambda+1)$, $\tilde V(r) = 4 V(2 \rho)$.

\subsection{Construction of a fundamental system for the homogeneous equation}\label{SubSecRes:HomEq} 
We consider the homogeneous version of Eq.~\eqref{Eq:ODEResolvent_Rescaled} which is
\begin{equation}\label{Eq:HomODE}
u''(r) + \frac{6}{r} u'(r) - 2 r u'(r) +  \tilde V(r)  u(r) - \tilde \lambda  u(r)  = 0,
\end{equation}
for $r  > 0$. This corresponds to Eq.~$(4.1)$ in  \cite{DonSch14a} for $d=7$ and $\ell = 0$.
Concerning the potential, we only need that $\tilde V(r) = \Oc(\langle r \rangle^{-2})$.
By setting $v(r) = r^{3} e^{-\frac{r^2}{2}} u(r)$, 
we obtain the equation
\begin{equation}\label{Eq:HomODENorm}
v''(r)  - r^2 v(r) - \frac{4\nu^2-1}{4 r^2} v(r) - \mu v(r)  +  \tilde V(r) v(r) =0
\end{equation}
for $\nu =\frac{5}{2}$ and $\mu := \tilde \lambda - 7 $. 
For $r \gg 1$, Eq.~\eqref{Eq:HomODENorm} resembles
a Weber equation, whereas for $r$ small we have a perturbed Bessel equation.
For the construction of a suitable fundamental system, we proceed 
along the lines of \cite{DonSch14a} with some simplifications due to the fact that $\nu$ is fixed in our case. 
For the rest of this section, we assume that 
\begin{align}\label{Eq:mu}
\mu = b + \I \omega
\end{align}
with $b > -4$ fixed, corresponding to $\mathrm{Re} \lambda > - \frac{1}{4}$. This is actually more than we need, since $\mathrm{Re} \lambda > - \frac{1}{75}$ and
thus $b > -\frac{229}{75} \sim -3.05$ would be sufficient. 
However, we expect that the operator $L$ has no eigenvalues with $\mathrm{Re} \lambda > - \frac{1}{4}$ and large imaginary parts, which motivates our choice. 

Finally, we note that most implicit constants depend on $b$ in the following. However, we do not indicate this dependence in order to improve readability.
All other dependencies are tracked during the calculations.

\subsubsection{A fundamental system away from the center}
To motivate the following, we write Eq.~\eqref{Eq:HomODENorm} as
\begin{align}\label{Eq:WeberInhom}
 v''(r)  - (r^2 + \mu) v(r)  = \Oc(r^{-2}) v(r)
\end{align}
for $r \geq 1$. In a first step, we construct a fundamental system for the equation  
\begin{equation}\label{Eq:Weberhom}
v''(r)  - (r^2 + \mu) v(r)  = 0.
\end{equation}
We note that Eq.~\eqref{Eq:Weberhom} can be solved in terms of parabolic cylinder functions. 
However, this is not very useful for our purpose since we need precise information on
the asymptotics for large imaginary parts. We use instead a Liouville-Green transform in combination
with perturbation theory. To motivate the following, we first assume that $\mu \in \R$, $\mu >0$. 
With $w(y)=v(\mu^{\frac12} y)$, Eq.~\eqref{Eq:Weberhom} transforms into 
\begin{equation}\label{Eq:Weberhom_T1}
w''(y)  - \mu^2 (1  + y^2) w(y)  = 0.
\end{equation}
We set 
\[h(\zeta(y)):=|\zeta'(y)|^{\frac12} w(y)\]
for some diffeomorphism $\zeta$ that has yet to
be defined. This leads to 
\[ h''(\zeta(y))-\frac{\mu^2(1+y^2)}{\zeta'(y)^2}
h(\zeta(y))
-\frac{q(y)}{\zeta'(y)^2}h(\zeta(y))=0
 \]with  
\[ q(y)=\frac12 \frac{\zeta'''(y)}{\zeta'(y)}-\frac34 \frac{\zeta''(y)^2}{\zeta'(y)^2}. \]
We require that $\zeta'(y) = \sqrt{1+y^2}$ and define
\begin{align}\label{Eq:DefZeta}
\begin{split}
\zeta(y) :=   \int_{10 \mu^{-\frac12} }^y  & \sqrt{1+z^2} dz   = F(y) - F(10 \mu^{-\frac12}),
\end{split}
\end{align}
where $F(z) = \frac12 \log(z+\sqrt{1+z^2})+ \frac12 z \sqrt{1+z^2}$. With this choice we obtain
\[ q(y) = \frac{2 - 3 y^2}{ 4(1+y^2)^2}. \]
Having this, we add $q w$ to Eq.~\eqref{Eq:Weberhom_T1}, i.e., we consider
\begin{equation}\label{Eq:Resc_plQ}
w''(y)  - \mu^2 (1  + y^2) w(y) + q(y) w(y) =0
\end{equation}
which transforms now into 
\[ h''(\zeta(y)) - \mu^2 h (\zeta(y)) = 0.\]
For this equation we can explicitly write down a fundamental system $\{h^{+}, h^{-} \}$, 
given by $h^{\pm} (\zeta(y)) = e^{\pm \mu \zeta(y)}$. 
Hence, Eq.~\eqref{Eq:Resc_plQ} has a fundamental system  $\{ w^{+},  w^{-} \}$,
\[w^{\pm}(y) =  \zeta'(y)^{-\frac12} e^{\pm \mu \zeta(y)}.\] 
By setting $r = \mu^{\frac12} y$ and $w(y)=v(\mu^{\frac12} y)$,  Eq.~\eqref{Eq:Resc_plQ} transforms into 
\[ v''(r)  - (r^2 + \mu) v(r) + \mu^{-1} q(\mu^{-\frac12} r) v(r) =0, \]
for which we now have the fundamental system $\{ v^+, v^-\}$, 
\[ v^{\pm}(r) = \zeta'(\mu^{-\frac12} r)^{-\frac12} e^{\pm \mu \zeta(\mu^{-\frac12} r)}.\] 

These calculations suggest to add $\mu^{-1} q(\mu^{-\frac12} r) v(r)$ to both sides of Eq.~\eqref{Eq:WeberInhom} 
and to treat the right hand side perturbatively.

\medskip
However, since we are interested in complex values of $\mu$, 
we have to extend the above quantities
to the complex plane. First, we define $\sqrt{\cdot} = (\cdot)^{\frac12}$ to
be the principal branch of the square root, which is holomorphic in $\C \setminus (-\infty, 0]$. 
More explicitly,
\begin{equation}
\label{eq:sqrt}
 \sqrt z=\frac{1}{\sqrt 2}\sqrt{|z|+\Re z}+\frac{\I \sgn(\Im z)}{\sqrt 2}\sqrt{|z|-\Re z},
 \end{equation}
for $z \in \C \setminus (-\infty, 0]$. We note that 
$\sqrt{z}^2=z$ and $|\sqrt{z}|=\sqrt{|z|}$ for all $z\in \C\backslash (-\infty,0]$.
Furthermore, 
\[ \sqrt{z}\sqrt{w}=\sqrt{zw} \] holds
at if $-\pi < \arg{z} + \arg{w} < \pi$.
The function $F$ is defined and holomorphic on $\C \setminus ([\I,\I\infty) \cup [-\I,-\I\infty))$. 
Since $\mathrm{Re}(\mu^{-\frac12} r) > 0$ for all $\mu  \in \C \setminus (-\infty,0]$, the function
$\mu \mapsto  \zeta(\mu^{-\frac12} r)$ is holomorphic on $\C \setminus (-\infty,0]$.

For the rest of this section we assume that $\mu = b + \I \omega$, where $b > -4$ is fixed and
$\omega > 0$. Furthermore, we assume that $r \geq 3$ for technical reasons. 
The next result is crucial for proving Proposition \ref{Prop:ResolventBounds}. 
\begin{lemma}\label{Le:RepXiQ}
For $r \in [3, \infty)$, we define 
\[ Q(r,\mu) := \mu^{-1} q(\mu^{-\frac12} r). \]
The function $Q(r,\cdot)$ is holomorphic and $Q(r,\mu) = \Oc(r^{-2} \omega^0)$ for $\omega \gg 1$. 
Furthermore, we define 
\begin{align*}
\begin{split}
\xi(r,\mu) :=  \zeta(\mu^{-\frac12} r),
\end{split}
\end{align*}
such that $\xi(r,\cdot)$ is holomorphic and 
 \[\partial_r \xi(r, \mu) = \mu^{-\frac{1}{2}} \sqrt{ 1 + \tfrac{r^2}{\mu}}.\]
For $\omega \gg 1$, we have the representations
\begin{align}\label{Eq:RepXi_large_alt}
\mathrm{Re} [\mu \xi(r,\mu)] & = \mathrm{Re}\mu^{\frac12}(r - 10) + \tilde \varphi(r,\omega),
\end{align}
and 
\begin{align}\label{Eq:RepXi_large}
\mathrm{Re} [ \mu \xi(r,\mu)] =  \tfrac12 r^2 + \tfrac{b}{2} \log \langle\omega^{-1/2} r \rangle  + \varphi(r,\omega),
\end{align}
where both  $\tilde \varphi(\cdot,\omega)$ and $\varphi(\cdot,\omega)$ are monotonically increasing functions.
Moreover, $|\varphi(10,\omega)| \lesssim 1$, $\tilde \varphi(10,\omega) = 0$ and $\tilde \varphi(r,\omega) = \mc O(r^3 \omega^{-\frac12})$ provided $r \omega^{-\frac12} \lesssim 1$. 
\end{lemma}

\begin{proof}
First, we note that for $\omega \gg 1$, we have $|\mu| \simeq \omega$.
To estimate the Liouville-Green potential we use that  
\begin{align}\label{Eq:EstLGPot}
| 1 + \mu^{-1} r^2|^2 \gtrsim  1+ |\mu|^{-2} r^4  \gtrsim   1+ \omega^{-2} r^4  \gtrsim  \langle \omega^{-\frac12} r \rangle^{4}  
\end{align}
for $\omega \gg 1$, see the proof of Lemma $4.1$ in \cite{DonSch14a}. Thus, 
\[ |Q(r,\mu)| \lesssim |\mu|^{-1} \frac{1 + |\mu|^{-1} r^2}{|1 + \mu^{-1} r^2|^2} \lesssim |\mu|^{-1} \langle \omega^{-\frac12} r \rangle^{-2} \lesssim r^{-2} \]
and we obtain $Q(r,\mu) = \Oc(r^{-2} \omega^0)$.

To see that Eq.~\eqref{Eq:RepXi_large_alt} holds, we use that 
$\mathrm{Re}[\mu \partial_r  \xi(r,\mu) ]  = \mathrm{Re} \sqrt{ \mu + r^2}.$ 
With Eq.~\eqref{eq:sqrt} we get  
\begin{align}\label{Eq:Realpart_DerXi}
\mathrm{Re} \sqrt{ \mu + r^2} = \frac{1}{ \sqrt{2}} \sqrt{ \sqrt{ \omega^2 + ( r^2 + b)^2 } + r^2 + b}.
\end{align}  
 Eq.~\eqref{Eq:RepXi_large_alt} is a consequence of the fact that $(r^2 + b)^2>  b^2$ for $b > -4$ and all $r \geq 3$, which implies 
\[ \partial_r \tilde \varphi(r,\omega) = \mathrm{Re}[\mu \partial_r  \xi(r,\mu) ] - \mathrm{Re} \mu^{\frac12}  =
 \mathrm{Re}[\mu \partial_r  \xi(r,\mu) ] - \frac{1}{ \sqrt{2}} \sqrt{ \sqrt{ \omega^2 + b^2 }  + b} > 0.\] 
Since $\xi(10,\mu) = 0 $,  $\tilde \varphi(10,\omega) = 0$ by definition.

If $|r\mu^{-\frac12}| \lesssim 1$, a Taylor expansion
of $F(z)$ around zero yields 
\[ \mu \xi(r,\mu) = \mu^{1/2}(r - 10) + \psi(r,\mu) - \psi(10,\mu), \]
where $\psi(r,\mu)  = \Oc(r^3 \omega^{-\frac12})$. By definition of $\tilde \varphi$, we have $\tilde \varphi(r,\omega) = \mathrm{Re}[\psi(R,\mu) - \psi(10,\mu)]$.
We now turn to Eq.~\eqref{Eq:RepXi_large}.
The right hand side of Eq.~\eqref{Eq:Realpart_DerXi} implies that 
\[\mathrm{Re}[\mu \partial_r  \xi(r,\mu) ]  \geq \sqrt{r^2+b} \]
and for $b \geq 0$,  the proof of Eq.~\eqref{Eq:RepXi_large} is along the lines of \cite{DonSch14a}, p.~2497.
For $-4 < b < 0$ we distinguish two cases. First, we assume that $3 \leq r \leq \sqrt{2 \omega}$.
Then, $4 |b| r^2 \leq 8 |b| \omega < \omega^2$ for $\omega \gg 1$. This implies 
\[ \omega^2 + (r^2- |b|)^2  = \omega^2 + r^4 - 2 r^2 |b| + |b|^2 > (r^2 + |b|)^2.\]
With this, we infer $\mathrm{Re}[\mu \partial_r  \xi(r,\mu) ]  > r $.
In particular,
\[\partial_r \varphi(r,\omega) = \mathrm{Re}[\mu \partial_r  \xi(r,\mu) ] - r + \frac{|b|}{2} \frac{r \omega^{-1}}{ 1 + \omega^{-1} r^2} > 0. \]
For $r > \sqrt{2 \omega}$ we have  $\frac{\omega}{r^2 - |b|} \leq 1$,
and since $\sqrt{1+x^2} \geq 1 + \frac{x^2}{4}$ for $0 \leq x \leq 1$ we infer
\begin{align}
 \mathrm{Re}[\mu \partial_r  \xi(r,\mu) ]  & \geq \sqrt{r^2 -|b| + \tfrac{\omega^2}{8(r^2-|b|)}}.
\end{align}
To see that Eq.~\eqref{Eq:RepXi_large} holds we show 
\[ \sqrt{r^2 -|b| + \tfrac{\omega^2}{8(r^2-|b|)}} >  r - \frac{|b| r \omega^{-1}}{2(1+ r^2 \omega^{-1})}. \]
It suffices to show
\[ r^2 -|b| + \frac{\omega^2}{8(r^2-|b|)} >
\left  (r - \frac{|b| r \omega^{-1}}{2(1+ r^2 \omega^{-1})} \right)^2 = r^2 + \frac{|b|^2 r^2 - 4 |b|r^4 - 4 |b|r^2 \omega}{4  (r^2 +\omega)^2 }. \]
Note that 
\begin{align*}
\frac{\omega^2}{8(r^2-|b|)}  & -|b| -\frac{|b|^2 r^2 - 4 |b|r^4 - 4 |b|r^2 \omega}{4  (r^2 +\omega)^2 }   > 
\frac{\omega^2}{8(r^2+ \omega)} - \frac{|b|^2 r^2 + 4 |b|r^2 \omega + 4 |b| \omega^2}{4  (r^2 +\omega)^2 } \\
& =  \frac{\omega^3 - 8 |b| \omega^2 + r^2(\omega^2 - 8 |b| \omega -2 |b|^2)}{8(r^2 +\omega)^2}  > 0
\end{align*}
for $\omega \gg 1$. This implies the claim.
\end{proof}

By adding $Q(r,\mu) v(r)$ to both sides we rewrite Eq.~\eqref{Eq:HomODENorm} as
\begin{equation}\label{Eq:InHomLGPot}
v''(r)  - r^2 v(r) - \mu v(r) + Q(r,\mu) v(r)   =  \Oc(r^{-2} \omega^0) v(r),
\end{equation}
and apply perturbation theory to obtain a fundamental system for $r \geq 4$. 
The proof of the next result is partially along the lines of \cite{DonSch14a}, Proposition $4.3$.

\begin{lemma}\label{Le:FundSysWeber}
Define 
\begin{align*}
 v_0^{\pm}(r,\omega)  := \tfrac{1}{\sqrt{2}} \mu^{-\frac14} (1 + \tfrac{r^2}{\mu})^{-\frac14} e^{ \pm \mu \xi(r,\mu) }
\end{align*}
with Wronskian $W(v_0^{-}(\cdot,\omega) ,v_0^+(\cdot,\omega) ) = 1$. For $\omega \gg 1$, Eq.~\eqref{Eq:HomODENorm} has a fundamental system $\{ v_{-},v_+\}$ of the form
\begin{align}\label{Eq:FundSystWeber}
v_{\pm}(r,\omega) = v_0^{\pm}(r,\omega)[1 + \Oc(r^{-1} \omega^{-\frac{1}{2}}) ],
\end{align}
for all $r \geq 4$.
\end{lemma}

\begin{proof}
To construct the solution $v_-$, we use the fundamental system for the homogeneous version of Eq.~\eqref{Eq:InHomLGPot}
and the variation of constants formula to obtain an integral equation for $v_-$ given by 
\begin{align*} v_-(r,\omega)& =v_0^-(r,\omega) 
+v_0^+(r,\omega)
\int_r^\infty v_0^-(s,\omega) \Oc(s^{-2} \omega^0) v_-(s,\omega)ds \\
& -v_0^-(r,\omega)
\int_r^\infty v_0^+(s,\omega) \Oc(s^{-2} \omega^0) v_-(s,\omega)ds.
\end{align*}
First, we assume that $r \geq 3$. Setting $h_-:=\frac{v_-}{v_0^-}$ yields 
\begin{equation}\label{Eq:VoltWeberh}
h_-(r,\omega)=1+\int_r^\infty K(r,s,\omega)h_-(s,\omega)ds
\end{equation}
where
\begin{align*}
  K(r,s,\omega) & :=\left [\frac{v_0^+(r,\omega)}{v_0^-(r,\omega)}
v_0^-(s,\omega)^2-v_0^- v_0^+(s,\omega)\right ] \Oc(s^{-2} \omega^0)  = g(s,\omega)  \left(e^{- 2(\mu \xi(s,\mu)- \mu \xi(r,\mu))}  -1 \right),
\end{align*}
and $g(s,\omega) = \tfrac12 \mu^{-\frac12} (1 + \tfrac{s^2}{\mu})^{-\frac12}\Oc(s^{-2} \omega^0) = \Oc(s^{-2} \omega^{-\frac12})$. 
Eq.~\eqref{Eq:RepXi_large_alt} implies that
\begin{align}\label{Eq:WeberKernel}
 |  K(r,s,\omega)| \lesssim s^{-2} \omega^{-\frac12} \big [  1 + e^{-2 \mathrm{Re}\mu^{1/2}(s-r)}  e^{-2(\tilde \varphi(s,\omega) - \tilde \varphi(r,\omega))} \big ] 
 \lesssim  s^{-2} \omega^{-\frac12},
\end{align} 
since $s > r$ and $\tilde \varphi(\cdot, \omega)$ is monotonically increasing on $[3,\infty)$. This yields the bound 
\[ \int_{r}^\infty \sup_{ r \in [3,s]} |K(r,s,\omega)| ds \lesssim 1. \] 
Thus, we can apply the standard result on Volterra equations, see e.g. Lemma $2.4$ in \cite{SchSoffStaub2010}, which yields the existence of a solution $h_-$ 
to Eq.~\eqref{Eq:VoltWeberh} satisfying the bound
\[|h_-(r,\omega)-1|\lesssim r^{-1} \omega^{-\frac12},\]
for all $r \geq 3$. The bound 
\begin{align}\label{Eq:SymbBehWeberhp}
|\partial_r^{k} [h_-(r,\omega)-1]|\lesssim_{k}  r^{-1-k}\omega^{-\frac12} 
\end{align}
for all $k \in \N_0$, follows from Eq.~\eqref{Eq:VoltWeberh} by applying standard arguments, cf.~Remark $4.4$ in  \cite{DonSch14a}.
Thus, we can write $h_-(r,\omega) = 1 + \Oc(r^{-1} \omega^{-\frac12})$, 
which yields Eq.~\eqref{Eq:FundSystWeber}.

The second solution $v_+$ is obtained by setting 
\[ v_+(r,\omega):=v_-(r,\omega)\left [\frac{v_0^+(3,\omega)}
{v_0^-(3,\omega)}
+\int_{3}^r v_-(s,\omega)^{-2}ds \right ]. \]

Following the lines of \cite{DonSch14a}, p.~2499-2500, we use the identity 
\[ v_0^-(r,\omega)\int_3^r v_0^-(s,\omega)^{-2}ds=v_0^+(r,\omega)
-\frac{v_0^+(3,\omega)}{v_0^-(3,\omega)}v_0^-(r,\omega) \]
to obtain an expression for $h_+ :=\frac{v_+}{v_0^+}$ given by
\begin{align*}
h_+(r,\omega)=1&+\Oc(r^{-1}\omega^{-\frac12}) +\frac{v_-(r,\omega)}{v_0^+(r,\omega)}\int_3^r 
v_0^-(s,\omega)^{-2}\Oc(s^{-1}\omega^{-\frac12})ds.
\end{align*}
We show by induction that for all $k \in \N_0$,
\begin{align}\label{Eq:Weberh+}
\begin{split}
\partial^k_r [h_+(r,\omega) - 1] &  = \Oc(r^{-1-k} \omega^{-\frac12} ) 
+ e^{-2[\mu\xi(r,\mu) -\mu\xi(3,\mu)]} \Oc(r^0 \omega^{-\frac12 + \frac{k}{2}})\Oc( \langle \omega^{-\frac12}r \rangle^k)  \\
& \phantom{!!!!!!!!!!!!} + \int_3^{r} H_k(r,s,\omega) ds,
\end{split}
\end{align}
with
\begin{align}\label{Eq:Weberh+_Kernel}
H_k(r,s,\omega) =   e^{-2[\mu\xi(r,\mu) -\mu\xi(s,\mu)]} 
\Oc(r^0 s^{-2-k} \omega^{-\frac12}) \Oc( \langle \omega^{-\frac12}r \rangle^k) \Oc( \langle \omega^{-\frac12}s \rangle^{-k}).
\end{align}
For $k=0$, we use that $v_0^-(s,\omega)^{-2} = \partial_{s}e^{2\mu \xi(s,\mu)}$ and integrate by parts to get 
\begin{align*}
h_+(r,\omega)&=1+\Oc(r^{-1}\omega^{-\frac12}) +e^{-2\mu\xi(r,\mu)}[1+\Oc(r^{-1}\omega^{-\frac12})] \\
&\quad \times \Big [e^{2\mu\xi(r,\mu)}\Oc (r^{-1}\omega^{-\frac12})
- e^{2\mu\xi(3,\mu)} \Oc (\omega^{-\frac12})  +\int_3^r e^{2\mu\xi(s,\mu)} \Oc (s^{-2}\omega^{-\frac12})ds \Big ] \\
&=1+ \Oc (r^{-1}\omega^{-\frac12}) + e^{-2[\mu\xi(r,\mu) -\mu\xi(3,\mu)]} \Oc (r^0\omega^{-\frac12}) \\
&\quad  + \int_3^r e^{-2[\mu\xi(r,\mu) - \mu\xi(s,\mu)]} \Oc (r^0 s^{-2}\omega^{-\frac12})ds.
\end{align*}
Assume that Eq.~\eqref{Eq:Weberh+} and Eq.~\eqref{Eq:Weberh+_Kernel} hold up to some $k \in \N$. 
With $\mu \partial_r \xi(r,\mu) = \mu^{\frac12}(1 + \mu^{-1} r^2)^{\frac12}$, we get 
\begin{align*}
\partial^{k+1}_r  & [h_+(r,\omega) - 1]  \\
& = \Oc(r^{-2-k} \omega^{-\frac12} ) 
+ e^{-2[\mu\xi(r,\mu) -\mu\xi(3,\mu)]} \Oc(r^0 \omega^{\frac{k}{2}})\Oc( \langle \omega^{-\frac12}r \rangle^{k+1} )  + \int_3^{r} \partial_r H_k(r,s,\omega) ds,
\end{align*}
where 
\begin{align*}%\label{Eq:Weberh+_KernelDiff}
\partial_r H_k(r,s,\omega) =   e^{-2[\mu\xi(r,\mu) -\mu\xi(s,\mu)]}
\Oc( r^0 s^{-2-k} \omega^{0}) \Oc( \langle \omega^{-\frac12}r \rangle^{k+1}) \Oc( \langle \omega^{-\frac12}s \rangle^{-k}).
\end{align*} 
Note that  $[\mu \partial_s\xi(s,\mu)]^{-1} = \mu^{-\frac12} \Oc(\ \langle \omega^{-\frac12} s \rangle^{-1})$, see Eq.~\eqref{Eq:EstLGPot}.
Thus, an integration by parts yields 
\begin{align*}
& e^{-2\mu\xi(r,\mu)} \Oc(r^0 \omega^{0})\Oc(\langle \omega^{-\frac12}r \rangle^{k+1})  \int_3^{r} \left (\partial_se^{2\mu\xi(s,\mu)}\right )
\Oc(s^{-2-k}\omega^{-\frac12}) \Oc( \langle \omega^{-\frac12}s \rangle^{-k-1}) ds \\
& =\Oc( r^{-2-k} \omega^{-\frac12})  + e^{-2[\mu\xi(r,\mu) -\mu\xi(3,\mu)]} \Oc(r^0 \omega^{-\frac12} )\Oc( \langle \omega^{-\frac12}r \rangle^{k+1}) 
+ \int_3^{r} H_{k+1}(r,s,\omega) ds,
\end{align*}
where 
\[ H_{k+1}(r,s,\omega)  = e^{-2[\mu\xi(r,\mu) -\mu\xi(s,\mu)]}\Oc(r^0 \omega^{-\frac12}) \Oc(s^{-3-k}\omega^{0})
 \Oc( \langle \omega^{-\frac12}r \rangle^{k+1}) \Oc( \langle \omega^{-\frac12}s \rangle^{-(k+1)}).
\]
This implies the claim. For  $4 \leq r \leq \omega^{\frac12}$, $\omega \gg 1$, we use Lemma \ref{Le:RepXiQ} to estimate
\begin{align*}
\omega^{-\frac12 + \frac{k}{2}} e^{-2\mathrm{Re}[\mu\xi(r,\mu) -\mu\xi(3,\mu)]} &  \lesssim_k  \omega^{-\frac12 + \frac{k}{2}} 
e^{-2 \mathrm{Re}\mu^{1/2}(r - 3) }  e^{-2(\tilde \varphi(r,\omega)-\tilde \varphi(3,\omega))}  \\
& \lesssim_k r^{-1-k} \omega^{k} e^{-2 \mathrm{Re}\mu^{1/2}} \lesssim_k r^{-1-k}  \omega^{-\frac12},
\end{align*}
and for $r > \omega^{\frac12}$ we get, 
\[ r^k \omega^{-\frac12 } e^{-2\mathrm{Re}[\mu\xi(r,\mu) -\mu\xi(3,\mu)]} \lesssim_k \omega^{-\frac12 }  r^{4 + k} e^{-r^2} \lesssim  r^{-1-k}  \omega^{-\frac12 }. \]
The integral kernels satisfy 
 \[ |H_k(r,s,\omega)| \lesssim_k r^{-2-k} \omega^{-\frac12}, \]
for all $3 \leq s \leq r$,  $\omega\gg 1$ and $k \in \N_0$, since  
\begin{align*}
|H_k(r,s,\omega)| & \lesssim_k \omega^{-\frac12} s^{-2-k} \langle \omega^{-\frac12}r \rangle^k \langle \omega^{-\frac12}s \rangle^{-k} e^{-2\mathrm{Re}[\mu\xi(r,\mu) -\mu\xi(s,\mu)]} \\ 
& \lesssim_k \omega^{-\frac12}  s^{-2-k}  \langle \omega^{-\frac12}r \rangle^{k-b} 
\langle \omega^{-\frac12}s \rangle^{-k+b} e^{-(r^2 - s^2)}  e^{-2[\varphi(r,\omega)-\varphi(s,\omega)]} \\
&\lesssim_k\omega^{-\frac12}  s^{-2-k} r^{k+4} s^{-k-4} e^{-(r^2 - s^2)} \\
& \lesssim_k \omega^{-\frac12}   r^{-2-k} r^{2+k} s^{-2-k} r^{k+4} s^{-k-4} e^{-(r^2 - s^2)}  \\
& \lesssim_k  r^{-2-k} \omega^{-\frac12}e^{-(r^2 - 2(k+3) \log r)} e^{s^2 - 2 (k+3) \log s}  \lesssim_k  r^{-2-k} \omega^{-\frac12}.
\end{align*}
With these estimates, we infer that 
\begin{align*}
|\partial^{k}_r [h_+(r,\omega) - 1]| \lesssim_k r^{-1-k} \omega^{-\frac12} 
\end{align*}
for all $r \geq 4$, $\omega \gg 1$ and $k \in \N_0$. This justifies the notation $h_+(r,\omega) = 1 + \Oc(r^{-1} \omega^{-\frac12})$ and implies Eq.~\eqref{Eq:FundSystWeber} 
for $v_+$. 
\end{proof}

\subsubsection{A fundamental system near the origin}
For small radii, Eq.~\eqref{Eq:HomODENorm} is written as 
\begin{equation}\label{Eq:Bessel}
v''(r)  - \frac{4\nu^2-1}{4 r^2} v(r) - \mu v(r)   = \Oc(\langle r \rangle^2) v(r).
\end{equation}
The right hand side is again treated perturbatively. 

\begin{lemma}\label{Le:FundSystemCenterBessel}
Choose $c > 1$. Define
\[\psi_0(r,\omega) := \sqrt{r} J_{5/2} ( \I \mu^{\frac12} r), 
\quad \tilde \psi_1(r,\omega) := \sqrt{r}  Y_{5/2}( \I \mu^{\frac12} r).\]
Eq.~\eqref{Eq:HomODENorm}  has a fundamental system $\{v_0,v_1\}$ of the form
\begin{align}\label{Eq:FundBessel}
\begin{split}
v_0(r,\omega) & = \psi_0(r,\omega) [ 1+ \Oc(r^2 \omega^{0})], \\
v_1(r,\omega) &= [\tilde \psi_1(r,\omega) +  \Oc(\omega^0) \psi_0(r,\omega)] [ 1+ \Oc(r^0\omega^{-\frac12})],
\end{split}
\end{align}
for all $\omega \gg c^2$ and $r \in (0, c \omega^{-\frac12}]$. Furthermore, 
%\begin{align}\label{Eq:FundSystemCenterBessel_Derivatives}
%\begin{split}
%\partial_r v_0(r,\omega) &= \partial_r \psi_0(r,\omega)[ 1+ \Oc(r^0\omega^{-\frac12})], \\
%\partial_r v_1(r,\omega) &= [\partial_r \tilde \psi_1(r,\omega) +  \Oc(\omega^0) \partial_r \psi_0(r,\omega)][ 1+ \Oc(r^0\omega^{-\frac12})],
%\end{split}
%\end{align}and 
$v_0(r,\omega) = \Oc(r^3 \omega^{\frac54})$, $v_1(r,\omega) = \Oc(r^{-2} \omega^{-\frac54})$.
\end{lemma}

\begin{proof}\label{Re:ExplicitBessel}
We note that we have the explicit expressions, 
\begin{align*}
\sqrt{z} J_{5/2} ( \I z) & = \alpha_0 z^{-2} [ (3+ z^2) \sinh(z)  - 3 z \cosh(z)],\\
\sqrt{z} Y_{5/2} ( \I z) & = \alpha_1 z^{-2} [ (3+ z^2) \cosh(z)  - 3 z \sinh(z)],
\end{align*}
for some $\alpha_0, \alpha_1 \in \C$. The construction of $\{v_0, v_1\}$ is along the lines of \cite{DonSch14a}, Lemma $4.5$, for fixed $\nu = \frac{5}{2}$,
so we only sketch the main steps. 
The solution $v_0$ is obtained via Volterra iteration as in the proof of Lemma \ref{Le:FundSysWeber}. This yields 
\[v_0(r,\omega)  = \psi_0(r,\omega) [ 1+ \Oc(r^2 \omega^0)].\]

To construct $v_1$, the fundamental system $\{ \psi_0, \tilde \psi_1 \}$ for the homogeneous version of Eq.~\eqref{Eq:Bessel} 
is not well suited since $\tilde \psi_1$ might have zeros in the domain of investigation.
Therefore, we construct another fundamental system  $\{ \psi_0,  \psi_1 \}$.
For this, we use that the zeros of $J_{5/2}(z)$ 
are real, cf. \cite{Olv97}, p.~245, Theorem 6.2, which implies that $|\psi_{0}(r,\omega)| > 0$ for $0 < r \leq c \omega^{-\frac12}$ and $\omega \gg 1$.
We set  
\[ \psi_1(r,\omega):=-\psi_0(r,\omega)\int_r^{c\omega^{-1/2}}\psi_0(s,\omega)^{-2}ds, \]
such that $|\psi_1(r,\omega)| > 0$ for $r \in [0,c \omega^{-\frac12})$, 
$\psi_1(c \omega^{-\frac12},\omega) =0$.
Furthermore, we have the representation
\[ \psi_1(r,\omega)=\tfrac{\pi}{2}\sqrt{r}Y_{5/2}(\I \mu^\frac12 r)
+\Oc(\omega^0)\sqrt{r}J_{5/2}(\I \mu^\frac12 r), \]
by calculating the respective Wronskians, see \cite{DonSch14a}.
Using that $|v_0(r,\omega)|>0$ for all $r \in (0,c\omega^{-\frac12}]$ and $\omega \gg 1$,
$v_1$ can be constructed by setting
\[ v_1(r,\omega)=-\frac{2}{\pi}v_0(r,\omega)\int_r^{c\omega^{-1/2}}v_0(s,\omega)^{-2}ds, \]
and investigating the corresponding equation for $h_1 = \frac{v_1}{\psi_1}$ and $r \in [0,c\omega^{-\frac12})$.
This yields Eq.~\eqref{Eq:FundBessel}.
One can easily check that $\psi_0(r,\omega)  = \Oc(r^3 \omega^{\frac54})$ and $\tilde \psi_1(r,\omega)  = \Oc(r^{-2} \omega^{-\frac54})$,
which implies that $v_0(r,\omega) = \Oc(r^3 \omega^{\frac54})$ and $v_1(r,\omega) = \Oc(r^{-2} \omega^{-\frac54})$.
\end{proof}

Finally, we construct a fundamental system for intermediate values of $r$. 

\begin{lemma}\label{Le:FundSystemCenterHankel}
Choose $c > 1$ sufficiently large and define 
\[\psi_{\pm}(r,\omega) := \sqrt r H_{5/2}^{\mp}(\I \mu^{\frac12} r), \]
where $ H_{5/2}^{\mp}= J_{5/2} \mp \I Y_{5/2}$.
Provided $\omega \gg c^2$, there exists a fundamental system $\{\tilde v_{-} ,\tilde v_{+}\}$ of Eq.~\eqref{Eq:HomODENorm}
given by
\begin{align}\label{Eq:FundHankel}
 \tilde v_\pm(r,\omega)= \psi_{\pm}(r,\omega)  [1+ \Oc(r^0\omega^{-\frac12})],
\end{align}
for $r \in [\frac{1}{2} c\omega^{-\frac12},40]$.
Furthermore,  
$ \tilde v_\pm(r,\omega) = \Oc(r^0 \omega^{-\frac14}) e^{\pm \mu^{1/2} r }$.
\end{lemma}

\begin{proof}
As in the proof of Lemma $4.6$ in \cite{DonSch14a} we define $\tilde \psi_{\pm}(r,\omega) := \alpha_{\pm}  \psi_{\pm}(r,\omega)$, for $\alpha_{\pm} \in \C \setminus \{0\}$ such that $W(\tilde \psi_{-}(\cdot, \omega), \tilde \psi_{+}(\cdot, \omega) )=1$.
Note that $\tilde \psi_{\pm}$ can be given in closed form by
\begin{align*}
\tilde \psi_{\pm}(r,\omega) = \tfrac{1}{\sqrt{2}} \mu^{-\frac14} e^{\pm \mu^{1/2} r } ( 1 \mp  3 r^{-1} \mu^{-\frac12}  + 3 r^{-2}\mu^{-1}  ).
\end{align*}
We choose $c > 1$ large enough such that $| \tilde \psi_{\pm}(r,\omega)| > 0$ for $r  \geq \frac{1}{4}c\omega^{-\frac12}$.
To construct $\tilde v_{-}$, we use the variation of constants formula, see \cite{DonSch14a}, Lemma $4.6$, 
which yields an equation for  $h_{-} =\frac{\tilde v_{-}}{\tilde \psi_{-}}$ given by 
\begin{align}\label{Eq:VolterraItHanel}
h_-(r,\omega)=1+\int_r^{R} K(r,s,\omega)h_-(s,\omega)ds, 
\end{align}
for some $R > 0$ and 
\[ K(r,s,\omega)=\left [\tilde\psi_-(s,\omega) \tilde \psi_+(s,\omega)
-\frac{\tilde \psi_+(r,\omega)}{\tilde \psi_-(r,\omega)}\tilde \psi_-(s,\omega)^2 \right ]\Oc (\langle s\rangle^2). \]
We set $R = 50$. By using the above explicit formulas we get that 
\[ |K(r,s,\omega)| \lesssim \omega^{-\frac12} + \omega^{-\frac12} e^{- 2\mathrm{Re} \mu^{1/2}(s-r)} \lesssim \omega^{-\frac12},\]
for all $\frac{1}{4} c\omega^{-\frac12} \leq r \leq s \leq 50$. The standard result on Volterra equations yields the existence of a solution $h_{-}$
satisfying
\[  |h_-(r,\omega) -1| \lesssim  \omega^{-\frac12}, \]
for all $r \in [\frac{1}{4} c\omega^{-\frac12},50]$ and $\omega \gg 1$.  In the following, we restrict ourselves to 
$\frac{1}{4} c\omega^{-\frac12} \leq r \leq 40$. With this, we obtain
\begin{align}\label{Eq:SymbolHankem}
|\partial^k_r  h_-(r,\omega)| \lesssim_k r^{-k} \omega^{-\frac12}
\end{align}
for all $k \in \N$. 
For the first derivative, we differentiate Eq.~\eqref{Eq:VolterraItHanel} and integrate by parts to get 
\begin{align*}
\partial_r  h_-(r,\omega)  & =   \int_r^{50}e^{- 2\mu^{1/2}(s-r)} \Oc(r^0 \omega^0 s^0) h_-(s,\omega)ds = 
g_1(r,\omega) + \int_r^{50}H_1(r,s,\omega)  \partial_s h_-(s,\omega) ds,
\end{align*}
where 
\[ g_1(r,\omega) = \Oc(r^0 \omega^{-\frac12}) + e^{- 2 \mu^{1/2}(50-r)} \Oc(r^0 \omega^{-\frac12}) 
  + \int_r^{50}e^{- 2 \mu^{1/2}(s-r)}\Oc(r^0 \omega^{-\frac12} s^{-1}) h_-(s,\omega) ds,\]
and $H_1(r,s,\omega) = e^{- 2 \mu^{1/2}(s-r)}\Oc(r^0 \omega^{-\frac12} s^0)$. It is easy to see that  
$|g_1(r,\omega)| \lesssim r^{-1} \omega^{-\frac12}$
and also $|H_1(r,s,\omega)| \lesssim \omega^{-\frac12}$ such that 
the standard result on Volterra equations 
implies that 
\[ |\partial_r  h_-(r,\omega)| \lesssim r^{-1} \omega^{-\frac12}.\]
Furthermore, we have the bound $| \partial_r g_1(r,\omega)| \lesssim r^{-2} \omega^{-\frac12}$,
since $\omega^{\frac12} e^{- 2 \mathrm{Re} \mu^{1/2}(50-r)} \lesssim 1$ for $r \leq 40$. For higher derivatives
we get terms with a similar structure, which yields Eq.~\eqref{Eq:SymbolHankem} as well as Eq.~\eqref{Eq:FundHankel}
for $\tilde v_-$ (up to a constant). The second solution is obtained by setting 
\[ \tilde v_{+}(r,\omega) := \tilde v_{-}(r,\omega) \left [ \frac{ \tilde \psi_{+} (a ,\omega)}{ \tilde \psi_{-}(a,\omega)}
+\int_{a}^r \tilde v_{-}(s,\omega)^{-2}ds \right ],\]
where $a := \frac{1}{4} c\omega^{-\frac12}$. 
As in the proof of Lemma \ref{Le:FundSysWeber} this yields an equation for $h_{+} = \frac{\tilde v_+}{\tilde \psi_{+}}$
given by 
\begin{align*}
& h_{+}(r,\omega) =  1+ \Oc(r^0 \omega^{-\frac12}) + 
 \int_{a}^r e^{- 2 \mu^{1/2} (r - s)} \Oc(r^ 0 s^0 \omega^{0}) ds. 
\end{align*}
Integration by parts yields
\begin{align*}
& h_{+}(r,\omega) =  1+ \Oc(r^0 \omega^{-\frac12}) + e^{- 2 \mu^{1/2} (r - a)} \Oc(r^0 \omega^{-\frac12})  
+  \int_{a}^r e^{-2 \mu^{1/2} (r- s)} \Oc(s^{-1} \omega^{-\frac12}) ds.
\end{align*}

From this it is obvious that  $|h_{+}(r,\omega) - 1| \lesssim \omega^{-\frac12}$. 
Restricting ourselves to $\frac{1}{2} c\omega^{-\frac12} \leq r \leq 40$, we get that 
$|\partial^k_r h_{+}(r,\omega)| \lesssim_k r^{-k} \omega^{-\frac12}$ for all $k \in\N$ and $\omega \gg 1$
provided $c > 1$ is chosen sufficiently large. 
Dividing both $\tilde v_\pm(r,\omega)$ by $\alpha_{\pm}$ yields Eq.~\eqref{Eq:FundHankel}. The explicit formula above shows that 
$ \tilde v_\pm(r,\omega) = \Oc(r^0 \omega^{-\frac14}) e^{\pm \mu^{1/2} r }$.
\end{proof}

\subsubsection{A global fundamental system}

The proof of the following result is along the lines of Corollary $4.7$ and Lemma $4.8$ in \cite{DonSch14a}
(with the only difference that the Wronskians $W(v_{-}(\cdot,\omega), \tilde v_{\pm}(\cdot, \omega))$ 
are obtained by evaluation at $r = 10$ instead of $r = 1$). From now on we fix $c > 1$ such that 
Lemma \ref{Le:FundSystemCenterHankel} holds.

\begin{lemma}\label{Le:FundSystemCenterRepresentations}
Provided $\omega \gg c^2$ we have the representation
\begin{align*}
 \tilde v_\pm(r,\omega)&=\Oc(\omega^0) v_0(r,\omega) \mp [\I+\Oc (\omega^{-\frac12})]
v_1(r,\omega), \\
v_0(r,\omega)&=[\alpha_- + \Oc(\omega^{-\frac12})]\tilde v_-(r,\omega)
+[ \alpha_{+} +\Oc (\omega^{-\frac12})]\tilde v_+(r,\omega),
\end{align*}
for all $r \in (0,40]$. Under the same assumptions, 
\begin{align*}
 v_-(r,\omega)=    \alpha e^{10 \mu^{1/2}  }[1+\Oc(\omega^{-\frac12})]\tilde v_-(r,\omega)  + e^{- 10  \mu^{1/2}  }\Oc(\omega^{-\frac12}) \tilde v_+(r,\omega),
\end{align*}
for all $r \in [\frac{1}{2} c\omega^{-\frac12}, \infty)$, where $\alpha_\pm, \alpha \in \C\backslash \{0\}$.
\end{lemma}

With this we can construct a global fundamental system for Eq.~\eqref{Eq:HomODENorm}. 

\begin{lemma}\label{Le:Rep_vmv0}
The functions $\{ v_0, v_{-} \}$ provide a fundamental system for Eq.~\eqref{Eq:HomODENorm} and we have the representations
\begin{align}\label{Eq:Rep_v0vmvp}
v_0(r,\omega) = e^{10 \mu^{1/2}  } \Oc(\omega^{-\frac12})v_-(r,\omega) 
+  \alpha_1e^{10 \mu^{1/2}  }[1+\Oc(\omega^{-\frac12})]v_+(r,\omega),
\end{align}
\begin{align} \label{Eq:Rep_vmv0v1}
v_-(r,\omega)&= e^{10 \mu^{1/2}  } \Oc(\omega^0)v_0(r,\omega) 
+ \alpha_2 e^{10 \mu^{1/2}  } [1+\Oc(\omega^{-\frac12})]
v_1(r,\omega), 
\end{align}
for all $r>0$, all $\omega \gg c^2$ and some constants $\alpha_j \in \C\backslash\{0\}$, $j=1,2$.
The Wronskian is given by
\begin{align} \label{Eq:Wronskian}
W(\omega):&=W(v_-(\cdot,\omega),v_0(\cdot,\omega)) = \alpha_1 e^{10 \mu^{1/2}  }[1+\Oc(\omega^{-\frac12})].
\end{align}
\end{lemma}

\begin{proof} 
The fundamental systems $\{ v_-,v_+\}$ and $\{ \tilde v_-,\tilde v_+\}$ are valid in the
intervals $r \in [4, \infty)$ and $r \in [\frac12 c \omega^{-\frac12}, 40]$,
respectively. The connection formula yields 
\begin{align}\label{Eq:ConnFormtildevpm}
\tilde v_{\pm}(r,\omega) =   \frac{W_{\pm}(\cdot,\omega)}{W(v_{-}(\cdot,\omega), v_{+}(\cdot,\omega))} v_-(r,\omega)   + 
 \frac{\tilde W_{\pm}(\cdot,\omega) }{W(v_{+}(\cdot,\omega), v_{-}(\cdot,\omega))}   v_+(r,\omega),
\end{align}
where 
\[  W_{\pm}(\omega) := W(\tilde v_{\pm}(\cdot,\omega) , v_+(\cdot,\omega)), \quad \tilde W_{\pm}(\omega) := W(\tilde v_{\pm}(\cdot,\omega) , v_-(\cdot,\omega)). \]
By Lemma \ref{Le:FundSystemCenterRepresentations}, $v_0$ can be expressed in terms of Hankel functions. In combination with Eq.~\eqref{Eq:ConnFormtildevpm}
this yields
\begin{align}\label{Eq:v0RepDetail}
\begin{split}
v_0(r,\omega) & =  [\alpha_{-} + \Oc( \omega^{-\frac{1}{2}})] W_{-}(\omega)  v_-(r,\omega) +   [\alpha_{+} + \Oc( \omega^{-\frac{1}{2}})] W_{+}(\omega) v_-(r,\omega)  \\ 
& + 
[\alpha_{-} + \Oc( \omega^{-\frac{1}{2}})] \tilde W_{-}(\omega)  v_+(r,\omega)   
+ [\alpha_{+} + \Oc( \omega^{-\frac{1}{2}})]\tilde W_{+}(\omega)   v_+(r,\omega), 
\end{split}
\end{align}
for some constants $\alpha_{\pm} \in \C\setminus \{0\}$ and all $r > 0$.
Evaluation at $r = 10$ shows that 
\begin{align*}
 \tilde W_{+}(\omega) & = \beta_{+}  e^{10 \mu^{1/2}} [1 + \Oc( \omega^{-\frac{1}{2}})],\\
 \tilde W_{-}(\omega)  & = e^{- 10 \mu^{1/2}}\Oc( \omega^{-\frac{1}{2}}),\\
 W_{-}(\omega)  & = \beta_{-} e^{-10 \mu^{1/2}} [1 + \Oc( \omega^{-\frac{1}{2}})], \\
 W_{+}(\omega) & =   \Oc(\omega^{-\frac12}) e^{10 \mu^{1/2}},
\end{align*}
for some  $\beta_{\pm}  \in \C\backslash\{0\}$.
This implies Eq.~\eqref{Eq:Rep_v0vmvp}. 
Eq.~\eqref{Eq:Rep_vmv0v1} is a consequence of Lemma \ref{Le:FundSystemCenterRepresentations}. The expression for the 
Wronskian in Eq.~\eqref{Eq:Wronskian} can be verified easily.
\end{proof}

\subsection{Explicit representation of the resolvent for large imaginary parts}\label{SubSecRes:InhomEq}

Having a global fundamental system for Eq.~\eqref{Eq:HomODENorm}, we can consider the inhomogenous equation, which yields an explicit formula for the resolvent for $\omega \gg 1$. 

\begin{lemma}\label{Le:Resolvent_ExplicitForm}
Fix $\alpha > - \frac{1}{75}$. Set $b = 4 \alpha - 3$ and choose $M_{\alpha}$ sufficiently large.
Let $\lambda = \alpha + \I \omega$, $ \omega > M_{\alpha}$,  and for   
$\tilde g \in \mc H$, $\tilde g = g(|\cdot|)$, set $f :=  4 g(2 \cdot)$.
Then the resolvent $R_{L}(\lambda): \mc H \to \mc D(L) \subset \mc H$ 
exists and
\[ [R_{L}(\lambda) \tilde g](\xi) = [\mc R (4 \omega)f]\left (|\xi|/2 \right),\]
where 
\begin{align}\label{Eq:ResolventExp}
[\mc R(\omega)f](r) :=  \int_r^{\infty} g_1(r,s,\omega) f(s)ds + \int_0^r g_2(r,s,\omega) f(s)ds,
\end{align}
with 
\begin{align}\label{Eq:KernelsResolvent}
\begin{split}
  g_1(r,s,\omega)& =\frac{r^{-3} e^{\frac{1}{2}(r^2-s^2)} s^3}{W(\omega)}  v_0(r,\omega)v_-(s,\omega), \\
   g_2(r,s,\omega) & = \frac{r^{-3} e^{\frac{1}{2}(r^2-s^2)} s^3}{W(\omega)}v_-(r,\omega)v_0(s,\omega).
\end{split}
\end{align}
\end{lemma}

\begin{proof}
First, we use the above results and the variation of constants formula to infer that solutions to the  equation
\begin{equation}\label{Eq:Res_ODEInhom_Normal}
v''(r)  - r^2 v(r) - \frac{4\nu^2-1}{4 r^2} v(r)   +  \tilde V(r) v(r) - (b + \I \tilde \omega) v(r) = - r^3 e^{-\frac{r^2}{2}} f(r)
\end{equation} 
are of the general form
\begin{align}\label{EQ:SolInhom}
\begin{split}
v(r) = c_0 v_0(r, \tilde  \omega) &  + c_{-} v_{-}(r,\tilde  \omega) -  v_0(r,\tilde \omega)   \int_{r_0}^{r} \frac{v_{-}(s,\tilde  \omega)}{W(\tilde \omega)} s^3 e^{-\frac{s^2}{2}} f(s) ds  \\
 & +   v_{-}(r,\tilde  \omega)  \int_{r_1}^{r} \frac{v_{0}(s,\tilde  \omega)}{W(\tilde \omega)} s^3 e^{-\frac{s^2}{2}} f(s) ds
 \end{split} 
\end{align}
for constants $c_0, c_{-} \in \C$, $r_0, r_1 \in \R$, provided $b > - 4$ and $\tilde  \omega \geq M_b$ for some $M_b > 0$ large enough.

For fixed $\alpha > - \frac{1}{75}$, we define $b := 4 \alpha - 3$, $M_{\alpha} := \frac14 M_b$. 
For $\lambda = \alpha + \I  \omega\in \Omega_{\alpha}$, we set $\tilde \omega = 4 \omega$. Since $\lambda \in \rho(L)$,
$(\lambda - L) R_{L}(\lambda) \tilde g = \tilde g$. With $\tilde w_{\lambda} := R_L(\lambda) \tilde  g$, $\tilde w_{\lambda} = w_{\lambda}(|\cdot|)$, we infer that 
$w_{\lambda}$ satisfies the equation
\begin{equation*}
\lambda w_{\lambda}(\rho) - w_{\lambda}''(\rho) - \frac{6}{\rho} w_{\lambda}'(\rho) + \frac{1}{2} \rho w_{\lambda}'(\rho) + w_{\lambda}(\rho) - V(\rho) w_{\lambda}(\rho)  = g(\rho).
\end{equation*}
Thus,  $v_{\lambda}(r) := r^3 e^{-\frac{1}{2} r^2} w_{\lambda}(2 r)$ solves Eq.~\eqref{Eq:Res_ODEInhom_Normal} for $f := 4 g(2 \cdot)$.
By definition, $b > -4$ and $\tilde \omega \geq M_b$ such that $v_{\lambda}$ is of the form stated in Eq.~\eqref{EQ:SolInhom}.
The fact that $R_{L}(\lambda) \tilde g  \in \mc H$ implies that $w_{\lambda}$ is continuous and decays at infinity according to
Lemma \ref{Le:PropH}. The resulting conditions on $v_{\lambda}$ can only be satisfied if 
\begin{align*}
c_0 =  \int_{r_0}^{\infty} \frac{v_{-}(s,\tilde \omega)}{W(\tilde \omega)} s^3 e^{-\frac{s^2}{2}} f(s) ds, 
\quad c_{-} = \int_{0}^{r_1} \frac{v_{0}(s,\tilde \omega)}{W(\tilde \omega)} s^3 e^{-\frac{s^2}{2}} f(s) ds.
\end{align*}
This shows that 
\begin{align}\label{Eq:Res_ODEInhom_Normal_Sol}
v_{\lambda}(r) =\frac{v_0(r,\tilde \omega)}{W(\tilde \omega)} \int_r^\infty v_-(s,\tilde \omega)
s^{3}e^{-\frac{1}{2}s^2} f(s)ds  +\frac{v_-(r,\tilde \omega)}{W(\tilde \omega)}\int_0^r v_0(s,\tilde \omega)
s^{3} e^{-\frac{1}{2}s^2}  f(s)ds,
\end{align}
which yields Eq.~\eqref{Eq:ResolventExp} and Eq.~\eqref{Eq:KernelsResolvent}.
\end{proof}

\subsection{Uniform bounds for large imaginary parts} 
In the following, we assume again that $b > -4$ is fixed and that $\omega \gg 1$.
It is sufficient to prove Eq.~\eqref{Eq:ResolventBounds} on a dense subset, hence
we restrict ourselves to $f \in C^{\infty}_{\mathrm{e},0}(\R)$.
Using Eq.~\eqref{Eq:KernelsResolvent} we define an integral operator $\mc T(\omega)$ by 
\[ [\mc T(\omega)f](r) := \int_r^{\infty} \partial_r g_1(r,s,\omega) f(s) ds + \int_0^{r} \partial_r g_2(r,s,\omega) f(s) ds\]
for $r > 0$. The following result is crucial.

\begin{lemma}\label{Le:ResEst_T}
Define $\delta \in [0,\frac12)$ by
\begin{align}\label{Def:delta}
\delta := \begin{cases} 
-\frac{3}{2} - \frac{b}{2}  &  \quad \mathrm{ for }~ -4 < b  \leq -3,\\
0  & \quad \mathrm{ for }~ b >  -3. 
 \end{cases} 
\end{align}
For $m=0,1$, 
\begin{align*}
r^{m+2}[\mc  T(\omega) f]^{(m)}(r)   &  =  O(r\omega^0)f(r)  +   O(r^2\omega^0) f'(r) \\
& + \sum_{k_1=0}^{2} \big [\mc  J^m_{1k_1}(\omega)(\cdot)^{k_1+1} f^{(k_1)} \big ](r)  
+   \sum_{k_2=0}^{3} \big [\mc  J^m_{2k_2}(\omega)(\cdot)^{k_2} f^{(k_2)} \big ](r) , 
\end{align*}
where 
\[ [\mc J^m_{ik_i}(\omega)f](r) = \int_0^{\infty}  J^m_{ik_i}(r,s,\omega) f(s) ds,\]
for $i = 1,2$. The integral kernels satisfy
\[ |J^m_{ik_i}(r,s,\omega) | \lesssim \min \{ r^{-1 + \delta} s^{-\delta},  r^{- \delta} s^{-1 + \delta} \},\]
for all $r,s > 0$ and $\omega \gg 1$. For $n = 0,\dots, 3$, 
\begin{align*}
& r^n [T(\omega)f]^{(n)}(r)  =  \sum_{j=0}^{2} O(r^{j+1}\omega^0) f^{(j)}(r) + \sum_{k=0}^{3} O(r^{k}\omega^0) f^{(k)}(r) 
 \\ & +  \sum_{j_1=0}^{2} \big [\mc K^n_{1j_1}(\omega)(\cdot)^{j_1+1} f^{(j_1)} \big ](r) 
 +  \sum_{j_2=0}^{3} \big [\mc K^n_{2j_2}(\omega)(\cdot)^{j_2} f^{(j_2)} \big](r) 
+ \sum_{j_3=1}^{4} \big [\mc K^n_{3j_3}(\omega)(\cdot)^{j_3-1} f^{(j_3)} \big ](r) ,
\end{align*}
where 
\[ [\mc K^n_{i j_i}(\omega)f](r) = \int_0^{\infty}  K^n_{i j_i}(r,s,\omega) f(s) ds\]
for $i = 1,2,3$, and 
\[ |K^n_{i j_i}(r,s,\omega) | \lesssim \min \{ r^{-1 + \delta} s^{-\delta},  r^{- \delta} s^{-1 + \delta} \},\]
for all $r,s > 0$ and $\omega \gg 1$.
\end{lemma}

\begin{proof}
The proof is based on the explicit representations of the kernel functions $g_1, g_2$ in different regimes.
We only sketch the argument here and refer the interested reader to Appendix \ref{App_ProofTLemma} for full details.

Consider the operator
$\int_0^r \partial_r g_2(r,s,\omega)f(s)ds$	, 
i.e., the case $s\leq r$.
From the representation \eqref{Eq:Rep_v0vmvp} of $v_0$ in terms of Weber functions we see that the most significant contribution to the kernel $g_2$ comes from the term
\begin{align*} \frac{1}{W(\omega)}&r^{-3}s^3 e^{\frac12(r^2-s^2)}v_-(r,\omega)\Oc(\omega^0)e^{10\mu^{1/2}}v_+(s,\omega) \\
&=e^{\frac12 r^2-\mu\xi(r,\mu)}e^{-[\frac12 s^2-\mu\xi(s,\mu)]}
\Oc(r^{-3}\langle\omega^{-\frac12}r\rangle^{-\frac12}s^3
\langle\omega^{-\frac12}s\rangle^{-\frac12}\omega^{-\frac12}),
\end{align*}
 see Lemma \ref{Le:FundSysWeber}, where we restrict ourselves to the regime $10\leq s\leq r$.
 With $\mu\partial_r \xi(r,\mu)=\mu^\frac12 \sqrt{1+\frac{r^2}{\mu}}$ one checks that 
 \begin{equation}
 \label{eq:dife} \partial_r e^{\frac12r^2-\mu\xi(r,\mu)}=\Oc(\langle \omega^{-\frac12}r\rangle^{-1}
 \omega^\frac12)e^{\frac12r^2-\mu\xi(r,\mu)} 
 \end{equation}
 and also
 \begin{equation}
 \label{eq:inte}
  e^{\frac12r^2-\mu\xi(r,\mu)}=\Oc(\langle \omega^{-\frac12}r\rangle
 \omega^{-\frac12})\partial_r e^{\frac12r^2-\mu\xi(r,\mu)}.
 \end{equation}
 Consequently, the dominant contribution to $\partial_r g_2(r,s,\omega)$ is of the form
 \[ e^{\frac12 r^2-\mu\xi(r,\mu)}e^{-[\frac12 s^2-\mu\xi(s,\mu)]}
\Oc(r^{-3}\langle\omega^{-\frac12}r\rangle^{-\frac32}s^3
\langle\omega^{-\frac12}s\rangle^{-\frac12}\omega^{0}) \]
and it suffices to consider the operator
\[ [T_W(\omega)f](r):=\int_0^r \chi(s) 
e^{\frac12 r^2-\mu\xi(r,\mu)}e^{-[\frac12 s^2-\mu\xi(s,\mu)]}
r^{-3}\langle\omega^{-\frac12}r\rangle^{-\frac32}s^3
\langle\omega^{-\frac12}s\rangle^{-\frac12}f(s)ds, \]
where $\chi$ is a smooth cut-off that localizes to $[10,\infty)$.
By Eq.~\eqref{eq:inte}, an integration by parts yields
\begin{align*} [T_W(\omega)f](r)=&\Oc(r^{-1}\omega^0)f(r) \\
&+\int_{0}^r \chi(s) 
e^{\frac12 r^2-\mu\xi(r,\mu)}e^{-[\frac12 s^2-\mu\xi(s,\mu)]}
\Oc(r^{-3}\langle\omega^{-\frac12}r\rangle^{-\frac32}s^0
\langle\omega^{-\frac12}s\rangle^{\frac12}\omega^{-\frac12})\partial_s[s^3f(s)]ds 
\end{align*}
where we have used $\Oc(\langle\omega^{-\frac12}r\rangle^{-1}\omega^{-\frac12})=\Oc(r^{-1}\omega^0)$ for the boundary term.
From Eqs.~\eqref{eq:dife} and \eqref{eq:inte} it follows that one may trade derivatives in $r$ for derivatives in $s$ at the expense of additional weights. 
More precisely, repeated integration by parts yields
\begin{align*}
\partial_r^n &[T_W(\omega)f](r)=\sum_{k=0}^n \Oc(r^{-1-n+k}\omega^0)f^{(k)}(r) \\
&+\int_{0}^r \chi(s)
e^{\frac12 r^2-\mu\xi(r,\mu)}e^{-[\frac12 s^2-\mu\xi(s,\mu)]}
\Oc(r^{-3}\langle\omega^{-\frac12}r\rangle^{-\frac32-n}s^0
\langle\omega^{-\frac12}s\rangle^{\frac12+n}\omega^{-\frac12})\partial_s^{1+n}[s^3f(s)]ds  
\end{align*}
for any $n\in \N_0$. From Eq.~\eqref{Eq:RepXi_large} we infer the bound
\[ \left |  e^{\frac12 r^2-\mu\xi(r,\mu)}e^{-[\frac12 s^2-\mu\xi(s,\mu)]} \right |
\lesssim \langle \omega^{-\frac12}r\rangle^{-\frac{b}{2}}\langle \omega^{-\frac12}s\rangle^{\frac{b}{2}}, \]
and it is straightforward to prove the stated estimate for the integral kernel.

In the case $r\leq s$, the most singular contribution to the kernel $g_1$ comes from the regime $0<r\leq s\lesssim \omega^{-\frac12}$ and the term 
\[ \frac{1}{W(\omega)}r^{-3}s^3 e^{\frac12(r^2-s^2)}v_0(r,\omega)\Oc(\omega^0)e^{10\mu^{1/2}}v_1(s,\omega)
=\Oc(r^{-3}s^3\omega^0)v_0(r,\omega)v_1(s,\omega), \]
see Eq.~\eqref{Eq:Rep_vmv0v1}.
From Lemma \ref{Le:FundSystemCenterBessel} we infer $\partial_r [r^{-3}v_0(r,\omega)]=\Oc(r\omega^\frac94)$
and this shows that the most important contribution to $\partial_r g_1(r,s,\omega)$ is of the form $\Oc(rs\omega)$. With this it is 
straightforward to prove the stated bounds.

The other cases are in some sense interpolates which can be treated similarly.
\end{proof}

\begin{lemma}\label{Le:EstT}
For $m = 0,1$, $n=0,\dots,3$, all $f \in C^{\infty}_{\mathrm{e},0}(\R)$ and all $\omega \gg 1$, we have the bounds 
\begin{align*}
\|(\cdot)^{m+2}  [\mc T(\omega)f]^{(m)} \|_{L^2(\R^+)} \lesssim \| f(|\cdot|) \| 
\quad \text{ and }  \quad \|(\cdot)^n  [\mc T(\omega)f]^{(n)} \|_{L^2(\R^+)} \lesssim \| f(|\cdot|)  \|.
\end{align*}
\end{lemma}

\begin{proof}
Choose $C = C_b >0$ sufficiently large such that the results of Lemma \ref{Le:ResEst_T} hold for all $\omega \geq C$. 
The integral operators $\mc J^m_{ik_i}(\omega)$ and $\mc K^n_{i j_i}(\omega)$ of Lemma \ref{Le:ResEst_T}
extend to bounded operators on $L^2(\R^+)$ by Lemma 5.5 in \cite{DonKri13}. Since all bounds are uniform in $\omega$, we infer
\[ \|\mc J^m_{ik_i}(\omega) f \|_{L^2(\R^+)}  \lesssim \| f \|_{L^2(\R^+)}, 
\quad \|\mc K^n_{i j_i}(\omega)f \|_{L^2(\R^+)}  \lesssim \| f \|_{L^2(\R^+)}, \]
for all $f \in  C^{\infty}_{\mathrm{e},0}(\R)$.
This yields 
\begin{align*}
 \|(\cdot)^{m+2}  [\mc T(\omega)f]^{(m)} \|_{L^2(\R^+)} & \lesssim \| (\cdot) f\|_{L^2(\R^+)}  + \| (\cdot)^2 f '\|_{L^2(\R^+)} \\
& +  \sum_{k_1=0}^{2} \|\mc  J^m_{1k_1}(\omega)(\cdot)^{k_1+1} f^{(k_1)} \|_{L^2(\R^+)}  
 +    \sum_{k_2=0}^{3} \| \mc  J^m_{2k_2}(\omega)(\cdot)^{k_2} f^{(k_2)} \|_{L^2(\R^+)}  \\
& \lesssim 
\sum_{k_1=0}^{2} \|(\cdot)^{k_1+1} f^{(k_1)} \|_{L^2(\R^+)} +  \sum_{k_2=0}^{3} \|(\cdot)^{k_2} f^{(k_2)} \|_{L^2(\R^+)}
 \lesssim \| f(|\cdot|) \|.
\end{align*}
for $m = 0,1$, where the last step follows from Lemma \ref{Le:AllBounds}. The second estimate can be derived analogously.
\end{proof}

\subsubsection{Proof of Proposition \ref{Prop:ResolventBounds}}
Fix $\alpha > - \frac{1}{75}$. Set $b = 4 \alpha - 3$ and choose $C_{\alpha} > 0$ sufficiently
large such that the above results hold for all $\omega \geq C_{\alpha}$. Set $M_{\alpha}:= \frac{1}{4} C_{\alpha}$.
Let $\lambda = \alpha + \I \omega$, $\omega \geq M_{\alpha}$. For $\tilde g \in C^{\infty}_{\mathrm{rad},0}(\R^7) $, 
$\tilde g = g(|\cdot|)$, we have that $g \in C^{\infty}_{\mathrm{e},0}(\R)$ and we set $f := 4 g(2 \cdot)$.
Furthermore, we define $\tilde \omega := 4 \omega$.
According to Lemma \ref{Le:Resolvent_ExplicitForm} we
have an explicit expression for the resolvent. By definition, rescaling and 
the fact that $[\mc R( \tilde\omega) f]'  = \mc T(\tilde \omega) f$, we can apply Lemma \ref{Le:EstT} to infer that 
\begin{align*}
 \| & R_L(\lambda ) \tilde g \|  \lesssim  
\|(\cdot)^2 [\mc R(\tilde \omega) f]'  \|_{L^2(\R^+)} + \|(\cdot)^3 [\mc R(\tilde \omega) f]''  \|_{L^2(\R^+)} 
 +  \sum_{k=1}^{4} \| (\cdot)^{k-1} [\mc R(\tilde \omega) f]^{(k)} \|_{L^2(\R^+)} \\
 & \lesssim   \|(\cdot)^{2} [\mc T(\tilde \omega) f] \|_{L^2(\R^+)} + \|(\cdot)^{3} [\mc T(\tilde  \omega) f]' \|_{L^2(\R^+)} + 
 \sum_{k=0}^{3} \| (\cdot)^{k} [\mc T(\tilde \omega) f]^{(k)}  \|_{L^2(\R^+)} 
  \lesssim  \|\tilde  g \|.
\end{align*}
The density of $C^{\infty}_{\mathrm{rad},0}(\R^7)$ in $\mc H$ implies that Eq.~\eqref{Eq:ResolventBounds}
holds for all $\tilde f\in \mc H$. For negative imaginary parts we write 
$\lambda = \alpha - \I \omega$, $ \omega \geq M_{\alpha}$. Since $L$ has real coefficients, the equation
$(\lambda -L) R_L(\lambda)\tilde g = \tilde g$ yields
\[ (\overline \lambda - L) \overline{ R_L(\lambda)\tilde g} = \overline{\tilde g}  \]
by complex conjugation. Hence, $\overline{ R_L(\lambda)\tilde g} = R(\overline \lambda) \overline{\tilde g}$.
By applying the above result we infer that Proposition \ref{Prop:ResolventBounds} holds.

\section{Growth estimates for $S(\tau)$}

Lemma \ref{Le:Spec} and the uniform bounds on the resolvent now allow us to 
derive growth bounds for the linearized time evolution. 

\begin{lemma}\label{Le:LinearTimeEvol}
There exists a projection $P \in \mc B(\mc H)$ with $\mathrm{rg} P = \mathrm{span}(\mb{\tilde g})$,
such that $P$ commutes with $S(\tau)$ for all $\tau \geq 0$ and
\[ \|P S(\tau) \tilde u \| = e^{\tau} \|P \tilde u  \| \]
for all $\tilde u  \in \mc H$. Moreover,
\begin{align}\label{Eq:StableEst}
 \|(1 - P)  S(\tau) \tilde u  \|  \leq C e^{-a \tau}  \|(1 -  P) \tilde u   \| 
\end{align}
for $a =  \frac{1}{150}$, all $ \tilde u  \in \mc H$, all $\tau \geq  0$ and some constant $C \geq 1$.
\end{lemma}

\begin{proof}
By Lemma \ref{Le:Spec}, the eigenvalue $\lambda = 1$ is isolated in the spectrum of $L$. Hence, we can define a spectral projection $P \in \mc B(\mc H)$ by
\begin{align*}
P=\frac{1}{2\pi \I}\int_\gamma R_{L}(\lambda) d\lambda,
\end{align*}
where $\gamma$ is a positively oriented circle around $1$ in the complex plane with radius $r_{\gamma} = \frac12$, cf.~\cite{kato}, p.~178. 
Note that $[ P, S(\tau) ] = 0$ for all $\tau \geq 0$.
Furthermore, $\mc H = \ker  P \oplus  \rg  P$ and $L$ is decomposed into the parts $L_\mc M$ and $L_\mc N$ in $\mc M:=\rg P$ and $\mc N:=\ker P$, respectively.
The respective spectra are given by
\begin{align*}
\sigma(L_\mc M) &= \{1\} \\
\sigma(L_\mc N) &= \sigma(L) \setminus \{1\}\subseteq \{z\in \C: \Re z\leq -\tfrac{1}{75}\},
\end{align*}
see Lemma \ref{Le:Spec}.
One always has that $\ker (1 - L) \subseteq \rg P$, see for example \cite{Sigal_Hislop}. We show that in our case also the reverse 
inclusion holds. First, we observe that $P$ has finite rank. This is a consequence of the invariance of the essential spectrum  
under relative compact perturbations, \cite{kato}, p.~239, Theorem 5.28, 
and the fact that $1 \not \in \sigma(L_0)$. We infer that the operator
$1 - L_\mc M: \mc M\to\mc M$ is finite-dimensional with zero as its only spectral point. Consequently, $1-L_\mc M$ is nilpotent and 
there exists a minimal $k \geq 1$ such that
$(1 - L_\mc M)^{k} \tilde u = 0$ for all $\tilde u  \in \mc M=\mathrm{rg} P$. If $k=1$, the desired inclusion $\rg P \subseteq \ker(1-L)$ follows immediately. 
Assume that $k \geq 2$. Then there exists a nontrivial
$\tilde v \in \mathrm{rg}(1 - L_\mc M) \cap \mathrm{ker} (1 - L_\mc M)$ and thus,  
a $\tilde u = u(|\cdot|) \in \mc M\subset \mc D(L)$ satisfying $(1 - L_\mc M) \tilde u = (1 - L)\tilde u = c \mb{\tilde g}$ for
some $c \in \C\setminus \{0\}$, see Lemma \ref{Le:Ker1mL}. Without loss of generality, we set $c = -1$. By introducing polar coordinates we see that $u$ solves
the  ordinary differential equation
\begin{align}\label{Eq:AlgMultODE}
u''(\rho) + \tfrac{6}{\rho} u'(\rho) - \tfrac12 \rho u'(\rho) + V(\rho) u(\rho) - 2 u(\rho) = \mb g(\rho).
\end{align}
Since $\tilde u \in \mc H$, we have $u \in C[0,\infty)$ and $\lim_{\rho \to \infty} \rho^{\frac32} |u(\rho)| = 0$, see Lemma \ref{Le:PropH}. 
We will see that there is no solution of Eq.~\eqref{Eq:AlgMultODE} having these properties.

For the homogeneous version of Eq.~\eqref{Eq:AlgMultODE} we have the fundamental system $\{ \mb g, \mb h \}$, see the proof
of Lemma \ref{Le:Ker1mL}, where $\mb h$ is given by Eq.~\eqref{Eq:SecSolSymmEigenvalEq}.
By the variation of constants formula, the general solution to Eq.~\eqref{Eq:AlgMultODE} 
is of the form
\begin{align*}
u(\rho) = c_0 \mb g(\rho) + c_1 \mb h(\rho) 
- \mb g(\rho) \int_{\rho_0}^{\rho}\mb h(s) \mb g(s)s^6 e^{-\frac{s^2}{4}}  ds 
+ \mb h(\rho) \int_{\rho_1}^{\rho} \mb g(s)^2 s^6 e^{-\frac{s^2}{4}}  ds. 
\end{align*}
Recall the representation $\mb h(\rho)=\rho^{-5}\langle\rho\rangle^2 e^{\frac{\rho^2}{4}}H(\rho)$ from the proof of Lemma \ref{Le:Ker1mL}.
To guarantee continuity of $u$ at zero, we are forced to choose $c_1 = \int_{0}^{\rho_1} \mb g(s)^2 s^6 e^{-\frac{s^2}{4}}  ds$
such that 
\begin{align*}
u(\rho) = c_0 \mb g(\rho)  - \mb g(\rho)  \int_{\rho_0}^{\rho}   s \langle s \rangle^{2} H(s) \mb g(s)  ds
+ \rho^{-5} \langle \rho \rangle^2 e^{\frac{\rho^2}{4}} H(\rho) \int_{0}^{\rho}  \mb g(s)^2  s^6 e^{-\frac{s^2}{4}}  ds.
\end{align*}
To obtain decay at infinity we must have
\[ \lim_{\rho \to \infty} \int_0^{\rho}  \mb g(s)^2 s^6 e^{-\frac{s^2}{4}}  ds = 0.\]
This, however, is impossible since the integrand is strictly positive on $\R^+$. 
By contradiction, we infer that $k=1$ and $\mc M=\mathrm{rg} P =\ker(1-L_\mc M)
\subseteq \ker(1-L)=\mathrm{span}(\mb{\tilde g})$.

It remains to prove the growth estimates on the stable and unstable subspaces. The fact that $\mb {\tilde g}$ is an eigenfunction of $L$ with eigenvalue $\lambda = 1$ yields 
$\|P S(\tau) \tilde u \| = \| S(\tau) P \tilde u \| = e^{\tau} \|P \tilde u \|$
for all $\tilde u \in \mc H$. 
To see that $(L_{\mc N}, \mathcal D(L)\cap \mc N)$ generates a $C_0$-semigroup $\{S_{\mc N}(\tau): \tau \geq 0 \}$ 
with $S(\tau)|_{\mc N} = S_{\mc N}(\tau)$, one argues for example as in \cite{DonSchAic12}, Lemma 4.17 and Corollary 4.18.
It is well-known that for all $\tau > 0$,
\[ r(S_{\mc N}(\tau)) = e^{\omega_{\mc N} \tau},\]
where $\omega_{\mc N}$ denotes the growth bound of the semigroup on $\mc N$ and $r(S_{\mc N}(\tau))$ is the spectral
radius of the bounded operator $S_{\mc N}(\tau): \mc N \to \mc N$, see \cite{engel}, p.~251.
To obtain Eq.~\eqref{Eq:StableEst}, it suffices to show that for each $\tau > 0$, 
\[ \Lambda_{\tau} := \{ z \in \C: |z| > e^{-\frac{1}{75}\tau} \} \subseteq   \rho(S_{\mc N}(\tau)). \]
Let $z \in \Lambda_{\tau}$ for some fixed $\tau > 0$ and assume that $z = e^{\lambda \tau}$ for some $\lambda \in \C$. 
Then $\mathrm{Re}\lambda = \frac{1}{\tau} \log |z| > -\frac{1}{75}$ and therefore, $\lambda \in \rho(L_{\mc N})$. Hence,
\[ \{ \lambda \in \C: z = e^{\lambda \tau} \} \subseteq \rho(L_{\mc N}). \]
By Proposition \ref{Prop:ResolventBounds} and the fact that $R_{L_{\mc N}}(\lambda) \tilde f = R_{L}(\lambda)|_{\mc N} \tilde f$ for $\tilde f \in \mc N$,
we infer 
\[ \| R_{L_{\mc N}}(\tfrac{1}{\tau} \log |z| + \tfrac{\I}{\tau} \arg z +\tfrac{2 \pi \I k }{\tau}  ) \tilde f \| \leq C  \|\tilde f\| \]
for all $\tilde f\in \mc N$, all $k\in \Z$, and an absolute constant $C  >0 $. In particular,
\[  \sup \{ \|R_{L_{\mc N}}(\lambda) \| : z = e^{\lambda \tau} \} < \infty \]
and by \cite{Pru84}, Theorem 3, $z \in \rho(S_{\mc N}(\tau))$. This implies $\omega_{\mc N} \leq - \frac{1}{75}$. 
Since the growth bound is defined as an infimum (that may not be attained),
we see that for each $\varepsilon > 0$ there exists a constant $C_{\varepsilon} \geq 1$ such that 
\[ \|(1-P)S(\tau) \tilde u \| = \|S_{\mc N}(\tau) (1-P) \tilde u \| \leq C_{\varepsilon} e^{-(\frac{1}{75} - \varepsilon) \tau} \|(1-P)\tilde u \|. \]
Choosing $\varepsilon = \frac{1}{150}$ implies the claim. 
\end{proof}

\section{Nonlinear Stability}\label{Sec:NonlinearStability}

In the sequel, we denote by $\mc B \subset \mc H$ the closed unit ball in $\mc H$.
For $\tilde u \in C^{\infty}_{\mathrm{rad},0}(\R^7)$ we define
\begin{align*}\label{Eq:DefNonlinearity}
  N(\tilde u) :=   f_1(|\cdot|) \tilde u^2  + f_2(|\cdot|)  \tilde u^3 , 
\end{align*}
where 
\[f_1(\rho) :=  - 9 ( 1+ \rho^2 \mb W(\rho)), \quad f_2(\rho) := -3 \rho^2. \]

\subsection{Estimates for the nonlinearity}

We need the following auxiliary result. 

\begin{lemma}\label{Le:NonlAux}
We have 
\begin{align}\label{Eq:AuxNonl_1}
\| f_1(|\cdot|)   \tilde u \tilde v \| \lesssim \| \tilde u \|  \| \tilde v \|,
\end{align}
and 
\begin{align}\label{Eq:AuxNonl_2}
 \| f_2(|\cdot|)  \tilde u \tilde v \tilde w \| \lesssim \| \tilde u \|  \| \tilde v \| \| \tilde w \|,
\end{align}
for all $\tilde u, \tilde v, \tilde w \in  C^{\infty}_{\mathrm{rad},0}(\R^7)$.
\end{lemma}

\begin{proof}
In the following, we make extensive use of Lemma \ref{Le:AllBounds} and as always, we write $\tilde u=u(|\cdot|)$ and analogously for $\tilde v$ and $\tilde w$.
To improve readability, we do not distinguish between $\Lap$ and $\Lap_{\mathrm{rad}}$.  Furthermore, we abbreviate $L^\infty:=L^{\infty}(\R^+)$ and $L^2:=L^{2}(\R^+)$.
We start with Eq.~\eqref{Eq:AuxNonl_1}.
First, observe that 
\begin{align*}
\begin{split}
|f_1^{(2k)}(\rho)| \lesssim_k \langle \rho \rangle^{-2k},  \quad |f_1^{(2k+1)}(\rho)| \lesssim_k \rho \langle \rho \rangle^{-2k-2} 
\end{split}
\end{align*}
for all $\rho \geq 0$ and $k\in \N_0$. 
Let 
$\tilde u,\tilde v \in C^{\infty}_{\mathrm{rad},0}(\R^7)$. By definition,
\begin{align}\label{Eq:LipschitzEst1}
\begin{split}
\| f_1(|\cdot|) \tilde u \tilde v \|^2 = & \left \| \Lap[f_1(|\cdot|) \tilde u \tilde v] \right  \|^2_{L^2(\R^7)}  
+ \left \|\Lap^2[f_1(|\cdot|) \tilde u \tilde v]  \right \|^2_{L^2(\R^7)}.
\end{split}
\end{align}
To control the first term, we use 
\begin{align*}
\|   f_1(|\cdot|)  \tilde u \Lap \tilde v\|_{L^2(\R^7)} \lesssim 
\|u \|_{L^{\infty}} \|\Lap \tilde v\|_{L^2(\R^7)} \lesssim \|\tilde u \| \|\tilde v\|.
\end{align*}
Since
\begin{align*}
\|(\cdot)^2  \Lap(f_1  u)  \|_{L^{\infty}}  & \lesssim \|(\cdot)^2  u \Lap f_1   \|_{L^{\infty}} + 
 \| (\cdot)^2  u'  f_1'   \|_{L^{\infty}} +  \| (\cdot)^2    f_1 \Lap u   \|_{L^{\infty}} \\
 & \lesssim \|u  \|_{L^{\infty}} + \| (\cdot)  u' \|_{L^{\infty}} + \| (\cdot)^2  \Lap u   \|_{L^{\infty}} \lesssim \|\tilde u \|,
\end{align*}
we also obtain 
\begin{align*}
\| \tilde v \Lap(f_1(|\cdot|)  \tilde u) \|_{L^2(\R^7)} \lesssim \|(\cdot)^2  \Lap(f_1 u)  \|_{L^{\infty}} \|(\cdot) v \|_{L^2} 
\lesssim \| \tilde u \| \| \tilde v\|. 
\end{align*}
Similarly,
\begin{align*}
  \|(\cdot)^3 (f_1 u)' v'\|_{L^2} & \lesssim  \|(\cdot) (f_1 u)' \|_{L^{\infty}} \|(\cdot)^2 v'\|_{L^2} 
   \lesssim  \left( \|u \|_{L^{\infty}} + \|(\cdot) u' \|_{L^{\infty}} \right)  \|(\cdot)^2 v')\|_{L^2} \lesssim \|\tilde u \| \|\tilde v \|.
\end{align*}
This implies  
\[ \left \| \Lap[f_1(|\cdot|) \tilde u \tilde v] \right  \|_{L^2(\R^7)}  \lesssim \|\tilde u \| \|\tilde v \|. \]
To estimate the fourth order term we use 
\begin{align*}
\| f_1(|\cdot|)  \tilde u  \Lap^2 \tilde v \|_{L^2(\R^7)} \lesssim \| u \|_{L^{\infty}} \|  \Lap^2 \tilde v \|_{L^2(\R^7)} 
\lesssim \|\tilde u \| \| \tilde v \|. 
\end{align*}
An explicit calculation shows that 
\begin{align*}
\|  \tilde v    \Lap^2(f_1(|\cdot|) \tilde u) \|_{L^2(\R^7)} & \lesssim   
\|v u \|_{L^2}  + \| v u' \|_{L^2}  + \|(\cdot) v u'' \|_{L^2}  +\|(\cdot)^2 v u^{(3)} \|_{L^2} + \|(\cdot)^3 v u^{(4)} \|_{L^2}   \\
& \lesssim \|v \|_{L^{\infty}} \left( \|u \|_{L^2} +  \|u'  \|_{L^2} + \|(\cdot) u'' \|_{L^2} + 
 \|(\cdot)^2 u^{(3)} \|_{L^2} + \|(\cdot)^3 u^{(4)} \|_{L^2} \right )\\
 & \lesssim \|\tilde u \| \| \tilde v \|.
\end{align*}

Since 
\[ \|(\cdot) \Lap u \|_{L^2} \lesssim \|(\cdot) u'' \|_{L^2}  + \|u' \|_{L^2}  \lesssim \|\tilde  u \|, \]
we also get  
\begin{align*}
\| \Lap(f_1(|\cdot|)  \tilde u)  & \Lap \tilde v \|_{L^2(\R^7)}  
 \lesssim \|(\cdot)^3 u \Lap v  \|_{L^2} + \|(\cdot)^3 \Lap u \Lap v \|_{L^2}  + \|(\cdot)^3 f_1' u'  \Lap v \|_{L^2} \\
 & \lesssim \|u\|_{L^{\infty}} \|(\cdot)^3  \Lap v \|_{L^2}  + \|(\cdot)^2\Lap u\|_{L^{\infty}} \|(\cdot) \Lap v \|_{L^2} 
 +  \|(\cdot)^2 u' \|_{L^{\infty}} \|(\cdot) \Lap v \|_{L^2} 
 \lesssim \|\tilde u \| \| \tilde v \|.
\end{align*}
Similar estimates for the remaining terms show that 
\begin{align*}
\left \|\Lap^2[f_1(|\cdot|) \tilde u \tilde v]  \right \|_{L^2(\R^7)} \lesssim \|\tilde u \| \|\tilde v \|,
\end{align*}
which yields Eq.~\eqref{Eq:AuxNonl_1}. To prove Eq.~\eqref{Eq:AuxNonl_2} we first convince ourselves that 
\begin{align}\label{Eq:AuxNonl_3}
 \big \|\Lap \big[f_2(|\cdot|) \tilde u \tilde v \tilde w \big] \big \|_{L^2(\R^7)}  \lesssim \| \tilde u \|  \| \tilde v \| \| \tilde w \|.
\end{align}
In fact, 
\begin{align*}
 \big \|\Lap \big[f_2(|\cdot|)\tilde u \tilde v \tilde w \big] \big \|_{L^2(\R^7)}  & 
  \lesssim   \|f_2(|\cdot|) \tilde u \tilde v  \Lap \tilde w   \|_{L^2(\R^7)}   
  +  \big \|\tilde w   \Lap \big[ f_2(|\cdot|) \tilde u \tilde v \big]  \big \|_{L^2(\R^7)}  \\
  & +  \big  \| \nabla \tilde w   \nabla \big[ f_2(|\cdot|) \tilde u \tilde v \big]   \big  \|_{L^2(\R^7)},
\end{align*}
where 
\begin{align*}
 \| f_2(|\cdot|)\tilde u \tilde v  \Lap \tilde w   \|_{L^2(\R^7)}   
 \lesssim \| (\cdot) u \|_{L^{\infty}} \| (\cdot) v \|_{L^{\infty}}  \| \Lap \tilde w   \|_{L^2(\R^7)} \lesssim  \| \tilde u \|  \| \tilde v \| \| \tilde w \|.
\end{align*}
Furthermore, 
\begin{align*}
 \|\tilde w  &  \Lap \big[ f_2(|\cdot|) \tilde u \tilde v \big]  \big \|_{L^2(\R^7)}    \lesssim \|(\cdot)^3 u v w \|_{L^2}   + \|(\cdot)^4 u v' w \|_{L^2}  
 + \|(\cdot)^4  u' v w \|_{L^2}   +   \|(\cdot)^5  u' v' w \|_{L^2}  \\
 & +   \|(\cdot)^5 u v'' w  \|_{L^2}   + \|(\cdot)^5 u'' v w \|_{L^2}   \lesssim  \| (\cdot) w \|_{L^2}  
 \bigg ( \| (\cdot) u \|_{L^{\infty}} \| (\cdot) v \|_{L^{\infty}}     +  \| (\cdot) u \|_{L^{\infty}} \| (\cdot)^2 v' \|_{L^{\infty}}   \bigg )\\
 & + \| (\cdot) w \|_{L^2}  \bigg( \| (\cdot)^2 u' \|_{L^{\infty}} \| (\cdot) v \|_{L^{\infty}}  + 
 \| (\cdot) u \|_{L^{\infty}} \|(\cdot)^3  v'' \|_{L^{\infty}} + \|(\cdot)^3  u'' \|_{L^{\infty}} \| (\cdot) v \|_{L^{\infty}} \
  \bigg)  \\
 & +  \| (\cdot) w \|_{L^2} \| (\cdot)^2 u' \|_{L^{\infty}}  \| (\cdot)^2 v' \|_{L^{\infty}}  \lesssim  \| \tilde u \|  \| \tilde v \| \| \tilde w \|,
\end{align*}
and
\begin{align*}
\big  \| \nabla \tilde w  &  \nabla \big[ f_2(|\cdot|) \tilde u \tilde v \big]   \big  \|_{L^2(\R^7)}    \lesssim  
\|(\cdot)^4 u v w'  \|_{L^2} + \|(\cdot)^5 u v' w'  \|_{L^2} + \|(\cdot)^5 u' v w'  \|_{L^2}  \\
 & \lesssim \|(\cdot)^2 w'  \|_{L^2} \bigg( \| (\cdot) u \|_{L^{\infty}}  \| (\cdot) v \|_{L^{\infty}} 
+\| (\cdot) u \|_{L^{\infty}} \| (\cdot)^2 v' \|_{L^{\infty}} + \| (\cdot)^2 u' \|_{L^{\infty}} \| (\cdot) v \|_{L^{\infty}}  \bigg) \\
& \lesssim  \| \tilde u \|  \| \tilde v \| \| \tilde w \| .
\end{align*} 
This proves Eq.~\eqref{Eq:AuxNonl_3}. To show that 
\begin{align}\label{Eq:AuxNonl_4}
 \big \|\Lap^2 \big[f_2(|\cdot|) \tilde u \tilde v \tilde w \big] \big \|_{L^2(\R^7)}  \lesssim \| \tilde u \|  \| \tilde v \| \| \tilde w \|.
\end{align}
we use 
\begin{align*}
 \| f_2(|\cdot|) \tilde u \tilde v   \Lap^2 \tilde w  \|_{L^2(\R^7)}  
 \lesssim \| (\cdot) u \|_{L^{\infty}} \| (\cdot) v \|_{L^{\infty}} \|(\cdot)^3 \Lap^2 w \|_{L^2} \lesssim  \| \tilde u \|  \| \tilde v \| \| \tilde w \|,
\end{align*}
and
\begin{align*}
 \| \tilde w   &  \Lap^2 [ f_2(|\cdot|) \tilde u \tilde v ]  \|_{L^2(\R^7)}   \lesssim \|(\cdot)^2 u v' w \|_{L^2} +  \|(\cdot)^2 u' v w \|_{L^2}  + 
 \|(\cdot)^3 u' v' w\|_{L^2} + 
 \|(\cdot)^3  u'' v w \|_{L^2}  \\
 & +  \|(\cdot)^3  u v'' w \|_{L^2}   + \|(\cdot)^4 u'' v' w \|_{L^2} +  \|(\cdot)^4 u' v'' w \|_{L^2} +  \|(\cdot)^4 u^{(3)} v w \|_{L^2}
 +  \|(\cdot)^4 u v^{(3)} w \|_{L^2} \\
 &+ \|(\cdot)^5 u^{(4)} v w \|_{L^2}  + \|(\cdot)^5 u v^{(4)} w \|_{L^2}     + \|(\cdot)^5 u^{(3)} v' w \|_{L^2}  
 + \|(\cdot)^5 u' v^{(3)} w \|_{L^2} +  \|(\cdot)^5 u'' v'' w \|_{L^2} 
  \\
 &  \lesssim \|(\cdot)w \|_{L^2}
 \bigg(  \| u \|_{L^{\infty}}  \| (\cdot) v' \|_{L^{\infty}} +  \| (\cdot) u' \|_{L^{\infty}}  \| v \|_{L^{\infty}} + \| (\cdot) u' \|_{L^{\infty}} \| (\cdot) v' \|_{L^{\infty}}
 +  \| (\cdot)^2 u'' \|_{L^{\infty}}  \| v \|_{L^{\infty}}   \\
  & +     \| u \|_{L^{\infty}}  \| (\cdot)^2 v'' \|_{L^{\infty}} + 
 \| (\cdot)^2 u'' \|_{L^{\infty}}   \| (\cdot) v' \|_{L^{\infty}} +
  \| (\cdot) u' \|_{L^{\infty}}   \| (\cdot)^2 v'' \|_{L^{\infty}}   \\
  & +
   \| (\cdot)^3  u^{(3)} \|_{L^{\infty}}   \| v \|_{L^{\infty}}  +      \| u \|_{L^{\infty}}  \| (\cdot)^3  v^{(3)} \|_{L^{\infty}} 
 \bigg) + \|(\cdot)^3 u^{(4)} \|_{L^2} \| (\cdot) v \|_{L^{\infty}}  \| (\cdot) w \|_{L^{\infty}}   \\
 & +  \| (\cdot) u \|_{L^{\infty}} \|(\cdot)^3 v^{(4)} \|_{L^2}  \| (\cdot) w \|_{L^{\infty}} +
 \|(\cdot)^2 u^{(3)} \|_{L^2} \| (\cdot)^2 v' \|_{L^{\infty}}  \| (\cdot) w \|_{L^{\infty}}  \\
 &  +\| (\cdot)^2 u' \|_{L^{\infty}} \|(\cdot)^2 v^{(3)} \|_{L^2}  \| (\cdot) w \|_{L^{\infty}} +
  \| (\cdot)^2 u'' \|_{L^{\infty}} \|(\cdot)^2 v'' \|_{L^2}  \| (\cdot) w \|_{L^{\infty}} 
  \lesssim \| \tilde u \|  \| \tilde v \| \| \tilde w \|.
\end{align*}
Furthermore, 
\begin{align*}
 \| \Lap \tilde w  &  \Lap [ f_2(|\cdot|) \tilde u \tilde v ]  \|_{L^2(\R^7)}   \lesssim \|(\cdot)^3 u v  \Lap w\|_{L^2} +  \|(\cdot)^4 u' v  \Lap w \|_{L^2} 
 + \|(\cdot)^4 u v'  \Lap w \|_{L^2}  \\
 &  +  \|(\cdot)^5 u' v' \Lap w\|_{L^2} \ + \|(\cdot)^5 u'' v \Lap w \|_{L^2} +\|(\cdot)^5 u v'' \Lap w \|_{L^2} \\
 & \lesssim  \|(\cdot)^3 \Lap w \|_{L^2} \bigg (  \|  u \|_{L^{\infty}} \|  v \|_{L^{\infty}} +  \| (\cdot) u' \|_{L^{\infty}} \|  u \|_{L^{\infty}}     + 
 \|  u \|_{L^{\infty}}  \| (\cdot) v' \|_{L^{\infty}} \\
 & + \| (\cdot) u' \|_{L^{\infty}} \| (\cdot)  v' \|_{L^{\infty}} + \| (\cdot)^2 u'' \|_{L^{\infty}} \| v \|_{L^{\infty}}
 +   \| u \|_{L^{\infty}} \| (\cdot)^2 v'' \|_{L^{\infty}} 
\bigg ) \\
 & \lesssim  \| \tilde u \|  \| \tilde v \| \| \tilde w \|.
\end{align*}
These bounds and similar estimates for the remaining terms prove Eq.~\eqref{Eq:AuxNonl_4}.
\end{proof}

\begin{lemma}\label{Le:Nonlinearity}
The nonlinearity extends to a continuous mapping $N: \mc H \to \mc H$ that satisfies 
\begin{align}\label{Eq:Nonlinear_Lip}
  \|  N(\tilde u_1 ) -  N(\tilde u_2) \| \lesssim ( \|\tilde u_1 \| + \|\tilde u_2 \|) \|\tilde u_1 - \tilde u_2\|
\end{align}
for all $\tilde u_1,\tilde u_2 \in \mc B \subset \mc H$. Furthermore, $N$ is differentiable at every 
$\tilde u \in \mc H$ with Fr\'echet-derivative $DN(\tilde u) \in \mc B(\mc H)$ and the mapping $DN : \mc H \to \mc B(\mc H)$
is continuous. 
\end{lemma}

\begin{proof}
Let $\tilde u_1,\tilde u_2 \in  C^{\infty}_{\mathrm{rad},0}(\R^7)$. By Lemma \ref{Le:NonlAux}
\begin{align*}
  \|  N(\tilde u_1 ) -  N(\tilde u_2) \| &  \lesssim \|f_1(|\cdot|)(\tilde u_1 - \tilde u_2)(\tilde u_1 + \tilde u_2) \|  
   + \|f_2(|\cdot|)(\tilde u_1 - \tilde u_2)(\tilde u_1^2 + \tilde u_2^2 + \tilde u_1 \tilde u_2) \|  \\
   & \lesssim \|\tilde u_1 - \tilde u_2\| \big( \|\tilde u_1 + \tilde u_2 \|  + \|\tilde u_1\|^2 +  \|\tilde u_2\|^2 + \|\tilde u_1 \| \| \tilde u_2\| \big).
\end{align*}
Hence, 
\begin{align}\label{Eq:NonlinLip}
  \|  N(\tilde u_1 ) -  N(\tilde u_2) \| \leq \gamma_1(\|\tilde u_1\|, \|\tilde u_2\|)   \|\tilde u_1 - \tilde u_2\| 
\end{align}
for $\gamma_1: [0,\infty) \times [0,\infty) \to [0,\infty)$ a continuous function. 
For $\tilde u_1,\tilde u_2 \in \mc B \cap C^{\infty}_{\mathrm{rad},0}(\R^7)$,   
 \[ \gamma_1(\|\tilde u_1\|, \|\tilde u_2\|) \lesssim \|\tilde u_1 \|  + \| \tilde u_2\|. \]
By density of $C^{\infty}_{\mathrm{rad},0}(\R^7)$ in $\mc H$, $N$ extends to a continuous mapping $N: \mc H \to \mc H$, see for example 
Lemma $3.2$ in \cite{Don11}, such that Eq.~\eqref{Eq:NonlinLip} holds for all $\tilde u_1, \tilde u_2 \in \mc H$. To see that the nonlinearity
is Fr\'echet-differentiable, we first assume 
that $\tilde u, \tilde v \in C^{\infty}_{\mathrm{rad},0}(\R^7)$ and set
\begin{align}\label{Eq:FDerivative}
DN(\tilde u) \tilde v := 2 f_1(|\cdot|) \tilde u \tilde v + 3 f_2(|\cdot|) \tilde u^2 \tilde v.
\end{align}
An application of Lemma \ref{Le:NonlAux} shows that  
\begin{align*}
 \| D N(\tilde u) \tilde v\| \lesssim_{\tilde u} \|\tilde v\|,
\end{align*}
for all $\tilde v \in C^{\infty}_{\mathrm{rad},0}(\R^7)$. By density, $D N(\tilde u)$ extends to a bounded, linear operator on $\mc H$.
Moreover, for $\tilde  u_1,\tilde  u_2, \tilde v  \in  C^{\infty}_{\mathrm{rad},0}(\R^7)$ 
\begin{align*}
  \| D N(\tilde u_1) \tilde v -   D N(\tilde u_2) \tilde v \| \lesssim (1 + \|\tilde u_1 + \tilde u_2\|) \| \tilde v\| \|\tilde u_1 - \tilde u_2\|,
\end{align*}
which extends to all $\tilde v \in \mc H$ by the boundedness of the operators $D N(\tilde u_j)$, $j=1,2$. Hence,
\begin{align}\label{Eq:LipDeriv}
   \| D N(\tilde u_1) -   D N(\tilde u_2) \|_{\mc B(\mc H)} \leq \gamma_2(\|\tilde u_1\|, \|\tilde u_2\|)   \|\tilde u_1 - \tilde u_2\| 
\end{align}
with $\gamma_2: [0,\infty) \times [0,\infty) \to [0,\infty)$ continuous.
Another application of Lemma $3.2$ in \cite{Don11} shows that $DN$ can be extended to a continuous mapping $DN : \mc H \to \mc B(\mc H)$ and that 
Eq.~\eqref{Eq:LipDeriv} remains valid for all $\tilde u_1,\tilde u_2, \tilde v  \in \mc H$. Finally, we convince ourselves that 
$DN(\tilde u)$ is indeed the Fr\'echet-derivative of $N$ at $\tilde u$. Let
$\tilde u, \tilde v \in C^{\infty}_{\mathrm{rad},0}(\R^7)$, $\| \tilde v\| \leq 1$. Then,
\begin{align*}
\| N(\tilde u + \tilde v) - N(\tilde u) - DN(\tilde u) \tilde v \| \lesssim \|f_1(|\cdot|) \tilde v^2 \| + \|f_2(|\cdot|) \tilde u \tilde  v^2 \| 
+ \|f_2(|\cdot|) \tilde  v^3 \| \lesssim_{\tilde u} \|\tilde v \|^2,
\end{align*}
by Lemma \ref{Le:NonlAux}. 
By approximation this bound holds in fact for all $\tilde v\in \mc B$ and any $\tilde u\in \mc H$.
Consequently, we find
\begin{align}\label{Eq:FDeriv}
\lim_{\|\tilde  v \| \to 0 } \frac{ \| N(\tilde u + \tilde v) - N(\tilde u) - DN(\tilde u) \tilde v \| }{\|\tilde  v \|} = 0
\end{align}
 and this implies the claim.
\end{proof}

 \subsection{Initial data operator}
 
In the following, we fix $T_0>0$ and set $\mb W(T, \rho) := \tfrac{T}{T_0} \mb W \left ( \tfrac{\sqrt{T} }{\sqrt{T_0}} \rho \right )$ with $\mb W$
given by Eq.~\eqref {Eq:GroundState}. We define 
 \begin{align*}
U(\tilde v,T) := T \tilde v ( T^{\frac12} \cdot ) + \mb W (T,|\cdot| ) - \mb W( |\cdot| ),
\end{align*}
for $\tilde v \in \mc H$, $T \in [T_0-\delta, T_0 + \delta ]$ and $0 < \delta \leq \frac{T_0}{2}$.

 \begin{lemma}\label{Le:InitialData}
For fixed $\tilde v \in \mc H$, the mapping $T \mapsto U( \tilde  v, T)$: $[T_0-\delta, T_0 + \delta ] \to \mc H$ is continuous.
Furthermore, if $\|\tilde   v \| \leq \delta$ then
\[  \|U( \tilde  v, T) \| \lesssim \delta \]
for all $T \in [T_0- \delta, T_0 + \delta]$. 
\end{lemma}

\begin{proof}
First, we note that $\mb W$ satisfies $|\Lap \mb W(\rho)| \lesssim \langle \rho \rangle^{-4}$, $| \Lap^2 \mb W(\rho)| \lesssim \langle \rho \rangle^{-6}$.
In particular, $\mb W(|\cdot|) \in \mc H$. We prove the statement for $T_0 = 1$ (the general case is analogous).
Fix $\tilde v \in \mc H$ and let $\delta \in (0, \frac12]$. For $T, \tilde T \in  [1 - \delta, 1 + \delta] \subseteq [\frac12, \frac32]$,
we have
\begin{align*}
 \|U( \tilde  v, T)  - U(\tilde  v, \tilde  T)  \| & \leq |T - \tilde T| \|\tilde v(T^{\frac12} \cdot) \| + \tilde T \|\tilde v(T^{\frac12} \cdot) - \tilde v(\tilde T^{\frac12} \cdot) \|
 +|T - \tilde T| \|\mb W(T^{\frac12} |\cdot|) \|  \\
 & \quad +  \tilde T \|\mb W(T^{\frac12} |\cdot|) - \mb W(\tilde T^{\frac12} |\cdot|) \| \\
 & \lesssim  |T - \tilde T| \|\tilde v\| +  |T - \tilde T|  \|\mb W(|\cdot|) \| 
 \|\mb W(T^{\frac12} | \cdot|) - \mb W(\tilde T^{\frac12}  |\cdot|) \| \\
 &  \quad +  \|\tilde v(T^{\frac12} \cdot) - \tilde v(\tilde T^{\frac12} \cdot) \|
\end{align*}
by rescaling. The first three terms tend to zero in the limit $T \to \tilde T$. For the last term, we use that for all $\tilde v, \tilde u \in \mc H$,
\begin{align*}
\| \tilde v(T^{\frac12} \cdot) - \tilde v( \tilde T^{\frac12} \cdot) \| & \leq \|\tilde v(T^{\frac12} \cdot) - \tilde u ( T^{\frac12} \cdot) \|
+ \| \tilde u (T^{\frac12} \cdot) - \tilde  u ( \tilde T^{\frac12} \cdot) \|
+ \| \tilde u( \tilde T^{\frac12} \cdot)  -  \tilde v( \tilde T^{\frac12} \cdot) \| \\
& \lesssim \|\tilde  v - \tilde u \| + \| \tilde u(T^{\frac12} \cdot) - \tilde u( \tilde T^{\frac12} \cdot) \|.
\end{align*}
Let $\varepsilon > 0$ be arbitrary. By density, there is a $\tilde u \in C^{\infty}_{\mathrm{rad},0}(\R^7)$, such that 
\[ \| \tilde v - \tilde u \| < \varepsilon. \]
Since $\tilde u$ is a smooth, compactly supported function, 
\[ \lim_{T \to \tilde T} \| \tilde u(T^{\frac12} \cdot) - \tilde u ( \tilde T^{\frac12} \cdot) \| = 0. \]
This implies the continuity of $U(\tilde v, \cdot): [1-\delta, 1 + \delta ] \to \mc H$. 

For $\tilde v \in \mc H$ satisfying $\| \tilde v \| \leq \delta$
we get 
\begin{align*}
 \| U(\tilde v, T ) \| &\leq T \| \tilde v(T^{\frac12} \cdot) \| + \|T\mb W (T^{\frac12} |\cdot|) - \mb W(|\cdot|) \| \lesssim \| \tilde v\| +  \|T \mb W (T^{\frac12}  |\cdot|) - \mb W(|\cdot|) \| \\
 &\lesssim \delta+|T-1| \lesssim \delta
\end{align*}
for all $T\in [1-\delta,1+\delta]$ since $\|T\mb W(T^\frac12 |\cdot|)-\mb W(|\cdot|)\|\lesssim |T-1|$.
\end{proof}

\subsection{Operator formulation of Eq.~\eqref{Eq:YM_AbstractPertubation}}

With the above definitions, Eq.~\eqref{Eq:YM_AbstractPertubation} can now be considered as an abstract
initial value problem on $\mc H$. The corresponding integral equation reads
 \begin{align}\label{Eq:IntegralEq}
 \Phi(\tau) = S(\tau) U(\tilde v, T )  + \int_0^{\tau} S(\tau - \tau') N ( \Phi(\tau')) d \tau',
 \end{align}
where $\{S(\tau): \tau \geq 0\}$ is the semigroup generated by $(L, \mc D(L))$, see Corollary \ref{Cor:SemigroupS}. 
We set $a=  \frac{1}{150}$ and introduce the Banach space
\begin{align}
\mc X := \{ \Phi \in C([0,\infty), \mc H) : \| \Phi \|_{\mc X} := \sup_{\tau \geq 0}  e^{a \tau } \|\Phi(\tau) \|   < \infty \}.
\end{align}
By $\mc X_{\delta}$ we denote the closed subspace 
\[ \mc X_{\delta} := \{ \Phi \in \mc X:  \| \Phi \|_{\mc X} \leq \delta \}. \]

From now on we proceed exactly as in our previous works on radial wave equations, see for example \cite{DonSch15a},
to prove Theorem \ref{Th:MainSim}. For convenience of the reader, we repeat the main arguments. 

\subsection{Correction of the unstable behavior}

We define
\begin{align*}
C(\Phi, \tilde u) :=  P \tilde u + \int_0^{\infty} e^{-\tau'} P N(\Phi(\tau')) d\tau',
\end{align*}
and set
\begin{align*}
 K(\Phi, \tilde u)(\tau) := S(\tau) \tilde u  + \int_0^{\tau}  S(\tau - \tau')  N(\Phi(\tau')) d\tau' - e^{\tau} C(\Phi, \tilde u).
\end{align*}

\begin{lemma} \label{Le:GlobalEx_ModEq}
Choose $\delta > 0$ sufficiently small and $c > 0$ sufficiently large.
For every $\tilde  u\in \mc H$ with $\| \tilde u \| \leq \frac{\delta}{c}$, there exists a unique $\Phi(\tilde u) \in \mc X_{\delta}$ that satisfies
\[\Phi(\tilde u)   = K (  \Phi(\tilde u) , \tilde u). \]
Furthermore, the mapping $ \tilde u \mapsto \Phi(\tilde u)$ is continuous.
\end{lemma}

\begin{proof}
We first convince ourselves that for fixed $\tilde u \in \mc H$, we have $K(\cdot, \tilde u) :  \mc X_{\delta} \to  \mc X_{\delta}$. 
For $\Phi \in \mc X_{\delta}$, the continuity of $K(\Phi, \tilde u)(\tau)$ in $\tau$ is a consequence of the strong continuity of the semigroup.
We write $K(\Phi, \tilde u)(\tau) = PK(\Phi, \tilde u)(\tau) + (1-P)K(\Phi, \tilde u)(\tau)$ and use Lemma \ref{Le:Nonlinearity} together with $N(0)=0$ to estimate
\begin{align*}
\| P K(\Phi, \tilde u)(\tau)\| & \lesssim  \int_{\tau}^{\infty} e^{ -(\tau' - \tau)} \| P N(\Phi(\tau')) \| d\tau'  \\
& \lesssim
\int_{\tau}^{\infty} e^{ -(\tau' - \tau)} \| \Phi(\tau') \|^2 d\tau'  
\lesssim e^{- 2 a \tau} \|\Phi \|^2_{\mc X} \leq
e^{-2a \tau} \delta^2 ,
\end{align*}
\begin{align*}
 \| (1 -P) K(\Phi, \tilde u)(\tau)\| \lesssim e^{- a \tau} \|\tilde  u \| +  \int_{0}^{\tau} e^{- a(\tau -\tau')} \|N(\Phi(\tau')) \| d\tau'  
\lesssim e^{- a \tau} ( \tfrac{\delta}{c} +  \delta^2).
\end{align*}
Thus, 
\[  \|K(\Phi, \tilde u)(\tau)\| \leq  e^{- a \tau} \delta \]
by choosing $\delta > 0$ sufficiently small and $c >0$ sufficiently large (independent of $\delta$). 
For $\Phi, \Psi \in  \mc X_{\delta}$ we use the fact that 
\[ \| N(\Phi(\tau)) - N(\Psi(\tau))  \|   \lesssim   \delta e^{-2 a \tau}  \| \Phi  - \Psi \|_{\mc X}  \]
for all $\tau > 0$ to infer that 
\begin{align*}
\| P[  K(\Phi,\tilde u)(\tau) - K(\Psi,\tilde u)(\tau)] \| &  \lesssim   \delta  e^{-2a \tau} \| \Phi - \Psi \|_{\mc X}  
\end{align*}
and 
\begin{align*}
\|(1 - P) [K(\Phi, \tilde u)(\tau) -  K(\Psi, \tilde u)(\tau)]\| & \lesssim \int_{0}^{\tau} e^{- a(\tau -\tau')}  \|  N(\Phi(\tau')) -  N(\Psi(\tau'))  \|  d\tau' \\
& \lesssim  \delta  e^{-a \tau} \| \Phi - \Psi \|_{\mc X}.
\end{align*}
This implies that there is a $0 < k < 1$ such that 
\[ \| K(\Phi, \tilde  u) -  K(\Psi, \tilde u)\|_{\mc X} \leq k \| \Phi - \Psi \|_{\mc X}, \]
provided $\delta > 0$ is sufficiently small. Since $\mc X_{\delta} \subset \mc X$ is closed, we can apply 
the Banach fixed point theorem to infer the existence of a unique solution $\Phi_{\tilde u}$ to the equation
\[\Phi  = K (\Phi ,\tilde u). \]
Standard arguments show that the mapping $\tilde  u \mapsto \Phi(\tilde u) := \Phi_{ \tilde u}$ is continuous. 
\end{proof}

\subsection{Proof of Theorem \ref{Th:MainSim}}

Let $\tilde v \in \mc H$ with $\|\tilde  v \| \leq \frac{\delta}{M^2}$, for $0 < \delta \leq \frac{T_0}{2}$. By Lemma \ref{Le:InitialData}, 
 \[ \| U(\tilde v, T)\| \leq   \tfrac{K \delta}{M}, \]
for all $T \in I_{M,\delta}:= [T_0 - \frac{\delta}{M}, T_0 + \frac{\delta}{M }]$ and some $K > 0$. By choosing $M$ sufficiently large, we obtain
 \[ \| U(\tilde v, T)\| \leq  \tfrac{\delta}{c}, \]
for all $T \in I_{M,\delta}$, where $c >0$ is the constant from Lemma \ref{Le:GlobalEx_ModEq}.
Let $\delta > 0$ be sufficiently small such that Lemma \ref{Le:GlobalEx_ModEq} applies.
Hence, for every $T \in I_{M,\delta}$ there exists a unique solution $\Phi_T := \Phi(U(\tilde v, T)) \in \mc X_{\delta}$ to the equation
\begin{align*}
\Phi_T(\tau) =  S(\tau) U(\tilde v, T)   + \int_0^{\tau}  S(\tau - \tau')   N( \Phi_T (\tau')) d\tau' - e^{\tau} C(\Phi_T,U( \tilde v, T)).
\end{align*}
Furthermore, since the mappings $T \mapsto U(\tilde v, T)$ and $\tilde u \mapsto \Phi(\tilde u)$ are continuous, we 
see that $T \mapsto \Phi_T$ is continuous. 
We show that there is a $T_{\tilde v} \in I_{M,\delta}$ such that 
$C (\Phi_{T_{\tilde  v}}, U(\tilde  v, T_{\tilde  v} )) = 0$, which implies that $\Phi_{T_{\tilde v}} \in \mc X_{\delta}$ is a solution to the 
original problem given by Eq.~\eqref{Eq:IntegralEq}.

Since $\mathrm{rg} P \subseteq \ \mathrm{span}(\mb{ \tilde  g})$, with
$\mb{ \tilde  g}(\xi) = \mb g(|\xi|)$ denoting the symmetry mode, see Lemma \ref{Le:LinearTimeEvol}, it suffices to show that
\begin{align}\label{Eq:CorrectionTerm}
\big (C (\Phi_{T_{\tilde  v}} , U(\tilde  v, T_{\tilde  v} )) \big |\mb{ \tilde  g} \big ) = 0,
\end{align}
for some $T_{\tilde  v} \in I_{M,\delta}$.
First, we estimate
\begin{align*}
 \left ( \int_0^{\infty} e^{-\tau} P N( \Phi_{T}(\tau))  d\tau \bigg |  \mb{ \tilde  g} \right ) & 
 \lesssim  \| \mb{ \tilde  g}\|  \int_0^{\infty} e^{-\tau} \|  P N( \Phi_{T}(\tau))\|  d\tau   
\lesssim   \int_0^{\infty} e^{-\tau} \| \Phi_{T}(\tau)\|^2  d \tau  \lesssim \delta^2 .
\end{align*}
Next, observe that 
\begin{align*}
\partial_T \mb W(T,\cdot)|_{T=T_0} =  \alpha \mb g
\end{align*}
for some $\alpha \in \R$, $\alpha \neq 0$. By the fundamental theorem of calculus we can write 
\begin{align*}
\mb W(T,\cdot) &  = \mb W(T_0,\cdot) +   (T-T_0) \partial_T \mb W(T,\cdot)|_{T=T_0}   + \int_{T_0}^{T} \left ( \int_{T_0}^{S}  \partial^2_s \mb W(s,\cdot) ds \right) dS 
\\ & =  \mb W +   \alpha (T-T_0) \mb  g + (T-T_0)^2 R(T, \cdot).
\end{align*}
It is easy to check that $\| R(T, |\cdot|) \| \lesssim 1$ for all $T \in I_{M,\delta}$ and that 
 $R$ depends continuously on $T$. Thus,
\[ U(\tilde  v,T) = T  \tilde v ( T^{\frac12} \cdot ) + \alpha (T-T_0) \mb{ \tilde g }+ (T-T_0)^2  R(T, |\cdot|),   \]
which implies 
\begin{align*}
(P U(\tilde  v,T)| \mb{ \tilde g })  & = T \big (P \tilde  v(T^{\frac12} \cdot)\big |\mb{ \tilde g }\big ) +  \alpha (T-T_0)  \|\mb{ \tilde g } \|^2
+  (T-T_0)^2 \big ( P R(T, |\cdot|) \big |\mb{ \tilde g } \big )\\
& = \alpha (T-T_0)  \|\mb{ \tilde g }\|^2 + f(T),
\end{align*}
with $|f(T)| \lesssim \frac{\delta}{M ^2}+\delta^2$. Consequently, we can write 
Eq.~\eqref{Eq:CorrectionTerm} as
\[ T = T_0 + F(T)\]
for a continuous function $F$ that satisfies $|F(T)| \lesssim \delta^2 +  \frac{\delta}{M ^2}$ for all $T \in I_{M,\delta}$.
Now, we choose $M$ sufficiently large (independent of $\delta$) and $\delta$ sufficiently small (now depending on $M$) to ensure that $|F(T)| \leq \frac{\delta}{M}$. 
Consequently, $T \mapsto T_0 + F(T)$ is a continuous function that maps the interval $I_{M,\delta}$ to itself. 
By Brouwer's fixed point theorem, there exists a $T_{\tilde v} \in [T_0 - \frac{\delta}{M }, T_0 + \frac{\delta}{M}]$ such that 
\[T_{\tilde v} = T_0 + F(T_{\tilde v}).\]
This proves that $\Phi_{T_{\tilde v}} \in \mc X_{\delta}$ satisfies Eq.~\eqref{Eq:IntegralEq}.
For the uniqueness of the solution in $C([0,\infty), \mc H)$, we refer the reader for example to the proof of 
Theorem 4.11 in \cite{DonSch12}.

\subsection{Theorem \ref{Th:MainSim} implies Theorem \ref{Th:Main}}\label{Subsec:FinalArg}
For fixed $T_0 > 0$, we choose $\delta, M > 0$ such that Theorem \ref{Th:MainSim} holds. Set $\delta' = \frac{\delta}{M}$. For $u_0 \in \mc E$ and  $\mb w_T$ defined as in Eq.~\eqref{Eq:WeinkoveSol}, we have
\begin{align*}
\|u_0(|\cdot|) - \mb w_{T_0}(0,|\cdot|) \| = C \|u_0 - \mb w_{T_0}(0,\cdot) \|_{\mc E} , 
\end{align*}
where the constant comes from the integration over $\mathbb{S}^6$. Assume that 
\begin{align}\label{Eq:Assump}
 \|u_0 - \mb w_{T_0}(0,\cdot) \|_{\mc E}   \leq \frac{\delta'}{K} 
\end{align}
for $K = CM$. We set  
\[ \tilde v_0 := u_0(|\cdot|) - \mb w_{T_0}(0,|\cdot|).\]
By definition of the initial data operator,
we have 
\[ U( \tilde v_0 , T) =  T u_0(T^{\frac12} |\cdot|) - \mb W(|\cdot|) =:\Phi^T_0.\]
Obviously, $\Phi^T_0 \in  C^{\infty}_{\mathrm{rad}}(\R^7)$ for all $T \in [T_0 - \delta', T_0 + \delta']$, 
$\| \Phi^T_0 \| < \infty$, and the decay of $\mb W$ implies
$\Phi^T_0\in \mc H$ by an approximation argument. It is also easy to check that 
$\Phi^T_0 \in \mc D(\tilde L_0) \subset \mc D(L)$.
Eq.~\eqref{Eq:Assump} implies that $\tilde v_0$ satisfies the assumptions of Theorem \ref{Th:MainSim}. 
Hence, there is a $T  \in [T_0 - \delta', T_0 + \delta']$ such that there exists a unique solution
$\Phi \in C([0,\infty),\mc H)$ to
\begin{align*}
\Phi(\tau) = S(\tau)\Phi^T_0  + \int_0^{\tau} S(\tau - \tau') N ( \Phi(\tau')) d \tau', \quad \tau \geq 0,
\end{align*}
where $\Phi(\tau)(\cdot) = \varphi(\tau,|\cdot|)$, $\varphi(\tau, \cdot) \in C[0,\infty) \cap C^2(0,\infty)$ and 
\[ \| \Phi(\tau) \| \lesssim e^{-\frac{1}{150} \tau} , \quad \forall \tau \geq 0. \]
Lemma \ref{Le:Nonlinearity} implies that $\Phi$ is also a classical solution, see Theorem 6.1.5 in \cite{Pazy}. This means that 
$\Phi: (0,\infty) \to \mc H$ is continuously differentiable, $\Phi(\tau) \in \mc D(L)$ for all $\tau > 0$ and 
\begin{align*}
\frac{d}{d\tau} \Phi(\tau) & = (L_0 + L') \Phi(\tau) + N(\Phi(\tau)) \quad \tau > 0, \\
\Phi(0) & = \Phi^T_0.
\end{align*}
Recall that $L_0$ acts as a classical 
differential operator on functions in $\mc D(L)$. 
By setting \[ \psi(\tau,\rho) := \mb W(\rho) + \varphi(\tau,\rho), \]
we obtain a classical solution to Eq.~\eqref{Eq:SelfSim} corresponding to the initial condition 
$\psi(0,\cdot) = T u_0(T^{\frac12} \cdot)$. As a consequence, 
\[ u(t, r) := (T-t)^{-1} \psi(- \log(T-t) + \log T , \tfrac{r}{\sqrt{T-t}}) \]
solves Eq.~\eqref{Eq:EquivarEq5} in a classical sense for all $0 < t < T$. Furthermore, $u(0,r) = u_0(r)$ for all $r \in [0,\infty)$.

Note that 
\[ \| \Lap \mb w_T(t,|\cdot|)\|_{L^2(\R^7)}  = c_1 (T-t)^{-\frac{1}{4}}, \quad  
\| \Lap^2 \mb w_T(t,|\cdot|)\|_{L^2(\R^7)} = c_2 (T-t)^{-\frac{5}{4}},   \]
for some constants $c_1,c_2 > 0$. By definition and rescaling, we get that for $k=1,2$, 
\begin{align*}
(T-t)^{-\frac{3}{4} + k} &  \|\Lap^k u(t,|\cdot|)   - \Lap^k \mb w_T(t,|\cdot|)\|_{L^2(\R^7)}  \\
& = \| \Lap^k \psi(-\log(T-t)+ \log T, |\cdot|) - \Lap^k \mb W(|\cdot|)\|_{L^2(\R^7)} \\
& = \|\Lap^k \Phi(-\log(T-t)+ \log T) \|_{L^2(\R^7)} \leq \| \Phi(-\log(T-t)+ \log T)\| \lesssim (T-t)^{a},
\end{align*}
which implies Eq. \eqref{Eq:ConvSolHnorm}. Furthermore, 
\[ \sup_{r > 0 }| \mb w_T(t,\cdot)| = \frac{1}{T-t} \sup_{r > 0 } \mb W\left (\tfrac{r}{\sqrt{T-t}} \right) = 
\frac{1}{T-t} \| \mb W \|_{L^{\infty}(\R^+)} = \frac{c_3}{T-t},  \]
for some $c_3 > 0$. By Lemma \ref{Le:PropH}
\begin{align*}
(T-t) \| u(t,\cdot) -  \mb w_T(t,\cdot)\|_{L^{\infty}(\R^+)}  & \lesssim \|  \psi(-\log(T-t)+ \log T, \cdot) - \mb W \|_{L^{\infty}(\R^+)} \\
& \lesssim \| \Phi(-\log(T-t)+ \log T)\| \lesssim (T-t)^{a},
\end{align*}
This finishes the proof.

 \appendix

\section{Hardy's inequality and proof of Lemma \ref{Le:AllBounds}}

\subsection{Hardy's inequality}

We give a proof of Hardy's inequality for integrals defined on $\R^+$. However, the same statement is true 
for arbitrary intervals $(0,a) \subset \R$, $a > 0$. 

\begin{lemma}\label{Le:Hardy}
Let $\alpha \in \N$. Assume that $f \in C^1[0,\infty)$ satisfies $\lim_{\rho \to 0}  \rho^{-2\alpha+1} |f(\rho)|^2 = 0$.
Then,
\[ \|  (\cdot)^{-\alpha} f \|_{L^2(\R^+)} \lesssim \|  (\cdot)^{-\alpha+1} f' \|_{L^2(\R^+)},\]
provided the right hand side is finite. 
\end{lemma}
\begin{proof}
For $f = 0$, the assertion is trivial. Let $f \neq 0$. Integration by parts and the Cauchy-Schwarz inequality yield
\begin{align*}
& \int_0^{a} \rho^{-2\alpha} |f(\rho)|^2 d \rho   \leq  \lim_{\rho\to 0} \left ( |2\alpha-1|^{-1} \rho^{-2\alpha+1} |f(\rho)|^2  \right)    \\
& + |2\alpha-1|^{-1} \int_0^{a}  \rho^{-2\alpha+1} ( f'(\rho) \overline{f(\rho)} + f(\rho) \overline{ f'(\rho)} )d\rho 
\lesssim   \int_0^{a}
\rho^{-2\alpha+1}  |f'(\rho)||\overline{f(\rho)}| d\rho    \\
& \lesssim  \left( \int_0^{a}  \rho^{-2\alpha} |f(\rho)|^2 d\rho\right)^{1/2}  \left( \int_0^{a} \rho^{-2\alpha+2} |f'(\rho)|^2  d\rho \right)^{1/2}.
\end{align*}
The claim follows by letting $a \to \infty$. 
\end{proof}

\subsection{Proof of Lemma \ref{Le:AllBounds}}\label{Sec:ProofAllBounds}
First, we note that $\tilde u = u(|\cdot|) \in C^4_{\mathrm{rad}}(\R^7)$ implies that $u'(0) = u'''(0)=0$. 
Our assumptions allow us to integrate by parts to obtain
\begin{align*}
\| \nabla \Lap \tilde  u \|^2_{L^2(\R^7)} \simeq \int_0^{\infty} \rho^{6} |(\Lap_{\mathrm{rad}} u)'(\rho)|^2 d\rho
\lesssim \int_0^{\infty} \rho^{6} |\Lap_{\mathrm{rad}} u(\rho)||\Lap^2_{\mathrm{rad}} u(\rho)| d\rho \lesssim  \| \tilde u \|^2.
\end{align*}
For all $\rho \geq 1$, we have
\begin{align*}
\begin{split}
 \rho^6 |(\Lap_{\mathrm{rad}}u)'(\rho)|^2 & \lesssim \int_1^{\infty} \rho^6 |(\Lap_{\mathrm{rad}}u)'(\rho)|^2 d\rho +
\int_1^{\infty} \rho^6 |(\Lap_{\mathrm{rad}}u)''(\rho)|^2 d\rho  \\
& \lesssim \| \nabla \Lap \tilde  u \|^2_{L^2(\R^7)}  
+ \| \Lap^2 \tilde  u \|^2_{L^2(\R^7)} \lesssim \| \tilde u \|^2.
\end{split}
\end{align*}
Similarly, we obtain 
\begin{align*}
\begin{split}
 \rho^6 |\Lap_{\mathrm{rad}}u(\rho)|^2 \lesssim \| \Lap \tilde  u \|^2_{L^2(\R^7)}  +  \| \nabla \Lap \tilde  u \|^2_{L^2(\R^7)}  
 \lesssim \| \tilde u \|^2,
\end{split}
\end{align*}
which implies that 
\begin{align}\label{Eq:EqSomebound}
\|(\cdot)^3 \Lap_{\mathrm{rad}}  u \|_{L^{\infty}(1,\infty)} \lesssim \| \tilde  u \|, 
\quad  \|(\cdot)^3 (\Lap_{\mathrm{rad}}  u)' \|_{L^{\infty}(1,\infty)} \lesssim \| \tilde  u \|.
\end{align}
Assume that in addition to Eq.~\eqref{Eq:CondLeAllB1} also
\[ \lim_{\rho \to \infty} \rho^{\frac32} |u(\rho)| = 0, \quad  \lim_{\rho \to \infty} \rho^{\frac52} |u'(\rho)| = 0. \]
Then 
\[ \lim_{\rho \to \infty} \rho^{3} |u''(\rho)| = 0, \quad  \lim_{\rho \to \infty} \rho^{3} |u^{(3)}(\rho)| = 0. \]
We prove Eq.~\eqref{Eq:L2Bounds} and 
first show that 
\[\| (\cdot)^2 u' \|_{L^2(\R^+)}  \lesssim \|\Lap \tilde u \|_{\rad}, \quad  \| (\cdot)^3 u'' \|_{L^2(\R^+)}  \lesssim \|\Lap \tilde u \|_{\rad}. \]
By Hardy's inequality, see Lemma \ref{Le:Hardy},
\begin{align*}
\int_0^{\infty} \rho^4 | u'(\rho) |^2 d\rho = \int_0^{\infty} \rho^{-8} |\rho^6 u'(\rho) |^2 d\rho \lesssim
\int_0^{\infty} \rho^{6} |\rho^{-6} \tfrac{d}{d\rho} (\rho^6 u'(\rho)) |^2 d\rho = \|\Lap \tilde  u \|^2_{\rad}.
\end{align*}
Furthermore,  
\[  \| (\cdot)^3 u'' \|^2_{L^2(\R^+)} \lesssim \|\Lap \tilde  u \|^2_{\rad} +  \| (\cdot)^2 u' \|^2_{L^2(\R^+)} \lesssim
\|\Lap \tilde  u \|^2_{\rad}  .\]
An integration by parts and the Cauchy-Schwarz inequality show that 
\begin{align*}
\int_0^{\infty} \rho^2 | u(\rho) |^2 d\rho \lesssim  \int_0^{\infty} \rho^3 |u (\rho) | |u'(\rho) | d\rho 
\lesssim   \left ( \int_0^{\infty} \rho^2 | u(\rho) |^2 d\rho \right)^{\frac12} \left ( \int_0^{\infty} \rho^4 | u'(\rho) |^2 d\rho \right)^{\frac12},
\end{align*}
and we get 
\[ \| (\cdot) u \|_{L^2(\R^+)} \lesssim \| \tilde u \|. \]

These bounds imply that 
\[ \| (\cdot) u' \|_{L^2(1,\infty)}  \lesssim \| \tilde  u \|, \quad \| (\cdot)^2 u'' \|_{L^2(1,\infty)}  \lesssim \| \tilde  u \|,\]
and 
\[ \| (\cdot)^3 u^{(3)} \|_{L^2(1,\infty)}  \lesssim  \| \tilde  u \| + \| \nabla \Lap \tilde  u \|_{L^2(\R^7)} \lesssim  \| \tilde  u \|. \]
Furthermore,
\[ \| (\cdot)^{j-1} u^{(j)} \|_{L^2(1,\infty)}  \lesssim \| \tilde  u \| \]
for $j=1,\dots,4$, since
\[  \| (\cdot)^{3} u^{(4)} \|_{L^2(1,\infty)}   \lesssim  \| \tilde  u \|  +  \|\Lap^2 \tilde  u \|_{L^2(\R^7)} \lesssim   \| \tilde  u \|. \]

To proceed, we define $Du$ by
\[Du(\rho) := \left (\rho^{-1} \frac{d}{d\rho}\right )^2 (\rho^5 u(\rho)) .\]
The assumptions on $\tilde u$ imply that $Du \in C^2[0,\infty)$, $Du(0) = [Du]''(0) =0$ and
\[ \lim_{\rho \to \infty} Du(\rho) = \lim_{\rho \to \infty} [Du]'(\rho) = 0.\]
Using integration by parts, a tedious, but straightforward calculation shows that 
\[ \|[Du]'' \|_{L^2(\R^+)} \simeq  \left \|\Lap^2 \tilde  u   \right \|_{\rad}. \]
We define,
\[Ku(\rho) := \rho^{-5} \mc K^2u(\rho), \quad \mc Ku(\rho) : = \int_0^{\rho} s u(s) ds.  \]
Since $\mathrm{ker} D = \mathrm{span}(\rho^{-3}, \rho^{-5})$, we have 
\[ D K u = K D u = u.\]
We set $w := Du$. Then, $[\mc K^2 w]^{(j)}(0)=0$ for $j=0,\dots, 4$ and $|[\mc K^2 w]^{(5)}(0)|=c|w'(0)|$ for some $c > 0$. 
We prove that 
\begin{align}\label{Eq:SomeEq1}
\int_0^{1}|u'(\rho)|^2 d\rho \lesssim \| \tilde u \|^2.
\end{align}
For this, we use the fundamental theorem of calculus and the fact that $w(0) = 0$ to write
\[ w(\rho) =  \rho w'(0) + \int_0^{\rho} \int_0^{s} w''(t) dt ds. \]
By setting $\mc V w(\rho):= \int_0^{\rho} w(s)ds$ and applying $K$, we get 
\begin{align*}
K w(\rho)  = K \mc V^2 w''(\rho) + c w'(0)
\end{align*}
for some constant $c > 0$. Repeated application of Hardy's inequality yields
\begin{align*}
\int_0^{1}   |(Kw)'(\rho)|^2 d\rho  & = \int_0^{1} \left |[K \mc V^2 w'']'(\rho) \right |^2 d\rho  \\
& \lesssim \int_0^{1}\rho^{-12} |[\mc K^2 \mc  V^2 w''](\rho)|^2 d\rho 
+ \int_0^{1} \rho^{-8} |\mc K \mc V^2w''(\rho)|^2 d\rho   \lesssim \int_0^{\infty}  | w''(\rho)|^2 d\rho ,
\end{align*}
since for $\rho \to 0$,  $[\mc K^2 \mc  V^2 w''](\rho) = \mc O(\rho^{6})$, $[\mc K \mc  V^2 w''](\rho) = \mc O(\rho^{4})$, $[\mc  V^2 w''](\rho) = \mc O(\rho^{2})$, $[\mc  V w''](\rho) = \mc O(\rho)$.
This implies Eq.~\eqref{Eq:SomeEq1}. A similar calculation shows that  
\[ \| (\cdot) u'' \|_{L^2(0,1)}  \lesssim \| \tilde  u \|, \quad \| (\cdot)^2 u^{(3)}\|_{L^2(0,1)}  \lesssim \|\tilde  u \|, \]
and thus, we also get that 
\[ \| (\cdot)^3 u^{(4)} \|_{L^2(0,1)} \lesssim \| \tilde  u \|.\]
Combining this with the  above bounds we infer that 
\[  \| (\cdot)^{j-1} u^{(j)} \|_{L^2(\R^+)} \lesssim \| \tilde u \| \]
for $j=1,\dots,4$. An integration by parts yields
\begin{align*}
\int_0^{\infty} |u(\rho)|^2 d\rho \lesssim  \int_0^{\infty} \rho |u(\rho)| |u'(\rho)| d\rho 
\lesssim \|(\cdot) u \|_{L^2(\R^+)}^2 + \| u' \|_{L^2(\R^+)}^2 \lesssim \| \tilde u \|^2.
\end{align*}
In summary, we get that
\[  \| (\cdot)^{j} u^{(j)} \|_{L^2(\R^+)} \lesssim \| \tilde u \| \]
for $j = 0, \dots, 3$. It is left to verify the $L^{\infty}$-bounds.  
We use the decay of $u$ and its derivative at infinity to write 
\begin{align*}
u(\rho) =  \int_{\rho}^{\infty}  \int_{s}^{\infty}  u''(t)  dt ds,
\end{align*}
for $\rho \geq 1$. 
Taking absolute values and applying the Cauchy-Schwarz inequality yields
\[ |u(\rho)| \lesssim    {\rho}^{-\frac{3}{2}}  \left(\int_{\rho}^{\infty} s^6  |u''(s)|^2  ds  \right)^{\frac12},  \] 
and in view of the above results we get 
$ \| (\cdot)^{\frac32} u \|_{L^{\infty}(1,\infty)} \lesssim \|\tilde u \|.$
Similarly,
\[  |u'(\rho)| \lesssim    {\rho}^{-\frac{5}{2}} \left(\int_{\rho}^{\infty} s^6 |u''(s)|^2  ds \ \right)^{\frac12} ,
\]
for all $\rho \geq 1$  such that 
$\| (\cdot)^{\frac52} u' \|_{L^{\infty}(1,\infty)} \lesssim \|\tilde u \|.$
The $L^2$-bounds and the Sobolev embedding $H^1(0,1) \hookrightarrow L^{\infty}(0,1)$ yield
\begin{align}\label{Eq:BoundZero2}
\| (\cdot)^{j} u^{(j)} \|_{L^{\infty}(0,1)} \lesssim \|\tilde  u \|, \quad j = 0, \dots, 3.
\end{align}
Thus,  
\[ \| (\cdot)^2 \Lap_{\mathrm{rad}}  u\|_{L^{\infty}(0,1)}   \lesssim \|\tilde  u \|, \quad 
\| (\cdot)^3 \Lap_{\mathrm{rad}}  u\|_{L^{\infty}(0,1)}   \lesssim \|\tilde  u \|. \]
These bounds together with Eq.~\eqref{Eq:EqSomebound} imply the estimates given in Eq.~\eqref{Eq:LinftyBounds}.

\section{Representation of $T(\omega)$ - Proof of Lemma \ref{Le:ResEst_T} }\label{App_ProofTLemma}

We write $\mc T(\omega) = \mc T_1(\omega) + \mc T_2(\omega)$ where 
\[ [\mc T_1(\omega)f](r) := \int_r^{\infty} \partial_r g_1(r,s,\omega) f(s) ds, \quad 
 [\mc T_2(\omega)f](r) := \int_0^{r} \partial_r g_2(r,s,\omega) f(s) ds.\]
 Lemma  \ref{Le:ResEst_T} is a direct consequence of the following two results that
 provide representations of $\mc T_1(\omega)$ and  $\mc T_2(\omega)$.

\subsection{Representation of $\mc T_1(\omega)$}

\begin{lemma}\label{Le:ResEst_T1}
Let $\mu = b + \I \omega$, with $b > -4$ fixed. Define $\delta \in [0,\frac12)$ as in Eq.~\eqref{Def:delta}.
Let $f \in C^{\infty}_{\mathrm{e},0}(\R)$. Then, for $m=0,1$, 
\begin{align}\label{Eq:LeT1_1}
\begin{split}
 r^{m+2} [\mc T_1(\omega)f]^{(m)}(r)     & =  O(r \omega^0)  f(r)  + O(r^2 \omega^0) f'(r) \\
& +\sum_{k_1=0}^{2} \int_r^{\infty} O(r^{-\delta} s^{-1 + \delta} \omega^0)   s^{k_1+1} f^{(k_1)}(s) ds \\
&+ \sum_{k_2=0}^{3} \int_r^{\infty} O(r^{-\delta} s^{-1 + \delta} \omega^0)  s^{k_2} f^{(k_2)}(s) ds,
\end{split}
\end{align}
for all $r > 0$ and all $\omega \gg 1$. Furthermore, for $m=0,\dots,3$,
\begin{align}\label{Eq:LeT1_2}
\begin{split}
 r^m [\mc T_1(\omega)f]^{(m)}(r)  & =  \sum_{j=0}^{2}  O(r^{j+1} \omega^0)  f^{(j)}(r) + \sum_{k=0}^{3}  O(r^{k} \omega^0)   f^{(k)}(r) \\
& + \sum_{j_1=0}^2 \int_r^{\infty} O(r^{-\delta} s^{-1 + \delta}\omega^0) s^{j_1+1} f^{(j_1)}(s) ds 
+ \sum_{j_2=0}^3 \int_r^{\infty} O(r^{-\delta} s^{-1 + \delta}\omega^0) s^{j_2} f^{(j_2)}(s) ds  \\
& + \sum_{j_3=1}^4 \int_r^{\infty} O(r^{-\delta} s^{-1 + \delta}\omega^0) s^{j_3-1} f^{(j_3)}(s) ds
\end{split}
\end{align}
for all $r > 0$ and all $\omega \gg 1$.
\end{lemma}

\begin{proof}

In the following, we denote by $\chi: \R^+ \to \R$ a smooth cut-off function
with $\chi(r) = 1$ for $r \leq \frac12$ and  $\chi(r) = 0$ for $r \geq 1$.
Let $c > 0$ be the constant from Lemma \ref{Le:FundSystemCenterHankel} (we always assume that $\omega \gg c^2$). We define 
\begin{align}\label{Def:Cutoffs}
\begin{split}
\chi_{B}(\omega^{\frac12} r) := \chi(c^{-1} \omega^{\frac12} r), \quad \chi_{H}(r,\omega) := \chi(r/8) -\chi(c^{-1} \omega^{\frac12} r) , \quad \chi_{W}(r) := 1- \chi(r/8)
\end{split}
\end{align}
and  $\chi_{W_1}(r) :=  \chi(r/40) - \chi(r/8)$, $\chi_{W_2}(r) :=  1 - \chi(r/40)$.
We write  $\mc T_1(\omega) =  \sum_{j=1}^7 T_{1,j}(\omega)$, where  
\[[T_{1,j}(\omega)f](r) :=  \int_r^{\infty} G_{1,j}(r,s,\omega) f(s) ds,\]
and 
\begin{align*}
& G_{1,1}(r,s,\omega)  := \chi_{B}(\omega^{\frac12} r) \chi_{B}(\omega^{\frac12} s)   \partial_r g_1(r,s,\omega), & \quad
& G_{1,2}(r,s,\omega)  :=  \chi_{B}(\omega^{\frac12} r)  \chi_{H}(s,\omega) \partial_r g_1(r,s,\omega), &  \\
& G_{1,3}(r,s,\omega)  := \chi_{B}(\omega^{\frac12} r)  \chi_W(s)  \partial_r g_1(r,s,\omega),& \quad 
& G_{1,4}(r,s,\omega)  :=  \chi_{H}(r,\omega) [1 - \chi_{W}(s)]  \partial_r g_1(r,s,\omega), & \\
& G_{1,5}(r,s,\omega)  :=  \chi_{H}(r,\omega) \chi_{W}(s) \partial_r g_1(r,s,\omega), & \quad
& G_{1,6}(r,s,\omega)  :=  \chi_{W_1}(r) \partial_r g_1(r,s,\omega) , &  \\
& G_{1,7}(r,s,\omega)  :=  \chi_{W_2}(r)  \partial_r g_1(r,s,\omega),&
\end{align*}
for $s \geq r > 0$. Recall that 
\[ g_1(r,s,\omega) =\frac{1}{W(\omega)} r^{-3} s^3 e^{\frac{1}{2}(r^2-s^2)}   v_0(r,\omega)v_-(s,\omega).\]
We consider the operators $T_{1,j}(\omega)$ individually and show that the corresponding weighted derivatives are of the form as stated in Lemma \ref{Le:ResEst_T1}.
For $j=1,\dots,6$, we only discuss the expressions $(\cdot)^{m}[T_{1,j}(\omega)f]^{(m)}$, for $m=0,\dots,3$. From this it will 
become obvious that $(\cdot)^{m+2}[T_{1,j}(\omega)f]^{(m)}$, $m=0,1$, can be written as in Eq.~\eqref{Eq:LeT1_2}.

%-------------T_{1,1}---------------------------
We start with $T_{1,1}(\omega)$. Here, we use the representation of $v_0$ and $v_-$ in terms of
perturbed Bessel functions. Let $\psi_0$ be as in  Lemma \ref{Le:FundSystemCenterBessel}.
By Taylor expansion, 
\[\psi_0(r,\omega) =  \alpha_0  r^3 \mu^{\frac54}[1  + \Oc(r^2\omega)],\]
for some constant $ \alpha_0 \in \C$  and all $0 < r \omega^{\frac12} \leq c$. Hence,  
\[\partial_r[e^{\frac{r^2}{2} }r^{-3} v_0(r,\omega)] = \Oc(r \omega^{\frac{9}{4}}).\]
In view of Eq.~\eqref{Eq:Rep_vmv0v1} and Eq.~\eqref{Eq:Wronskian} we infer that 
\begin{align}\label{Eq:BBKernel}
\begin{split}
G_{1,1}(r,s,\omega)  & =  \chi_{B}(\omega^{\frac12} r) \chi_{B}(\omega^{\frac12} s) 
\Oc(r \omega^{\frac{9}{4}}) [\Oc(s^3 \omega^{\frac{5}{4}}) + \Oc(s^{-2} \omega^{-\frac{5}{4}})] s^3  \\
&  = \chi_{B}(\omega^{\frac12} r) \chi_{B}(\omega^{\frac12} s)\Oc(r^{0} s^{0} \omega^{0} ).
\end{split}
\end{align}
Note that $\partial_r \chi_{B}(\omega^{\frac12} r) = \chi'_B(\omega^{\frac12} r) \omega^{\frac{1}{2}} = \Oc(r^{-1} \omega^0)$, since
 $\omega^{\frac12} \simeq r^{-1}$ on the support of $\chi'_B$. Thus, we can write $\chi_{B}(\omega^{\frac12} r) = \Oc(r^0 \omega^0)$. 
Consequently, the claim follows from
 \begin{align}\label{Eq:FormDerint}
\partial^m_r \int_r^{\infty} G_{1,1}(r,s,\omega) f(s) ds  & = \sum_{k=0}^{m-1} \Oc(r^{k-m+1})f^{(k)}(r) 
+ \int_r^{\infty} \Oc(r^{-m}s^0\omega^0) f(s) ds.
\end{align}

%-------------T_{1,2}---------------------------
For $T_{1,2}(\omega)$, we use Lemma \ref{Le:FundSystemCenterRepresentations} and Lemma \ref{Le:FundSystemCenterHankel} to see that 
\begin{align*}
  \chi_H(s,\omega)v_{-}(s,\omega)  & =   \chi_H(s,\omega) [ \Oc(\omega^0) e^{10\mu^{1/2}} \tilde v_{-}(s,\omega) + \Oc(\omega^{-\frac12}) e^{- 10\mu^{1/2}}  \tilde v_{+}(s,\omega) ] \\
  & = \Oc(\omega^{-\frac{1}{4}} s^{0}) e^{\mu^{1/2}(10-s)} + \Oc(\omega^{-\frac{3}{4}} s^{0})e^{-\mu^{1/2}(10 -s)}.
\end{align*}
Hence, 
\begin{align*}
 G_{1,2}(r,s,\omega)   = \chi_{B}(\omega^{\frac12} r) \chi_{H}(s,\omega)
 [\Oc( r^{0} s^{3} \omega ) e^{-\mu^{1/2}(20 -s)} + \Oc( r^{0} s^{3} \omega^{\frac32} ) e^{-\mu^{1/2}s}] . 
\end{align*}
We estimate
\[ | G_{1,2}(r,s,\omega)| \lesssim  \omega  e^{-12 \mathrm{Re} \mu^{1/2}} + (s \omega^{\frac12})^3 e^{-\mathrm{Re} \mu^{1/2}s} \lesssim 1\]
for all $0 < r \leq s$, since $\mathrm{Re} \mu^{1/2} \simeq \omega^{\frac12}$.
On the diagonal, i.e., for $s = r$, the support of $G_{1,2}$ is contained in the interval $[ \frac{1}{2} c \omega^{-\frac12}, c \omega^{-\frac12}]$.
There, $\omega^{\frac12} \simeq r^{-1}$, such that $e^{\pm \mu^{1/2}r} = \Oc(r^0 \omega^0)$ and 
$G_{1,2}(r,r,\omega) = \Oc(r^0 \omega^0).$ Hence,
\[ r  [T_{1,2}(\omega)f]'(r)   =  \Oc(r \omega^0) f(r)  +  \int_r^{\infty} r \partial_r G_{1,2}(r,s,\omega) f(s) ds.\]
By the same arguments as above, $| r \partial_r G_{1,2}(r,s,\omega)| \lesssim 1$.
 Since $[\partial^m_r  G_{1,2}(r,s,\omega)]_{s=r} = \Oc(r^{-m} \omega^{0})$, we 
get that 
\begin{align*}
r^m  [T_{1,2}(\omega)f]^{(m)}(r)  & = \sum_{j=0}^{m-1} \Oc(r^0 \omega^0)  r^j f^{(j)}(r)  +  \int_r^{\infty} r^{m} \partial^m_r G_{1,2}(r,s,\omega) f(s) ds,
\end{align*}
for $m=0,\dots 3$, and 
\[ |r^{m} \partial^m_r G_{1,2}(r,s,\omega) | \lesssim 1.\]
This implies the claim for $T_{1,2}(\omega)$. 

%-------------T_{1,3}---------------------------
For the next term we get that   
\begin{align*}
r^m [T_{1,3}(\omega)f]^{(m)}(r) = \int_r^{\infty} r^m \partial_r^m G_{1,3}(r,s,\omega) f(s) ds,
\end{align*}
since the kernel as well as its derivatives vanish for $s = r$. Furthermore,
\begin{align*}
r^m \partial_r^m  & G_{1,3}(r,s,\omega) = e^{-10 \mu^{1/2}} \Oc(r  \omega^{\frac{9}{4}}) \chi_W(s) e^{-\frac{s^2}{2}} s^3 v_{-}(s,\omega),  
\end{align*}
where $0 < r \leq  c \omega^{-\frac12}$. By Lemma \ref{Le:FundSysWeber}, 
$|v_{-}(s,\omega)| \lesssim \omega^{-\frac14} \langle \omega^{-\frac12} s \rangle^{-\frac12} e^{-\mathrm{Re} \mu \xi(s,\mu)}$,
for $s \geq 4$. We use Eq.~\eqref{Eq:RepXi_large_alt} together with
the fact that $\tilde \varphi(\cdot, \omega)$ is monotonically increasing and that $|\tilde \varphi(4, \omega)| \lesssim 1$, to estimate 
\begin{align*}
|r^m \partial_r^m  G_{1,3}(r,s,\omega)| & \lesssim  e^{- \mathrm{Re} \mu^{1/2}s } (\omega^{\frac12} s)^3 e^{-\frac{s^2}{2}} 
e^{-\tilde \varphi(s, \omega)}  \lesssim e^{-\frac{s^2}{2}}  e^{-\tilde \varphi(4, \omega)} s^{1-\delta} s^{-1 + \delta} \lesssim  r^{-\delta} s^{-1 + \delta} 
\end{align*}
for $m = 0, \dots, 3$, all $s \geq r > 0$ and $\omega \gg 1$.

%-------------T_{1,4}---------------------------
We now turn to the regime $\frac12 c \omega^{-\frac12} \leq r \leq s \leq 8$, where both $v_0$ and $v_{-}$ 
are represented in terms of Hankel functions. 
By Lemma \ref{Le:FundSystemCenterRepresentations}  and Lemma \ref{Le:FundSystemCenterHankel}, 
$v_0(r,\omega) = \Oc(r^{0} \omega^{-\frac{1}{4}}) e^{- \mu^{1/2} r} +  \Oc(r^{0} \omega^{-\frac{1}{4}}) e^{\mu^{1/2} r}$
such that 
\begin{align}\label{Eq:ResEst_DerHankel}
\begin{split}
\partial^j_r [ e^{\frac{r^2}{2}} r^{-3} v_0(r,\omega) ]  & = 
\Oc(r^{-3} \omega^{-\frac{1}{4} + \frac{j}{2}}) e^{- \mu^{1/2} r} +  \Oc(r^{-3} \omega^{-\frac{1}{4} + \frac{j}{2}}) e^{\mu^{1/2} r}.
\end{split}
\end{align}
We split the integral kernel into $G_{1,4} =  G^-_{1,4} + G^+_{1,4}$, where
\begin{align*}
 G^-_{1,4}(r,s,\omega) & =   \Oc(r^{-3} s^3 \omega^{-\frac12}) e^{-\mu^{1/2} r}  e^{-\mu^{1/2}(20 -s)} 
 +  \Oc(r^{-3} s^3 \omega^{-\frac12})  e^{-\mu^{1/2}(10 -r)}  e^{-\mu^{1/2}(10 -s)} \\ 
 & \phantom{!!!!!!!!!!!!!!!!!!!} +\Oc(r^{-3} s^3 \omega^0) e^{-\mu^{1/2} r}  e^{-\mu^{1/2} s}, 
\end{align*}
and $G^+_{1,4}(r,s,\omega)  =   \Oc(r^{-3} s^3 \omega^0)  e^{-\mu^{1/2}(s -r)}$. 
Accordingly, we write $T_{1,4}(\omega) = T^{-}_{1,4}(\omega) + T^{+}_{1,4}(\omega)$. 
We estimate
\begin{align}
| G^-_{1,4}(r,s,\omega)| \lesssim \omega e^{-12 \mathrm{Re} \mu^{1/2}} + \omega e^{-4 \mathrm{Re} \mu^{1/2}} 
+ (s\omega^{\frac12})^3 e^{- \mathrm{Re} \mu^{1/2}s} \lesssim 1 \lesssim  r^{-\delta} s^{-1 + \delta},
\end{align}
which shows that $T^{-}_{1,4}(\omega)$ is of the claimed form.
For derivatives of $G^-_{1,4}$, we get
\begin{align*}
\partial^j_r G^-_{1,4}(r,s,\omega) & =   \Oc(r^{-3} s^3 \omega^{-\frac12+\frac{j}{2}}) e^{-\mu^{1/2} r}  e^{-\mu^{1/2}(20 -s)} \\
& +  \Oc(r^{-3} s^3 \omega^{-\frac12+\frac{j}{2}})  e^{-\mu^{1/2}(10 -r)}  e^{-\mu^{1/2}(10 -s)} 
+\Oc(r^{-3} s^3 \omega^{\frac{j}{2}}) e^{-\mu^{1/2} r}  e^{-\mu^{1/2} s}.
\end{align*}
Since
\[G^-_{1,4}(r,r,\omega) = \Oc(r^0 \omega^{-\frac12}) e^{-20 \mu^{1/2}} + \Oc(r^0 \omega^{-\frac12})  e^{-2 \mu^{1/2}(10 -r)} + \Oc(r^0 \omega^0)  e^{-2 \mu^{1/2} r},\]
on the diagonal,  we infer that 
\[ r  [T^{-}_{1,4}(\omega)f]'(r)   =   d^{-}_0(r,\omega) f(r)  +  \int_r^{\infty} r \partial_r G^-_{1,4}(r,s,\omega) f(s) ds,\]
with $| d^{-}_0(r,\omega)| \lesssim 1$ and $|r \partial_r G^-_{1,4}(r,s,\omega)| \lesssim  1$ by the same arguments as above. 
By exploiting the exponential decay of the  diagonal terms, it is not hard to see that for $m=0,\dots,3$,
\[ r^m  [T^{-}_{1,4}(\omega)f]^{(m)}(r)   =   \sum_{k=0}^{m-1} d^-_{mk}(r,\omega) r^{k} f^{(k)}(r)  +  \int_r^{\infty} r^m \partial^m_r G^-_{1,4}(r,s,\omega) f(s) ds,\]
with  $|d^-_{mk}(r,\omega) | \lesssim 1$ and 
\[ |r^m \partial^m_r G^-_{1,4}(r,s,\omega)| \lesssim  \omega^{1+\frac{m}{2}} e^{-12 \mathrm{Re} \mu^{1/2}}
+ \omega^{1+\frac{m}{2}} e^{-4 \mathrm{Re} \mu^{1/2}} + (s\omega^{\frac12})^{3+\frac{m}{2}} e^{- \mathrm{Re} \mu^{1/2}s} \lesssim 1
\lesssim r^{-\delta} s^{-1 + \delta},\]
for all $s \geq r >0$ and $\omega \gg 1$.
These considerations imply also that the expressions $(\cdot)^{m+2}  [T^{-}_{1,4}(\omega)f]^{(m)}(r)$ are of the required form for $m=0,1$.  
Next, we consider $T^{+}_{1,4}(\omega)$ and estimate
\begin{align*}
 |G^+_{1,4}(r,s,\omega)| & \lesssim  r^{-3} s^3 e^{-\mathrm{Re} \mu^{1/2}(s -r)} 
\lesssim  (\mathrm{Re} \mu^{1/2} r)^{-3} (\mathrm{Re} \mu^{1/2} s)^3 e^{-\mathrm{Re} \mu^{1/2}(s -r)} \\
& \lesssim  e^{-\mathrm{Re} \mu^{1/2} s + 3 \log(\mathrm{Re} \mu^{1/2} s )} e^{\mathrm{Re} \mu^{1/2} r - 3 \log(\mathrm{Re} \mu^{1/2} r )} \lesssim 1,
\end{align*}
since $r\mapsto \mathrm{Re} \mu^{1/2} r - 3 \log(\mathrm{Re} \mu^{1/2} r )$ is monotonically increasing.
For the first derivative, we get  
\[r  [T^{+}_{1,4}(\omega)f]'(r) = \Oc(r^0 \omega^0) f(r) + \int_r^{\infty} \Oc(r^{-2} s^3 \omega^{\frac12})  e^{-\mu^{1/2}(s -r)}f(s)ds.\]
where the integrand is supported on $\frac12 c \omega^{-\frac12} \leq r \leq s \leq 8$. One integration by parts shows that 
\[ r  [T^{+}_{1,4}(\omega)f]'(r) = \Oc(r^0 \omega^0) f(r) + \int_r^{\infty}   e^{-\mu^{1/2}(s -r)}[\Oc(r^{-2} s^2 \omega^{0})f(s) + \Oc(r^{-2} s^2 \omega^{0}) s f'(s)] ds,\]
where $|r^{-2} s^2 e^{-\mu^{1/2}(s -r)}| \lesssim 1$ by the same reasoning as above. For the
second derivative, we use the above form of $[T^{+}_{1,4}(\omega)f]'$ and apply similar arguments. 
In summary, it is easy to check that 
\begin{align*}
r^m  [T^{+}_{1,4}(\omega)f]^{(m)}(r) = \sum_{j=0}^{m-1} \Oc(r^0 \omega^0) r^{j} f^{(j)}(r) 
+   \sum_{k=0}^{m} \int_r^{\infty} j_{mk}^{+}(r,s,\omega) s^{k} f^{(k)}(s) ds,
\end{align*}
 for $m=0,\dots,3$,  with 
$j_{mk}^{+}(r,s,\omega) =e^{-\mu^{1/2}(s -r)} \sum_{k=0}^{m} \Oc(r^{-{3+k}} s^{3-k} \omega^{0}).$
Thus,  
\[ |j_{mk}^{+}(r,s,\omega)| \lesssim 1 \lesssim r^{-\delta} s^{-1 + \delta},\]
for $m = 0,\dots,3$, $k = 0, \dots,m $, all $s \geq r >0$ and all $\omega \gg c^2$. 

%-------------T_{1,5}---------------------------
Next, we study $T_{1,5}(\omega)$, where the corresponding integral kernel is supported on $\frac12 c \omega^{-\frac12} \leq r \leq 8$, $s \geq 4$ and $s \geq r$. 
Recall that $v_{-}(s,\omega) = \Oc(\omega^{-\frac{1}{4}} s^0)\Oc(\langle \omega^{-\frac12} s \rangle^{-\frac12}) e^{-\mu \xi(s,\mu)}$. Hence,
\begin{align*}
 G_{1,5}(r,s,\omega) =  & \chi_{H}(r,\omega) \chi_{W}(s)e^{- 10 \mu^{1/2}}
[ \Oc(r^{-3} s^3 \omega^0)  e^{\mu^{1/2} r} + \Oc(r^{-3} s^3 \omega^0)  e^{-\mu^{1/2} r} ] \\
& \times e^{-\mu \xi(s,\mu)} e^{-\frac{s^2}{2}} \Oc(\langle \omega^{-\frac12} s \rangle^{-\frac12}).
\end{align*}

In view of Eq.~\eqref{Eq:RepXi_large_alt}, we can use the same arguments as above to estimate
\begin{align*}
| G_{1,5}(r,s,\omega)|  & \lesssim r^{-3} s^3 e^{-\frac{s^2}{2}} (e^{-\mathrm{Re} \mu^{1/2} s} e^{-\mathrm{Re} \mu^{1/2} r}   + e^{-\mathrm{Re} \mu^{1/2} (s-r)}) \\
& \lesssim
e^{-\frac{s^2}{2}}  (  e^{-\mathrm{Re} \mu^{1/2} s} \omega^{\frac32} s^3 +  r^{-3} s^3  e^{-\mathrm{Re} \mu^{1/2} (s-r)})
\lesssim e^{-\frac{s^2}{2}} \lesssim s^{-1 + \delta} r^{-\delta} .
\end{align*}
For the derivative of $[T_{1,5}(\omega)f]$, we use that  on the diagonal
\[G_{1,5}(r,r,\omega) =  \chi_{H}(r,\omega) \chi_{W}(r) [\Oc(r^0 \omega^0) e^{- 10 \mu^{1/2}} e^{\mu^{1/2} r -\mu \xi(r,\mu)} 
+ \Oc(r^0 \omega^0) e^{- 10 \mu^{1/2}}  e^{-(\mu^{1/2} r + \mu \xi(r,\mu))} ],
 \]
which implies that $|r G_{1,5}(r,r,\omega)| \lesssim 1$. We use that 
\[(s +\mu \partial_s \xi(s,\mu) )^{-1}  =  (s + (\mu + s^2)^{\frac12})^{-1} = \mu^{-\frac12} \Oc(\langle \omega^{-\frac12} s \rangle^{-1}),\] 
see Eq.~\eqref{Eq:EstLGPot}, such that 
\[e^{-\tfrac{s^2}{2} - \mu \xi(s,\mu)} = \left (\partial_s e^{-\tfrac{s^2}{2} - \mu \xi(s,\mu)} \right )  \mu^{-\frac12} \Oc(\langle \omega^{-\frac12} s \rangle^{-1}).\]
Integration by parts yields
\begin{align}\label{Eq:ResEst_Eq8}
\begin{split}
\int_r^{\infty} &  \chi_W(s)  e^{-\tfrac{s^2}{2}  - \mu \xi(s,\mu)} \Oc( \langle \omega^{-\frac12} s \rangle^{-\frac12}) \Oc(s^3 \omega^{0} ) f(s) ds \\
&   = \chi_W(r) e^{-\tfrac{r^2}{2}  - \mu \xi(r,\mu)}  \Oc(r^3 \omega^{-\frac{1}{2}}) \Oc(\langle \omega^{-\frac12} r \rangle^{-\frac{3}{2}}) f(r) \\
& + \sum_{k \in \{0,1\}} \int_r^{\infty}   e^{-\tfrac{s^2}{2}  - \mu \xi(s,\mu)} 
\Oc(\langle \omega^{-\frac12} s \rangle^{-\frac{3}{2}}) \Oc(s^2 \omega^{-\frac{1}{2}}) s^k f^{(k)}(s) ds,
\end{split}
\end{align}
and with this, it is easy to see that 
\begin{align*}
 r  [T_{1,5}(\omega)f]'(r) =  d(r,\omega) f(r) + \sum_{k \in \{0,1\}} \int_r^{\infty}  j_k(r,s,\omega) s^{k} f^{(k)}(s) ds,
\end{align*}
where 
\[d(r,\omega) =  \chi_{H}(r,\omega) \chi_{W}(r) [\Oc(r^0 \omega^0) e^{- 10 \mu^{1/2}} e^{\mu^{1/2} r -\mu \xi(r,\mu)} 
+ \Oc(r^0 \omega^0) e^{- 10 \mu^{1/2}}  e^{-(\mu^{1/2} r + \mu \xi(r,\mu))} ],
 \]
and 
\[ j_k(r,s,\omega)  = e^{- 10 \mu^{1/2}}
[ \Oc(r^{-2} s^2 \omega^0)  e^{\mu^{1/2} r} + \Oc(r^{-2} s^2 \omega^0)  e^{-\mu^{1/2} r} ]  e^{-\tfrac{s^2}{2} 
- \mu \xi(s,\mu)}\Oc(\langle \omega^{-\frac12} s \rangle^{-\frac32}).\]
Thus, we infer that $|d(r,\omega) | \lesssim 1$ and 
$|j_k(r,s,\omega)| \lesssim s^{-1 + \delta} r^{-\delta}$. To control derivatives of the boundary terms we exploit that for $4 \leq r \leq 8$,
\begin{align} \label{Der:ExpHankelPlusWeberMin}
\partial_r e^{\mu^{1/2} r -\mu \xi(r,\mu)} = e^{\mu^{1/2} r -\mu \xi(r,\mu)} \mu^{\frac12}
[1 - (1+\mu^{-1} r^2)^{\frac12}] = \Oc(r^0 \omega^{-\frac12}) e^{\mu^{1/2} r -\mu \xi(r,\mu)}.
\end{align}
Based on the above representation of $ [T_{1,5}(\omega)f]'$, repeated integration by parts shows that
for $m=0, \dots 3$, 
\begin{align*}
r^m  [T_{1,5}(\omega)f]^{(m)}(r) = \sum_{j=0}^{m-1} \tilde d_{mj}(r,\omega)  r^{j} f^{(j)}(r) 
+  \sum_{k=0}^{m}  \int_r^{\infty} \tilde j_{mk}(r,s,\omega)  s^{k} f^{(k)}(s) ds,
\end{align*}
with $|\tilde d_{mj}(r,\omega)| \lesssim 1$ for $j = 0, \dots, m-1$, $|\tilde j_{mk}(r,s,\omega)  | \lesssim  s^{-1 + \delta} r^{-\delta}$ for $k = 0,\dots,m$, 
all $s \geq r > 0$ and $\omega \gg c^2$. 

%-------------T_{1,6}---------------------------
For $T_{1,6}(\omega)$, i.e., the regime $4 \leq r \leq 40$, $s \geq 4$ and $s \geq r$, we write $v_0$ again in terms of Hankel functions.
This is similar to the previous case and one immediately gets that 
\begin{align*}
| G_{1,6}(r,s,\omega)|  & \lesssim r^{-3} s^3 e^{-\frac{s^2}{2}} (e^{-\mathrm{Re} \mu^{1/2} s} e^{-\mathrm{Re} \mu^{1/2} r}   + e^{-\mathrm{Re} \mu^{1/2} (s-r)}) 
\lesssim s^3 e^{-\frac{s^2}{2}}  \lesssim s^{-1 + \delta} r^{-\delta} .
\end{align*}
Higher derivatives can be handled using integration by parts. This is analogous to $T_{1,5}(\omega)$ and we omit further details.

%-------------T_{1,7}---------------------------
It is left to discuss $T_{1,7}(\omega)$, i.e., the regime  $20 \leq r \leq s$. 
In view of Lemma \ref{Le:FundSysWeber}, we get that
\begin{align}\label{Eq:Der_Expvpm}
\partial^m_r[r^{-3} e^{\frac12 r^2} v_{\pm}(r,\omega)]= \Oc(r^{-3} \omega^{-\frac{1}{4} + \frac{m}{2}})
\Oc(\langle \omega^{-\frac12} r \rangle^{-\frac12 \pm m}) e^{\frac12 r^2\pm \mu \xi(r,\mu)} 
\end{align}
for $m \in \N_0$, since $r \pm \mu \partial_r  \xi(r,\mu) = \Oc(\omega^{\frac12} \langle \omega^{-\frac12} r \rangle^{\pm 1}).$
An application of Lemma \ref{Le:Rep_vmv0} yields
\begin{align*}
\frac{\chi_{W_2}(r)}{W(\omega)} \partial_r [r^{-3} e^{\frac{r^2}{2}}v_0(r,\omega)]
 & =  \Oc(r^{-3} \omega^{-\frac14}) \Oc(\langle \omega^{-\frac12} r \rangle^{-\frac{3}{2}}) e^{\frac12 r^2 -  \mu \xi(r,\mu)} \\
 & + \Oc(r^{-3} \omega^{\frac14}) \Oc(\langle \omega^{-\frac12} r \rangle^{\frac12}) e^{\frac12 r^2 +  \mu \xi(r,\mu)},
\end{align*}
such that 
\begin{align*}
G_{1,7}(r,s,\omega) & = \Oc(r^{-3} s^3 \omega^{-\frac12})  e^{\frac12 r^2 -  \mu \xi(r,\mu)}  e^{-\frac12 s^2 - \mu \xi(s,\mu)} \Oc(\langle \omega^{-\frac12} r \rangle^{-\frac{3}{2}}) 
\Oc(\langle \omega^{-\frac12} s \rangle^{-\frac{1}{2}}) \\
& +  \Oc(r^{-3} s^3 \omega^{0})  e^{\frac12 r^2 +  \mu \xi(r,\mu)}  e^{-\frac12 s^2 - \mu \xi(s,\mu))} \Oc(\langle \omega^{-\frac12} r \rangle^{\frac{1}{2}}) 
\Oc(\langle \omega^{-\frac12} s \rangle^{-\frac{1}{2}}).
\end{align*}
An integration by parts, cf.~ Eq.~\eqref{Eq:ResEst_Eq8}, shows that   
\begin{align}\label{Eq:ResEst_T17}
[T_{1,7}(\omega)f](r) = p(r,\omega) r f(r) +  \sum_{k \in \{0,1\}} \int_r^{\infty} h_k(r,s,\omega) s^{k+1} f^{(k)}(s) ds, 
\end{align}
where 
\begin{align*}
 p(r,\omega)  & = \Oc( r^{-1} \omega^{-1}) \Oc(\langle \omega^{-\frac12} r \rangle^{-3}) e^{- 2 \mu \xi(r,\mu)}
 + \Oc( r^{-1} \omega^{-\frac12})  \Oc(\langle \omega^{-\frac12} r \rangle^{-1}) \\
 & =  \Oc( r^{-2} \omega^{-\frac12}) e^{- 2 \mu \xi(r,\mu)} +  \Oc(r^{-2} \omega^0),
\end{align*}
and 
\begin{align*}
h_k(r,s,\omega)  &  = \Oc(r^{-3} s \omega^{-1})  e^{\frac12 r^2 -  \mu \xi(r,\mu)}  e^{-\frac12 s^2 - \mu \xi(s,\mu)} \Oc(\langle \omega^{-\frac12} r \rangle^{-\frac{3}{2}}) 
\Oc(\langle \omega^{-\frac12} s \rangle^{-\frac{3}{2}}) \\
& +  \Oc(r^{-3} s \omega^{-\frac12})  e^{\frac12 r^2 +  \mu \xi(r,\mu)}  e^{-\frac12 s^2 - \mu \xi(s,\mu)} \Oc(\langle \omega^{-\frac12} r \rangle^{\frac{1}{2}}) 
\Oc(\langle \omega^{-\frac12} s \rangle^{-\frac{3}{2}})
\end{align*}
for $k=1,2$. We consider the expression  $(\cdot)^2[T_{1,7}(\omega)f]$. 
First, we observe that $| r^2 p(r,\omega)| \lesssim 1$, since by Eq.~\eqref{Eq:RepXi_large_alt}
\begin{align}\label{Eq:ResESt_T17Exp}
\begin{split}
 e^{- 2 \mathrm{Re} \mu \xi(r,\mu)} & = e^{-2 \mathrm{Re} \mu^{1/2} (r - 10)} e ^{-2\tilde \varphi(r,\omega)} 
\leq e^{-2 \mathrm{Re} \mu^{1/2} (r - 10)} e ^{-2\tilde \varphi(10,\omega)} \\
& \lesssim e^{-2 \mathrm{Re} \mu^{1/2} (r - 10)} \lesssim e^{-20 \mathrm{Re} \mu^{1/2}}
\end{split}
\end{align}
for all $r \geq 20$. For the derivative we get
\[\partial_r p(r,\omega) =  [\Oc( r^{-3} \omega^{-\frac12}) + \Oc( r^{-2} \omega^{-\frac12}) \sqrt{\mu + r^2}] e^{- 2 \mu \xi(r,\mu)} + \Oc(r^{-3} \omega^0) .\]
For $20 \leq r < \omega^{\frac12}$, 
\[ |\Oc( r \omega^{-\frac12}) \sqrt{\mu + r^2} e^{- 2 \mu \xi(r,\mu)}| \lesssim \omega^{\frac12} 
e^{- 2 \mathrm{Re} \mu \xi(r,\mu)} \lesssim \omega^{\frac12} e^{-20 \mathrm{Re} \mu^{1/2}} \lesssim 1.\]
For $r \geq \omega^{\frac12}$, we use Eq.~\eqref{Eq:RepXi_large}. 
In fact, since $b > -4$, we have the bound $\langle \omega^{-\frac12} r \rangle^{-\frac{b}{2}} \lesssim r^2$.
With this, we infer that 
\[ |\Oc( r \omega^{-\frac12}) \sqrt{\mu + r^2} e^{- 2 \mu \xi(r,\mu)}| \lesssim  r^2 e^{- \frac{r^2}{2}} \langle \omega^{-\frac12} r \rangle^{-\frac{b}{2}} e^{- \varphi(r,\omega)} 
\lesssim  e^{- \frac{r^2}{2}}  r^{4} e^{- \varphi(10,\omega) } \lesssim 1,\]
since $\varphi(\cdot,\omega)$ is monotonically increasing for all $r \geq 3$ and $|\varphi(10,\omega)| \lesssim 1$. 
Hence, $|r^3 \partial_r p(r,\omega)| \lesssim 1$. In a similar manner one can show that $|r^3 \partial^2_r p(r,\omega)| \lesssim 1$.

To estimate the integral kernels, which are supported for $s \geq r \geq 20$,
we use the bound$\langle \omega^{-\frac12} r \rangle^{-\frac{3}{2}- \frac{b}{2}}  \lesssim r^{\frac12}$ to infer that 
%$r s^{-1} \lesssim \langle \omega^{-\frac12} r \rangle \langle \omega^{-\frac12} s \rangle^{-1}$, to obtain
\begin{align}\label{Eq:ResEst_KernelT17}
\begin{split}
| r^2 h_k(r,s,\omega) |  &  \lesssim 
\omega^{-1} r^{-1} s \langle \omega^{-\frac12} r \rangle^{-\frac{3}{2}- \frac{b}{2}} 
\langle \omega^{-\frac12} s \rangle^{-\frac{3}{2}-\frac{b}{2}} e^{-s^2} e^{-\varphi(s,\omega)} e^{-\varphi(r,\omega)} \\
& + \omega^{-\frac12 } r^{-1} s  \langle \omega^{-\frac12} r \rangle^{\frac{1}{2} + \frac{b}{2}} \langle \omega^{-\frac12} s \rangle^{-\frac{3}{2}-\frac{b}{2}}  
e^{-(\varphi(s,\omega) - \varphi(r,\omega))} e^{-(s^2-r^2)}  \\
& \lesssim \omega^{-1} r^{-\frac12} s^{\frac{3}{2}} e^{-s^2}   + \omega^{-\frac12}  r^{-2} s^2 e^{-(s^2-r^2)} 
\langle \omega^{-\frac12} r \rangle^{\frac{3}{2} + \frac{b}{2}} \langle \omega^{-\frac12} s \rangle^{-\frac{5}{2}-\frac{b}{2}} 
\end{split}
\end{align}
for $k=0,1$. We readily estimate $ \omega^{-1} r^{-\frac12} s^{\frac{3}{2}} e^{-s^2} \lesssim  s^{-1}\lesssim r^{-\delta}s^{-1+\delta}$. 
For the second term we first consider the case $b>-3$ which yields 
\[  \omega^{-\frac12}  r^{-2} s^2 e^{-(s^2-r^2)} 
\langle \omega^{-\frac12} r \rangle^{\frac{3}{2} + \frac{b}{2}} \langle \omega^{-\frac12} s \rangle^{-\frac{5}{2}-\frac{b}{2}} \lesssim \omega^{-\frac12}  \langle \omega^{-\frac12} s \rangle^{-1}  \lesssim s^{-1}. \]
If $b \in (-4,-3]$, we distinguish three cases. For $20 \leq r \leq s < \omega^{\frac12}$ we use that $\omega^{-\frac12} s <1$ to get 
\[  \omega^{-\frac12}  r^{-2} s^2 e^{-(s^2-r^2)} 
\langle \omega^{-\frac12} r \rangle^{\frac{3}{2} + \frac{b}{2}} \langle \omega^{-\frac12} s \rangle^{-\frac{5}{2}-\frac{b}{2}}  \lesssim s^{-1} \lesssim  r^{\frac32 +\frac{b}{2}}s^{-\frac{5}{2} - \frac{b}{2}},\] 
since $r^{-\frac32 -\frac{b}{2}} \leq s^{-\frac32 -\frac{b}{2}}$ 
for $-4 < b \leq -3$.
For $20 \leq r <  \omega^{\frac12} \leq s$, 
\begin{align*}
 \omega^{-\frac12}  r^{-2} s^2 e^{-(s^2-r^2)} 
\langle \omega^{-\frac12} r \rangle^{\frac{3}{2} + \frac{b}{2}} \langle \omega^{-\frac12} s \rangle^{-\frac{5}{2}-\frac{b}{2}} 
\lesssim   r^{-3} s^3 e^{-(s^2-r^2)} s^{-1} \lesssim s^{-1} \lesssim r^{\frac32 +\frac{b}{2}}s^{-\frac{5}{2} - \frac{b}{2}}. 
\end{align*}
Finally, for $\omega^{\frac12} \leq r \leq s$, 
\begin{align*}
 \omega^{-\frac12}  r^{-2} s^2 e^{-(s^2-r^2)} 
\langle \omega^{-\frac12} r \rangle^{\frac{3}{2} + \frac{b}{2}} \langle \omega^{-\frac12} s \rangle^{-\frac{5}{2}-\frac{b}{2}} 
\lesssim   r^{\frac{3}{2} + \frac{b}{2}} s^{-\frac{5}{2}-\frac{b}{2}}.
\end{align*}

We conclude that 
\[| r^2 h_k(r,s,\omega)| \lesssim r^{-\delta} s^{-1 + \delta},\]
for $k=1,2$, all $s \geq r > 0$ and $\omega \gg c^2$.   
This proves that $(\cdot)^{2}[ T_{1,7}(\omega)f]$ is of the required form. Next, we consider $(\cdot)^{3}[ T_{1,7}(\omega)f]'$, starting from Eq.~\eqref{Eq:ResEst_T17}.
Evaluation of the integral kernels yields 
\[ h_k(r,r,\omega) = \Oc( r^{-3} \omega^{-\frac12}) e^{- 2 \mu \xi(r,\mu)} +  \Oc(r^{-3} \omega^0), \]
such that  Eq.~\eqref{Eq:ResESt_T17Exp} yields $|r^{3} h_k(r,r,\omega)| \lesssim 1$ for $k=1,2$. 
As in Eq.~\eqref{Eq:Der_Expvpm},
\begin{align*}
\partial_r h_k(r,s,\omega)  &  = \Oc(r^{-3} s \omega^{-\frac12})  e^{\frac12 r^2 -  \mu \xi(r,\mu)}  e^{-\frac12 s^2 - \mu \xi(s,\mu)} \Oc(\langle \omega^{-\frac12} r \rangle^{-\frac{5}{2}}) 
\Oc(\langle \omega^{-\frac12} s \rangle^{-\frac{3}{2}}) \\
& +  \Oc(r^{-3} s \omega^{0})  e^{\frac12 r^2 +  \mu \xi(r,\mu)}  e^{-\frac12 s^2 - \mu \xi(s,\mu)} \Oc(\langle \omega^{-\frac12} r \rangle^{\frac{3}{2}}) 
\Oc(\langle \omega^{-\frac12} s \rangle^{-\frac{3}{2}}).
\end{align*}
An integration by parts, cf.~\eqref{Eq:ResEst_Eq8}, shows that 
\begin{align}\label{Eq:ResEst_T17prime}
[T_{1,7}(\omega)f]'(r) = p_1(r,\omega) r f(r) + p_2(r,\omega) r^2 f'(r) + \sum_{k=0}^{2} \int_r^{\infty}  \tilde h_{k}(r,s,\omega) s^{k+1} f^{(k)}(s) ds,
\end{align}
with $p_{2}(r,\omega) = \Oc( r^{-3} \omega^{-\frac12}) e^{- 2 \mu \xi(r,\mu)} +  \Oc(r^{-3} \omega^0)$,
\begin{align*}
p_{1}(r,\omega) = \partial_r p(r,\omega) + \Oc( r^{-3} \omega^{-\frac12}) e^{- 2 \mu \xi(r,\mu)} +  \Oc(r^{-3} \omega^0).
\end{align*}
Hence, $|r^3 p_j(r,\omega)| \lesssim 1$, for $j=1,2$. The integral kernels are of the form 
\begin{align}\label{Eq:EstT17_Eq1}
\tilde  h_k(r,s,\omega)  &  = \Oc(r^{-3} s^0 \omega^{-1})  e^{\frac12 r^2 -  \mu \xi(r,\mu)}  e^{-\frac12 s^2 - \mu \xi(s,\mu)} \Oc(\langle \omega^{-\frac12} r \rangle^{-\frac{5}{2}}) 
\Oc(\langle \omega^{-\frac12} s \rangle^{-\frac{5}{2}}) \\
& +  \Oc(r^{-3} s^0 \omega^{-\frac12})  e^{\frac12 r^2 +  \mu \xi(r,\mu)}  e^{-\frac12 s^2 - \mu \xi(s,\mu)} \Oc(\langle \omega^{-\frac12} r \rangle^{\frac{3}{2}}) 
\Oc(\langle \omega^{-\frac12} s \rangle^{-\frac{5}{2}}),
\end{align} 
such that 
\begin{align*}
| r^3 \tilde  h_k(r,s,\omega) |  &  \lesssim 
\omega^{-1}  \langle \omega^{-\frac12} r \rangle^{-\frac{5}{2}- \frac{b}{2}} 
\langle \omega^{-\frac12} s \rangle^{-\frac{5}{2}-\frac{b}{2}} e^{-s^2} e^{-\varphi(s,\omega)} e^{-\varphi(r,\omega)} \\
& + \omega^{-\frac12 }   \langle \omega^{-\frac12} r \rangle^{\frac{3}{2} + \frac{b}{2}} \langle \omega^{-\frac12} s \rangle^{-\frac{5}{2}-\frac{b}{2}}  
e^{-(\varphi(s,\omega) - \varphi(r,\omega))} e^{-(s^2-r^2)}  \\
& \lesssim \omega^{-1}  e^{-s^2}   + \omega^{-\frac12} e^{-(s^2-r^2)} 
\langle \omega^{-\frac12} r \rangle^{\frac{3}{2} + \frac{b}{2}} \langle \omega^{-\frac12} s \rangle^{-\frac{5}{2}-\frac{b}{2}}  \lesssim r^{-\delta} s^{-1 + \delta},
\end{align*}
by applying the same arguments as above. This proves that the expressions
$(\cdot)^{m+2} [T_{1,7}(\omega)f]^{(m)}$ can be written as claimed for $m=0,1$.
From this, we obtain the respective statement for $(\cdot)^{m} [T_{1,7}(\omega)f]^{(m)}$, $m=0,1$.
It is left to consider $(\cdot)^{m} [T_{1,7}(\omega)f]^{(m)}$ for $m = 2,3$.
From Eq.~\eqref{Eq:EstT17_Eq1} we get that  
\begin{align*}
\partial_r \tilde  h_k(r,s,\omega)  &  = \Oc(r^{-3} s^0 \omega^{-\frac12})  e^{\frac12 r^2 -  \mu \xi(r,\mu)}  e^{-\frac12 s^2 - \mu \xi(s,\mu))}
\Oc(\langle \omega^{-\frac12} r \rangle^{-\frac{7}{2}}) 
\Oc(\langle \omega^{-\frac12} s \rangle^{-\frac{5}{2}}) \\
& +  \Oc(r^{-3} s^0 \omega^{0})  e^{\frac12 r^2 +  \mu \xi(r,\mu)}  e^{-\frac12 s^2 - \mu \xi(s,\mu))} \Oc(\langle \omega^{-\frac12} r \rangle^{\frac{5}{2}}) 
\Oc(\langle \omega^{-\frac12} s \rangle^{-\frac{5}{2}}).
\end{align*} 
Integrating by parts shows that  
\begin{align}\label{Eq:ResEst_T17dprime}
[T_{1,7}(\omega)f]''(r) =  \sum_{k=0}^{2} q_k(r,\omega) r^{k}f^{(k)}(r)  +  \sum_{j=0}^{3} \int_r^{\infty} k_{j}(r,s,\omega) s^{j} f^{(j)}(s) ds 
\end{align}
with $|r^2 q_k(r,\omega) | \lesssim 1$, since for example 
\[ \tilde  h_k(r,r,\omega) = \ \Oc( r^{-4} \omega^{-\frac12}) e^{- 2 \mu \xi(r,\mu)} +  \Oc(r^{-4} \omega^0),\]
such that $|r^3 \tilde  h_k(r,r,\omega)| \lesssim 1$ for $k=1,2,3$.  Furthermore, $| r^2 p_j(r,\omega)| \lesssim 1$ and $| r^3 \partial_r p_j(r,\omega)| \lesssim 1$ for $j=1,2$. 
For the integral kernels we have 
\begin{align*}
k_{j}(r,s,\omega) &  = \Oc(r^{-3} s^0 \omega^{-1})  e^{\frac12 r^2 -  \mu \xi(r,\mu)}  e^{-\frac12 s^2 - \mu \xi(s,\mu)} \Oc(\langle \omega^{-\frac12} r \rangle^{-\frac{7}{2}}) 
\Oc(\langle \omega^{-\frac12} s \rangle^{-\frac{7}{2}}) \\
& +  \Oc(r^{-3} s^0 \omega^{-\frac12})  e^{\frac12 r^2 +  \mu \xi(r,\mu)}  e^{-\frac12 s^2 - \mu \xi(s,\mu)} \Oc(\langle \omega^{-\frac12} r \rangle^{\frac{5}{2}}) 
\Oc(\langle \omega^{-\frac12} s \rangle^{-\frac{7}{2}}),
\end{align*} 
and by the same arguments as above, 
\[ | r^2 k_{j}(r,s,\omega) | \lesssim r^{-1} \omega^{-1}  e^{-s^2}  +
 r^{-1} \omega^{-\frac12 }  \langle \omega^{-\frac12} r \rangle^{\frac{5}{2} + \frac{b}{2}} \langle \omega^{-\frac12} s \rangle^{-\frac{7}{2}-\frac{b}{2}}  e^{-(s^2-r^2)}  
\lesssim  r^{-\delta} s^{-1 + \delta}. \] 
From this, it is obvious that taking the derivative of Eq.~\eqref{Eq:ResEst_T17dprime} and performing another integration by parts yields
\[ (\cdot)^3 [T_{1,7}(\omega)f]''' =  \sum_{k=0}^{3} \tilde q_k(r,\omega) r^{k}f^{(k)}(r)  + \int_r^{\infty} \tilde k_{0}(r,s,\omega) f(s) ds + 
\sum_{j=1}^{4}  \int_r^{\infty} \tilde k_{j}(r,s,\omega) s^{j-1} f^{(j)}(s) ds, \]
with $|\tilde q_k(r,\omega)| \lesssim 1$ and $| \tilde k_{j}(r,s,\omega)| \lesssim r^{-\delta} s^{-1 + \delta}$ for $j = 0,\dots 4$.

Since $\mc T_1(\omega) =  \sum_{j=1}^7 T_{1,j}(\omega)$, we sum up the individual contributions
to see that Eq.~\eqref{Eq:LeT1_1} and Eq.~\eqref{Eq:LeT1_2} hold.
\end{proof}

\subsection{Representation of $\mc T_2(\omega)$}\label{AppT2}

\begin{lemma}\label{Le:ResEst_T2}
Let $\mu = b + \I \omega$, with $b > -4$ fixed. Define $\delta \in [0,\frac12)$ as in Eq.~\eqref{Def:delta}.
Let $f \in C^{\infty}_{\mathrm{e},0}(\R)$. Then, for $m=0,1$, 
\begin{align}\label{Eq:LeT2_1}
\begin{split}
 r^{m+2} [\mc T_2(\omega)f]^{(m)}(r)     & =   O(r \omega^0)  f(r)  +  O(r^2 \omega^0)  f'(r) \\
& +\sum_{k_1=0}^{2} \int_0^r  O( r^{-1 + \delta} s^{-\delta} \omega^0) s^{k_1+1} f^{(k_1)}(s) ds \\
& + \sum_{k_2=0}^{3} \int_0^r O( r^{-1 + \delta} s^{-\delta} \omega^0)s^{k_2} f^{(k_2)}(s) ds
\end{split}
\end{align}
for all $r > 0$ and $\omega \gg 1$. Furthermore, for $m=0,\dots,3$,
\begin{align}\label{Eq:LeT2_2}
\begin{split}
 r^m [\mc T_2(\omega)f]^{(m)}(r)  & =  \sum_{j=0}^{2}  O(r^{j+1}  \omega^0)  f^{(j)}(r) + \sum_{k=0}^{3} 
 O( r^{k}   \omega^0) f^{(k)}(r) \\
& + \sum_{j_1=0}^2 \int_0^r O( r^{-1 + \delta} s^{-\delta} \omega^0) s^{j_1+1} f^{(j_1)}(s) ds 
+ \sum_{j_2=0}^3 \int_0^r O( r^{-1 + \delta} s^{-\delta} \omega^0)  s^{j_2} f^{(j_2)}(s) ds  \\
& + \sum_{j_3=1}^4 \int_0^r O( r^{-1 + \delta} s^{-\delta} \omega^0) s^{j_3-1} f^{(j_3)}(s) ds
\end{split}
\end{align}
for all $r > 0$ and $\omega \gg 1$.
\end{lemma}

\begin{proof}
With the notation of Eq.~\eqref{Def:Cutoffs}  
we write $\mc T_2(\omega)$ as $\mc T_2(\omega)=  \sum_{j=1}^7 T_{2,j}(\omega)$, where  
\[[T_{2,j}(\omega)f](r) :=  \int_0^{r} G_{2,j}(r,s,\omega) f(s) ds,\]
and 
\begin{align*}
  & G_{2,1}(r,s,\omega)  := \chi_{B}(\omega^{\frac12} r)  \partial_r g_2(r,s,\omega), &  \quad 
 & G_{2,2}(r,s,\omega)  := \chi_{H}(r,\omega)\chi_{B}(\omega^{\frac12} s)\partial_r g_2(r,s,\omega), &    \\
 & G_{2,3}(r,s,\omega)  := \chi_{H}(r,\omega)[ 1 -\chi_{B}(\omega^{\frac12} s)]\partial_r g_2(r,s,\omega),  &
\quad   & G_{2,4}(r,s,\omega)  :=   \chi_{W}(r)\chi_{B}(\omega^{\frac12} s) \partial_r g_2(r,s,\omega), &  \\
 &  G_{2,5}(r,s,\omega)  :=  \chi_{W}(r) \chi_{H}(s,\omega)  \partial_r g_2(r,s,\omega),  & 
\quad  &  G_{2,6}(r,s,\omega)  :=  \chi_{W_1}(r) \chi_W(s) \partial_r g_2(r,s,\omega), &  \\
& G_{2,7}(r,s,\omega)  := \chi_{W_2}(r) \chi_W(s)  \partial_r g_2(r,s,\omega), & 
\end{align*} 
for $0 < s \leq r$ and 
\[ g_2(r,s,\omega) = \frac{1}{W(\omega)} r^{-3} s^3 e^{\frac{1}{2}(r^2-s^2)}  v_-(r,\omega)v_0(s,\omega). \]
For $j=1,2,3$, it suffices to consider the expressions $(\cdot)^{(m)}[T_{2,j}(\omega)f]^{(m)}$, where $m=0,\dots,3$.

%-----------------T_{2,1}----------------------------------------------------
We start with $T_{2,1}(\omega)$ corresponding to the regime  $0 < s \leq r \leq c \omega^{-\frac12}$. As in the proof of Lemma \ref{Le:ResEst_T1}, we use
Lemma \ref{Le:Rep_vmv0} and Lemma \ref{Le:FundSystemCenterBessel} to see that 
\[ \frac{ \chi_{B}(\omega^{\frac12} r)}{W(\omega)}\partial_r [ e^{\frac{r^2}{2}} r^{-3} v_{-}(r,\omega)] = \Oc(r \omega^{\frac{9}{4}}) + \Oc(r^{-6} \omega^{-\frac{5}{4}}).\]
Since $v_0(s,\omega) = \Oc(s^{3} \omega^{\frac{5}{4}})$, we get 
\[G_{2,1}(r,s,\omega) = \Oc(r s^6 \omega^{\frac{7}{2}}) + \Oc(r^{-6} s^6 \omega^{0}) = \Oc(r^0 s^0 \omega^0).\]
This implies 
\begin{align*}
r^m [T_{2,1}(\omega)f]^{(m)}(r)  & =  \sum_{j=0}^{m-1} \Oc(r^{0}\omega^0) r^j f^{(j)}(r)  + \int_0^{r} \Oc(r^0 s^0 \omega^0) f(s)ds,
\end{align*}
cf.~\eqref{Eq:FormDerint}.
%-----------------T_{2,2}----------------------------------------------------
For $\frac12 c \omega^{-\frac12}  \leq  r \leq  8$, we have 
\[ \frac{1}{W(\omega)}\partial_r [ e^{\frac{r^2}{2}} r^{-3} v_{-}(r,\omega)] = 
\Oc(r^{-3}  \omega^{\frac{1}{4}}) e^{-\mu^{1/2} r } + e^{-20 \mu^{1/2}} \Oc(r^{-3} \omega^{-\frac{1}{4}})e^{\mu^{1/2} r },
\]
by Lemma \ref{Le:FundSystemCenterRepresentations} and Lemma \ref{Le:FundSystemCenterHankel}. Since $\chi_H(r,\omega) = \Oc(r^0 \omega^0)$,
\[ G_{2,2}(r,s,\omega) = \Oc(r^0 s^0 \omega^{0}) e^{-\mu^{1/2} r } + \Oc(r^0 
s^0 \omega^{-\frac12})  e^{-\mu^{1/2}(20 - r)}\]
and 
\[ r^m \partial_r^m G_{2,2}(r,s,\omega) = \Oc(r^m s^0  \omega^{\frac{m}{2}}) e^{-\mu^{1/2} r } + \Oc(r^m 
s^0 \omega^{\frac{m-1}{2}})  e^{-\mu^{1/2}(20 - r)}\]
such that  
\[|r^m \partial_r^m G_{2,2}(r,s,\omega)| \lesssim r^m \omega^{\frac{m}{2}} e^{-\mathrm{Re} \mu^{1/2} r }  
+ \omega^{\frac{m-1}{2}} e^{-12 \mathrm{Re} \mu^{1/2}} \lesssim 1 \lesssim 
r^{-1 + \delta} s^{-\delta}. \]
On the diagonal the kernel is supported for $\frac12 c \omega^{-\frac12} \leq r \leq c \omega^{-\frac12}$, which implies that 
$G_{2,2}(r,r,\omega) = \Oc(r^0 \omega^0)$. The same is true for derivatives of the kernel and we conclude that
\begin{align*}
r^m [T_{2,2}(\omega)f]^{(m)}(r)& = \sum_{k=0}^{m-1} \Oc(r^0 \omega^0) r^k f^{(k)}(r)   + \int_0^{r} r^m \partial_r^m  G_{2,2}(r,s,\omega)   f(s) ds.
\end{align*}

%-----------------T_{2,3}----------------------------------------------------
Next, we consider $T_{2,3}(\omega)$, i.e., the regime  $\frac12 c  \omega^{-\frac12} \leq  s \leq  r \leq  8$. Here, we use that 
 \begin{align}\label{Eq:v0Hankel}
 v_0(s,\omega) = \Oc(s^0 \omega^{-\frac{1}{4}}) e^{-\mu^{1/2} s } + \Oc(s^0 \omega^{-\frac{1}{4}}) e^{\mu^{1/2} s } , 
 \end{align}
by Lemma \ref{Le:FundSystemCenterRepresentations}. We split the kernel into $G_{2,3} = G^{-}_{2,3} + G^{+}_{2,3}$, where 
\begin{align*}
G^{-}_{2,3}(r,s,\omega)  & = \Oc(r^{-3} s^3 \omega^{-\frac{1}{2}}) e^{-\mu^{1/2}(20- r)} e^{-\mu^{1/2} s } 
+ \Oc(r^{-3} s^3 \omega^{-\frac{1}{2}}) e^{-\mu^{1/2}(10- r)} e^{-\mu^{1/2}(10 - s) }   \\
&  \phantom{!!!!!!!!!!!!!!!!!} + \Oc(r^{-3} s^3 \omega^0) e^{-\mu^{1/2} s } e^{-\mu^{1/2} r }, 
\end{align*}
and $G^{+}_{2,3}(r,s,\omega) = \Oc(r^{-3} s^3 \omega^0) e^{-\mu^{1/2} (r -s)}$. The operators $T^{\pm}_{2,3}(\omega)$ are defined accordingly. 
First, we observe that since $r^{-3} s^3 \leq 1$, 
\begin{align*}
|r^m \partial_r^{m} G^{-}_{2,3}(r,s,\omega)| \lesssim r^{m} \omega^{\frac{m}{2}}  ( e^{- 4 \mathrm{Re} \mu^{1/2}}  + e^{-\mathrm{Re} \mu^{1/2} r } ) \lesssim 1.
\end{align*}
On the diagonal, 
\begin{align*}
G^{-}_{2,3}(r,r,\omega)  & = \Oc(r^0 \omega^{-\frac{1}{2}}) e^{-20 \mu^{1/2}}
+ \Oc(r^0 \omega^{-\frac{1}{2}}) e^{-2 \mu^{1/2}(10- r)}  + \Oc( r^0 \omega^0) e^{-2 \mu^{1/2} r }.
\end{align*}
For derivatives of the boundary terms we exploit the exponential decay in $\omega$ and one can easily check that for $m=0,\dots,3$, 
\begin{align*}
r^m [T^{-}_{2,3}(\omega)f]^{(m)}(r) = \sum_{k=0}^{m-1} d_{mk}(r,\omega) r^{k} f^{(k)}(r) + \int_0^{r} r^m \partial_r^{m} G^{-}_{2,3}(r,s,\omega)  f(s) ds,
\end{align*}
with $|d_{mk}(r,\omega)| \lesssim 1$ and $|r^m\partial_r^{m} G^{-}_{2,3}(r,s,\omega)| \lesssim  r^{-1 + \delta} s^{-\delta} $.
Obviously, $|G^{+}_{2,3}(r,s,\omega)| \lesssim 1$, which settles things for $m=0$.
Furthermore, $G^{+}_{2,3}(r,r,\omega) = \Oc(r^0 \omega^0)$. Since
\[r \partial_r G^{+}_{2,3}(r,s,\omega) = \Oc(r^{-2} s^3\omega^{\frac12} )  e^{-\mu^{1/2} (r -s)},\]
an integration by parts yields
\begin{align*}
\int_0^{r} r \partial_r G^{+}_{2,3}(r,s,\omega) f(s) ds & = \Oc(r^0 \omega^0) r f(r)  \\
&  + \int_0^{r} e^{-\mu^{1/2} (r -s)} [\Oc(r^{-2} s^2 \omega^0) f(s) + \Oc(r^{-2} s^2 \omega^0) s f'(s) ] ds,
\end{align*}
which shows that  $(\cdot)[T_{2,3}(\omega)f]'$ is of the required form.  For $m=2,3$, we perform one, respectively, two
additional integrations by parts and conclude that 
\begin{align*}
r^m [T^{+}_{2,3}(\omega)f]^{(m)}(r) = \sum_{j=0}^{m-1} \tilde d_{mj}(r,\omega) r^{j+1} f^{(j)}(r) + 
\sum_{k=0}^{m}  \int_0^{r} j_{mk}(r,s,\omega)  s^{k} f^{(k)}(s) ds,
\end{align*}
where $|d_{mj}(r,\omega) | \lesssim 1$ and
\[ |j_{mk}(r,s,\omega) | \lesssim 1 \lesssim r^{-1 + \delta} s^{-\delta},\]
for all $0 < s \leq r$ and all $\omega \gg c^2$.

%-----------------T_{2,4}----------------------------------------------------
We proceed with $T_{2,4}(\omega)$, i.e.,
the case $0 < s \leq c \omega^{-\frac12}$ and $r \geq 4$.
Here, we study the expressions $(\cdot)^{n+2} [T_{2,4}(\omega)f]^{(n)}$ for $n=0,1$, and
$(\cdot)^{m} [T_{2,4}(\omega)f]^{(m)}$ for $m=2,3$. 
By Eq.~\eqref{Eq:Der_Expvpm} and Lemma \ref{Le:FundSystemCenterBessel} we have
\[ \partial^m_r G_{2,4}(r,s,\omega)  = e^{-10 \mu^{1/2}}  e^{\frac{r^2}{2}
- \mu \xi(r,\mu)} \Oc(r^{-3} s^0 \omega^{-\frac{3}{2}+\frac{m}{2}}) \Oc( \langle \omega^{-\frac12} r \rangle^{-\frac32 - m}).  \]
Note that the kernel and its derivatives vanish identically on the diagonal. Hence, 
\begin{align*}
r^{\alpha_m}  [T_{2,4}(\omega)f]^{(m)}(r) =   \int_0^{r}  r^{\alpha_m}  \partial^m_r G_{2,4}(r,s,\omega)  f(s) ds,
\end{align*}
for
\begin{align}\label{Eq:ResEst_T2_alpha}
\alpha_m =
\begin{cases}
	  2 \quad \text{ for } m = 0 ,2, \\
	   3 \quad \text{ for } m = 1 ,3.
\end{cases}
\end{align}
We show that
\begin{align}\label{Eq:ResEst_T24_1}
| r^{\alpha_m} \partial^m_r  G_{2,4}(r,s,\omega)|  \lesssim r^{-1 + \delta} s^{-\delta},
\end{align}
for all $r \geq s > 0$ and $\omega \gg c^2$, by distinguishing different regimes. For $4 \leq r < 10$, we use  
Eq.~\eqref{Eq:RepXi_large_alt} to get 
\begin{align*}
| r^{\alpha_m} \partial^m_r  G_{2,4}(r,s,\omega)| \lesssim  \omega^{-\frac{3}{2} + \frac{m}{2}}  
e^{-\mathrm{Re}\mu^{1/2}r}  e^{- \tilde \varphi(r, \omega) } \lesssim e^{- \tilde \varphi(4, \omega) }
\lesssim 1 \lesssim r^{-1 + \delta} s^{-\delta}.
\end{align*}
For $10 \leq r \leq \omega^{\frac12}$, Eq.~\eqref{Eq:RepXi_large} yields
\begin{align*}
|   r^{\alpha_m} \partial^m_r  G_{2,4}(r,s,\omega)| \lesssim  r^{-4 + {\alpha_m}} e^{-10\mathrm{Re}  \mu^{1/2}}  e^{- \varphi(r, \omega) }
\omega^{-1 + \frac{m}{2}} \lesssim r^{-1}  e^{- \varphi(10, \omega) } \lesssim  r^{-1 + \delta} s^{-\delta}.
\end{align*}
For $r >  \omega^{\frac12}$, and $m=0,2$,
\begin{align*}
| r^{2} \partial^m_r  G_{2,4}(r,s,\omega)| \lesssim r^{-1} 
\langle \omega^{-\frac12} r \rangle^{-\frac32 -\frac{b}{2} - m} e^{-10\mathrm{Re}  \mu^{1/2}}  e^{- \varphi(r, \omega) } \lesssim  r^{-1 + \delta} s^{-\delta},
\end{align*}
since $\langle \omega^{-\frac12} r \rangle^{-\frac32 -\frac{b}{2}} \lesssim 1$ if $b > -3$ and
$\langle \omega^{-\frac12} r \rangle^{-\frac32 -\frac{b}{2}} \lesssim r^{-\frac32 -\frac{b}{2}}$
if $-4 < b \leq -3$. For $m=1,3$, 
\begin{align*}
| r^{3} \partial^m_r  G_{2,4}(r,s,\omega)| & \lesssim 
\langle \omega^{-\frac12} r \rangle^{-\frac32 -\frac{b}{2} - m} e^{-10\mathrm{Re}  \mu^{1/2}}  e^{- \varphi(r, \omega) } \\
& 
\lesssim  r^{-1 + \delta}   (\omega^{-\frac12} r)^{1-\delta} \langle \omega^{-\frac12} r \rangle^{-1 + \delta} \omega^{\frac12} e^{-10\mathrm{Re}  \mu^{1/2}} \lesssim  
r^{-1 + \delta} s^{-\delta}.
\end{align*}
This implies Eq.~\eqref{Eq:ResEst_T24_1}.

%-----------------T_{2,5}----------------------------------------------------
Next, we consider $T_{2,5}(\omega)$. Here, the corresponding kernel is 
supported for $r \geq 4$, $\frac12 c \omega^{-\frac12} \leq s \leq 8$, $r \geq s$. 
We use Eq.~\eqref{Eq:v0Hankel} and write
$G_{2,5} =G^{-}_{2,5}  + G^{+}_{2,5}$, where 
\begin{align*}
G^{\pm}_{2,5}(r,s,\omega) = e^{-10 \mu^{1/2}}  \Oc(r^{-3} s^3 \omega^0) 
\Oc( \langle \omega^{-\frac12} r \rangle^{-\frac32})  e^{\frac{r^2}{2} - \mu \xi(r,\mu)} e^{\pm \mu^{1/2} s }.
\end{align*}
The operator $T_{2,5}(\omega) =T^{-}_{2,5}(\omega)  + T^{+}_{2,5}(\omega)$ splits accordingly. For the derivatives, we get 
\begin{align*}
\partial^m_r G^{\pm}_{2,5}(r,s,\omega) = e^{-10 \mu^{1/2}}  \Oc(r^{-3} s^3 \omega^{\frac{m}{2}})
\Oc( \langle \omega^{-\frac12} r \rangle^{-\frac32-m})  e^{\frac{r^2}{2} - \mu \xi(r,\mu)} e^{\pm \mu^{1/2} s }.
\end{align*}
Taking derivatives and using the definition Eq.~\eqref{Eq:ResEst_T2_alpha} yields
\begin{align}\label{Eq:ResEst_T25m}
r^{\alpha_m}[T^{-}_{2,5}(\omega)f]^{(m)}(r) = \sum_{k=0}^{m-1} p_{mk}(r,\omega) r^{k+1} f^{(k)}(r) 
+ \int_0^{r } r^{\alpha_m} \partial^m_r G^{-}_{2,5}(r,s,\omega) f(s) ds,
\end{align}
for $m=0,\dots,3$. We show that $|p_{mk}(r,\omega) | \lesssim 1$, and  
\[ |r^{\alpha_m} \partial^m_r G^{-}_{2,7}(r,s,\omega)| \lesssim r^{-1 + \delta} s^{-\delta}.\] 
For $4 \leq r < 10$,  we estimate
\[|r^{\alpha_m} \partial^m _r G^{-}_{2,5}(r,s,\omega)| \lesssim r^{-3 + \alpha_m} e^{-\mathrm Re \mu^{1/2} r } \omega^{\frac{m}{2}}  e^{- \tilde \varphi(4, \omega) }  
\lesssim 1 \lesssim r^{-1 + \delta} s^{-\delta},\]
for $10 \leq r \leq \omega^{\frac12}$, 
\[ 
 |r^{\alpha_m} \partial^m _r G^{-}_{2,5}(r,s,\omega)| \lesssim r^{-4 + \alpha_m} e^{-10 \mathrm Re \mu^{1/2}} \omega^{\frac{m+1}{2}}  e^{- \varphi(10, \omega) } 
 \lesssim r^{-1} \lesssim r^{-1 + \delta} s^{-\delta}, 
\]
and for $r > \omega^{\frac12}$, 
\[ 
 |r^{\alpha_m} \partial^m _r G^{-}_{2,5}(r,s,\omega)| \lesssim r^{-3 + \alpha_m} e^{-10 \mathrm Re \mu^{1/2}} \omega^{\frac{m}{2}}
 \langle \omega^{-\frac12} r \rangle^{-\frac32-\frac{b}{2} - m} e^{- \varphi(10, \omega) }  \lesssim r^{-1 + \delta} s^{-\delta}, 
\] 
by applying similar arguments as above. On the diagonal, the kernel and derivatives thereof are supported for $r \in [4,8]$, hence 
\[ | \partial_r^m G^{-}_{2,5}(r,s,\omega)|_{s=r}| \lesssim \omega^\frac{m}{2}e^{- 2 \mathrm{Re} \mu^{1/2}r} e^{- \tilde \varphi(4, \omega) }   \lesssim 1. \]
Adding the respective weights does not change anything. The exponential decay in $\omega$ guarantees that 
derivatives of kernels evaluated as $s=r$ behave well and we conclude that Eq.~\eqref{Eq:ResEst_T25m} holds. 
Next, we convince ourselves that  
\[| r^2 G^{+}_{2,5}(r,s,\omega) |  \lesssim r^{-1 + \delta} s^{-\delta} .\]
If fact, for $4 \leq r < 10$, 
\[ | r^2 G^{+}_{2,5}(r,s,\omega) |  \lesssim e^{-\mathrm{Re} \mu^{1/2}(r-s)} e^{- \tilde \varphi(4,\omega)} \lesssim 1, \]
since $r > s$. For $10 \leq r \leq  \omega^{\frac12}$,
\[ | r^2 G^{+}_{2,5}(r,s,\omega) |  \lesssim r^{-1} e^{-\mathrm{Re} \mu^{1/2}(10-s)} 
 e^{- \varphi(10,\omega)}   \lesssim r^{-1},\]
since $s \leq 8$ and for  $r >  \omega^{\frac12}$, 
\[ | r^2 G^{+}_{2,5}(r,s,\omega) |  \lesssim r^{-1} e^{-\mathrm{Re} \mu^{1/2}(10-s)} 
 e^{- \varphi(10,\omega)}   \langle \omega^{-\frac12} r \rangle^{-\frac32-\frac{b}{2}}  \lesssim r^{-1+\delta}. \]
For higher derivatives, we integrate by parts. For example, 
\[ \int_0^{r} r^3 \partial_r G^{+}_{2,5}(r,s,\omega) f(s)ds = p_0(r,\omega) f(r) + 
\int_0^{r} [ h_1(r,s,\omega) f(s) + h_2(r,s,\omega)  sf'(s)]ds  \]
with the support of 
$p_0(\cdot,\omega)$ contained in the interval $[4,8]$ and 
\[ p_0(r,\omega) = \Oc(r^0 \omega^0) e^{-10 \mu^{1/2}} e^{\mu^{1/2} r - \mu \xi(r,\mu)}.\]
By Eq.~\eqref{Eq:RepXi_large_alt}, we infer that 
$|p_0(r,\omega) | \lesssim 1$. 
The integral kernels are of the form 
\[h_j(r,s,\omega) = e^{- \mu^{1/2}(10 - s)} e^{\frac{r^2}{2} - \mu \xi(r,\mu)} \Oc(r^0 s^2 \omega^0) \Oc( \langle \omega^{-\frac12} r \rangle^{-\frac52}) , \]
for $j=1,2$.
By using similar argument as above it is not hard to see that 
$|h_j(r,s,\omega)| \lesssim r^{-1+\delta} s^{-\delta}$. In summary, we get that 
\begin{align*}
r^3[T^+_{2,5}(\omega)f]'(r) =  \tilde p_0(r,\omega) r f(r) +  \int_0^{r} [ h_1(r,s,\omega) f(s) + h_2(r,s,\omega)  sf'(s)]ds  
\end{align*}
where $\tilde p_0(r,\omega) =\Oc(r^0 \omega^0) e^{-10 \mu^{1/2}} e^{\mu^{1/2} r - \mu \xi(r,\mu)}$ such that  $|\tilde p_0(r,\omega)| \lesssim 1$.
Derivatives of second and third order can be handled by exploiting Eq.~\eqref{Der:ExpHankelPlusWeberMin} to control derivatives of boundary terms
and using integration by parts.
With this, we infer that for $m=2,3$, 
\begin{align*}
r^m[T^{+}_{2,5}(\omega)f]^{(m)}(r) = \sum_{j=0}^{m-1} \tilde p_{mj}(r,\omega) r^j f^{(j)}(r) 
+  \sum_{k=0}^{m} \int_0^{r} \tilde h_{mk}(r,s,\omega)  s^{k} f^{(k)}(s) ds,
\end{align*}
where $|\tilde p_{mj}(r,\omega)| \lesssim 1$,
\[ \tilde h_{mk}(r,s,\omega) = e^{- \mu^{1/2}(10 - s)} e^{\frac{r^2}{2} - \mu \xi(r,\mu)}
  \Oc(r^{-3+m} s^{3 -m}\omega^0) \Oc( \langle \omega^{-\frac12} r \rangle^{-\frac32-m }),
\]
such that $|\tilde h_{mk}(r,s,\omega) | \lesssim r^{-1+\delta} s^{-\delta}$.

%-----------------T_{2,6}----------------------------------------------------
Next, we have $4 \leq s \leq r \leq 40$. Again, we express $v_0$ in terms of Hankel functions and $v_{-}$ in terms of Weber functions 
and write $G_{2,6} = G^+_{2,6} +G^-_{2,6}$, $T_{2,6}(\omega) = T^{-}_{2,6}(\omega) + T^{+}_{2,6}(\omega)$, where  
\begin{align*}
G^{\pm}_{2,6}(r,s,\omega) =  \chi_{W_1}(r) \chi_W(s)  e^{-10 \mu^{1/2}}  \Oc(r^{0} s^0 \omega^0)  e^{- \mu \xi(r,\mu)} e^{\pm \mu^{1/2} s}.
\end{align*}
In this regime, weights can be added arbitrarily in the estimates, which simplifies matters. 
We claim that for $m = 0, \dots, 3$, 
\begin{align*}
r^m[T_{2,6}(\omega)f]^{(m)}(r) = \sum_{k=0}^{m-1} \hat p_{mk}(r,\omega) r^k f^{(k)}(r) 
+   \sum_{j=0}^{m}  \int_0^{r} \hat h_{mj}(r,s,\omega)   s^{j} f^{(j)}(s) ds,
\end{align*}
where $| \hat h_{mj}(r,s,\omega)   | \lesssim r^{-1+\delta} s^{-\delta}$, $|\hat p_{mk}(r,\omega)| \lesssim 1$. For example, we immediately get that 
\[|G^{\pm}_{2,6}(r,s,\omega)| \lesssim e^{- \mathrm{Re} \mu^{1/2} ( r - s)} + e^{- \mathrm{Re} \mu^{1/2}r} e^{- \mathrm{Re} \mu^{1/2} s}\lesssim 1,\]
which implies the claim for $m = 0$. For higher derivatives, we use that 
\[ \partial^m_r G^{-}_{2,6}(r,s,\omega) = e^{-10 \mu^{1/2}} \Oc(r^0 s^0 \omega^{\frac{m}{2}}) e^{- \mu \xi(r,\mu)} e^{- \mu^{1/2} s}, \]
which yields $|r^m \partial^m_r G^{-}_{2,6}(r,s,\omega)| \lesssim  \omega^{\frac{m}{2}} e^{- 8 \mathrm{Re} \mu^{1/2}  } \lesssim 1$.
On the diagonal $G^{-}_{2,6}$ and its derivatives behave well and it is obvious that $T^{-}_{2,6}(\omega)$ is of the required form. 
To handle $T^{+}_{2,6}(\omega)$, we use integration by parts and exploit Eq.~\eqref{Der:ExpHankelPlusWeberMin} to control derivatives of boundary terms.
This is analogous to the above calculations and we omit further details.

%-----------------T_{2,7}----------------------------------------------------
It is left to discuss $T_{2,7}(\omega)$, i.e., the case $s \geq 4$, $r \geq 20$, $r \geq s$. 
According to the proof of Lemma \ref{Le:Rep_vmv0} we can write $v_0$ as  
\begin{align*}
v_0(s,\omega) = e^{10 \mu^{1/2}} \Oc(\omega^0) v_{+}(s,\omega) +
e^{-10 \mu^{1/2}} \Oc(\omega^0) v_{-}(s,\omega) +   \Oc(\omega^0)  W_{+}(\omega)  v_{-}(s,\omega),
\end{align*}
where $W_{+}(\omega) = \Oc(\omega^{-\frac12}) e^{10 \mu^{1/2}}$.
In the following, we exploit the fact that by definition 
\begin{align}\label{Eq:WronskianPlus}
W_{+}(\omega) =  W(\tilde v_{+}(\cdot,\omega) , v_+(\cdot,\omega))  = 
\Oc(s^0 \omega^{-\frac12}) e^{\mu^{1/2}s}  e^{\mu \xi(s,\mu)},
\end{align}
for $s \in [4,10]$. We write $G_{2,7} = G^{+}_{2,7}  +   G^{-}_{2,7}$
where 
\begin{align*}
G^{+}_{2,7}(r,s,\omega)  & = \chi_{W_2}(r) \chi_W(s)  \Oc(r^{-3} \omega^0)   e^{\frac{r^2}{2}- \mu \xi(r,\mu)}  \Oc( \langle \omega^{-\frac12} r \rangle^{-\frac32 }) \\
 & \phantom{!!!!!!!!!!!!!!!}  \times e^{- \frac{s^2}{2} + \mu \xi(s,\mu)} \Oc( s^3 \omega^0)  \Oc( \langle \omega^{-\frac12} s \rangle^{-\frac12}) ,
\end{align*} 
$G^{-}_{2,7}(r,s,\omega)  = K^{-}_1(r,s,\omega) + K^{-}_2(r,s,\omega)$ with 
\begin{align*}
 K^{-}_1(r,s,\omega) & =  \chi_{W_2}(r) \chi_W(s)   \Oc(r^{-3}\omega^0)  e^{\frac{r^2}{2}- \mu \xi(r,\mu)}  \Oc( \langle \omega^{-\frac12} r \rangle^{-\frac32 }) \\
   & \phantom{!!!!!} \times   e^{-\frac{s^2}{2} -\mu \xi(s,\mu)} e^{-20 \mu^{1/2}}    \Oc(s^3 \omega^0)   \Oc( \langle \omega^{-\frac12} s \rangle^{-\frac12}),
\end{align*} 
\begin{align*}
 K^{-}_2(r,s,\omega) & =  \chi_{W_2}(r) \chi_W(s)   \Oc(r^{-3}\omega^0)  e^{\frac{r^2}{2}- \mu \xi(r,\mu)}  \Oc( \langle \omega^{-\frac12} r \rangle^{-\frac32 })\\
&  \phantom{!!!!!} \times  e^{-\frac{s^2}{2} -\mu \xi(s,\mu)}  e^{-10 \mu^{1/2}}   W_{+}(\cdot,\omega) \Oc(s^3 \omega^0)   \Oc( \langle \omega^{-\frac12} s \rangle^{-\frac12}).
\end{align*} 
%\begin{align*}
%\tilde G^{-}_{2,9}(r,s,\omega)  & = \chi_{W_2}(r) \chi_W(s)   \Oc(r^{-3} \omega^0)  e^{-10 \mu^{1/2}} e^{\frac{r^2}{2}- \mu \xi(r,\mu)} \Oc( \langle \omega^{-\frac12} r \rangle^{-\frac32 })  \\
%& \phantom{!!!!!!!!!!!!!!!!!} \times \Oc(s^3 \omega^0)   W_{+}(\cdot,\omega)  e^{-\frac{s^2}{2} -\mu \xi(s,\mu)}  \Oc( %\langle \omega^{-\frac12} s \rangle^{-\frac12}).
%\end{align*} 
Accordingly, we write $T_{2,7}(\omega)  =T^{+}_{2,7}(\omega)  + T^{-}_{2,7}(\omega)$, where 
\[ [T^{\pm}_{2,7}(\omega)f](r) = \int_0^{r} G^{\pm}_{2,7}(r,s,\omega) f(s) ds. \]
First, we consider
\[ r^2 [T^{-}_{2,7}(\omega)f](r) = \int_0^{r} r^2 G^{-}_{2,7}(r,s,\omega) f(s) ds, \]
and convince ourselves that $|r^2 G^{-}_{2,7}(r,s,\omega)| \lesssim  r^{-1+\delta} s^{-\delta}$.  By Eq.~\eqref{Eq:WronskianPlus},
\[W_{+}(s,\omega)  e^{-\mu \xi(s,\mu)}  =   \Oc(s^0 \omega^{-\frac12}) e^{\mu^{1/2}s},\]
for $4 \leq s \leq 10$. Using Eq.~\eqref{Eq:RepXi_large} and Eq.~\eqref{Eq:RepXi_large_alt}, we get that for $4 \leq s \leq 10$, $r \geq 20$,
\begin{align*}
 |r^2 G^{-}_{2,7}(r,s,\omega)|& \lesssim r^{-1} \langle \omega^{-\frac12} r \rangle^{-\frac32 - \frac{b}{2}}
 \big( e^{-20\mathrm{Re} \mu^{1/2}}   e^{-\mathrm{Re}  \mu \xi(s,\mu)} + e^{- \mathrm{Re} \mu^{1/2}(10 -s)} \big) \\
& \lesssim r^{-1} \langle \omega^{-\frac12} r \rangle^{-\frac32 - \frac{b}{2}} \big(
e^{-10\mathrm{Re} \mu^{1/2}} e^{-\mu^{1/2}s}   + e^{- \mathrm{Re} \mu^{1/2}(10 -s)} \big) \lesssim r^{-1+\delta} s^{-\delta}.
\end{align*}
In the last step, we used that for $r > \omega^{\frac12}$, 
$r^{-1} \langle \omega^{-\frac12} r \rangle^{-\frac32 - \frac{b}{2}} \lesssim r^{-1}$
if  $b \geq -3$, and that 
$r^{-1} \langle \omega^{-\frac12} r \rangle^{-\frac32 - \frac{b}{2}} \lesssim r^{-\frac{5}{2}-\frac{b}{2}}$
if $- 4  < b < -3$. For $s > 10$, $r \geq 20$, we use the value of the Wronskian and Eq.~\eqref{Eq:RepXi_large} to estimate 
\begin{align*}
 |r^2  G^{-}_{2,7}(r,s,\omega)|& \lesssim r^{-1} s^3 e^{-s^2} \langle \omega^{-\frac12} r \rangle^{-\frac32 - \frac{b}{2}}
 \langle \omega^{-\frac12} s \rangle^{-\frac12 - \frac{b}{2}}   \\
 & \lesssim r^{-1}\langle \omega^{-\frac12} r \rangle^{-\frac32 - \frac{b}{2}} 
 s^{- \delta} s^{\frac{9}{2} + \delta} e^{-s^2} \lesssim  r^{-1+\delta} s^{-\delta},
\end{align*}
since $\langle \omega^{-\frac12} s \rangle^{-\frac12 - \frac{b}{2}}  \lesssim s^{\frac{3}{2}}$ for $b > -4$. 
Next, we show that 
\begin{align}\label{Eq:T29m1}
r^3 [T^{-}_{2,7}(\omega)f]'(r) =    q^{-}_0(r,\omega) r f(r) + \int_0^{r} r^3 [k^{-}_0(r,s,\omega)  s f(s) +k^{-}_1(r,s,\omega)  s^2 f'(s)] ds 
\end{align}
where $|  q^{-}_0(r,\omega)| \lesssim 1$ and  $|r^3 k^{-}_j(r,s,\omega)| \lesssim r^{-1+\delta} s^{-\delta}$ for $j=0,1$.
First, we observe that 
\begin{align*}
R^{-}_{2,7}(r,\omega)  & :=G^{-}_{2,7}(r,r,\omega) = e^{- 2\mu \xi(r,\mu)} \Oc( \langle \omega^{-\frac12} r \rangle^{-2})[
e^{-20 \mu^{1/2}}  \Oc(r^0 \omega^0)   + \Oc(r^0 \omega^{-\frac12})],
\end{align*}
where $R^{-}_{2,7}(\cdot,\omega)$ has support for $20 \leq r < \infty$. Hence,
\[|r^2 R^{-}_{2,7}(r,\omega)| \lesssim  r^2  e^{-r^2}  \langle \omega^{-\frac12} r \rangle^{-2-b} \lesssim 1.\] 
Furthermore, for $r \geq 20$,  $4 \leq s \leq 10$,  
\begin{align*}
|r^3 \partial_r K^{-}_1(r,s,\omega)|  & \lesssim \omega^{\frac12} \langle \omega^{-\frac12} r \rangle^{-\frac52-\frac{b}{2} }
e^{-10\mathrm{Re} \mu^{1/2}} e^{-\mathrm{Re} \mu^{1/2}s} e^{-\varphi(10, \omega)} e^{-\tilde \varphi(4, \omega)}   \\
&  \lesssim  r^{-1+ \delta}  (\omega^{-\frac12} r)^{1 - \delta} \langle \omega^{-\frac12} r \rangle^{-1+\delta} 
\omega e^{-10\mathrm{Re} \mu^{1/2}} \lesssim r^{-1+ \delta} s^{-\delta},
\end{align*}
and for $s > 10$, 
\begin{align*}
|r^3 s^{-1} \partial_r K^{-}_1(r,s,\omega)|  & \lesssim \omega^{\frac12}   s^2 e^{-s^2} e^{-20\mathrm{Re} \mu^{1/2}}  
\langle \omega^{-\frac12} r \rangle^{-\frac52-\frac{b}{2}} \langle \omega^{-\frac12} s \rangle^{-\frac12 - \frac{b}{2}}  \lesssim r^{-1+ \delta} s^{-\delta}.
\end{align*}
For $K^{-}_2$ we integrate by parts. Using that $[r + \mu \partial_r \xi(r,\mu)]^{-1} = \mu^{-\frac12} \Oc( \langle \omega^{-\frac12} r \rangle^{-1})$ we obtain
\begin{align}\label{Eq:g2m} 
\begin{split}
\int_0^{r} r^3 \partial_r K^{-}_2(r,s,\omega) f(s)ds  & = \Oc(r^2 \omega^{-\frac12}) e^{-2 \mu \xi(r,\mu)}  \Oc( \langle \omega^{-\frac12} r \rangle^{-2}) r f(r) \\
  & \phantom{!!!!!!!!!!!!!!!!!!!!!}  \int_0^r [  r^3 \tilde k^{-}_0(r,s,\omega)  s f(s) +  r^3 \tilde k^{-}_1(r,s,\omega) s^2 f'(s)] ds,
  \end{split}
\end{align}
where for $j=0,1$, 
\[\tilde k^{-}_j(r,s,\omega) = \Oc(r^{-3} s \omega^0)    \Oc( \langle \omega^{-\frac12} r \rangle^{-\frac52 })
 \Oc( \langle \omega^{-\frac12} s \rangle^{-\frac32})
e^{-10 \mu^{1/2}}e^{\frac{r^2}{2}- \mu \xi(r,\mu)}  e^{-\frac{s^2}{2} -\mu \xi(s,\mu)}  W_{+}(\cdot,\omega). 
\]
The boundary term can be estimated easily.  
For the integral kernels we proceed as above: For $4 \leq s \leq 10$, we apply Eq.~\eqref{Eq:WronskianPlus} 
to get
\[ |r^3 \tilde k^{-}_j(r,s,\omega)| \lesssim  \langle \omega^{-\frac12} r \rangle^{-\frac52-\frac{b}{2}} e^{-10\mathrm{Re} \mu^{1/2}} e^{-\mathrm{Re} \mu^{1/2}s}
\lesssim r^{-1+ \delta} s^{-\delta},\]
and for $s >10$, we use the value of $W_+$ together with 
 Eq.~\eqref{Eq:RepXi_large} to obtain
\[  |r^3 \tilde k^{-}_j(r,s,\omega)|  \lesssim \omega^{-\frac12} \langle \omega^{-\frac12} r \rangle^{-\frac52-\frac{b}{2}} 
\langle \omega^{-\frac12} s \rangle^{-\frac32 - \frac{b}{2}} 
e^{-s^2} \lesssim r^{-1 + \delta} s^{-\delta}. \]
In the last step, we used that for $r > \omega^{\frac12}$, 
\[ \omega^{-\frac12}  \langle \omega^{-\frac12} r \rangle^{-\frac52 - \frac{b}{2} }  \lesssim 
 \omega^{-\frac12}  \langle \omega^{-\frac12} r \rangle^{-1 + \delta} 
 \lesssim r^{-1 + \delta} \omega^{-\frac{\delta}{2}}  \lesssim  r^{-1 + \delta}. \]
This shows that Eq.~\eqref{Eq:T29m1} holds. Next, we show that for $m=2,3$, 
\begin{align}\label{Eq:T29m2}
r^{m} [T^{-}_{2,7}(\omega)f]^{(m)}(r) =   \sum_{k=0}^{m-1} q^{-}_{mk}(r,\omega)  r^{k} f^{(k)}(r) +  \sum_{j=0}^{m} \int_0^{r} 
 k^{-}_{mj} (r,s,\omega)  s^{j} f^{(j)} (s) ds,
\end{align}
with $|q^{-}_{mk}(r,\omega) | \lesssim 1$ and $|k^{-}_{mj} (r,s,\omega) | \lesssim r^{-1 + \delta} s^{-\delta}$. 
We start with $m=2$ using Eq.~\eqref{Eq:T29m1}. To estimate derivatives of the boundary terms note that for example
\[r^2 \partial_r R^{-}_{2,7}(r,\omega) =  e^{- 2\mu \xi(r,\mu)}(\mu + r^2)^{\frac12}  \Oc(r^2 \omega^0) 
\Oc( \langle \omega^{-\frac12} r \rangle^{-2}).\] 
For $20\ \leq r \leq \omega^{\frac12}$, 
\[  \omega^{\frac32} e^{- 2 \mathrm{Re} \mu \xi(r,\mu)} \lesssim 
\omega^{\frac32}  e^{- \mathrm{Re} \mu^{1/2}(r -10)} e^{-\tilde \varphi(10,\omega)} \lesssim 1, \]
whereas for $r  > \omega^{\frac12}$,  $ r^3 e^{- 2 \mathrm{Re} \mu \xi(r,\mu)}  \lesssim r^3  e^{-r^2} 
\langle \omega^{-\frac12} r \rangle^{-b} e^{- \varphi(10,\omega)} \lesssim 1.$
In fact, all boundary terms that occur in Eq.~\eqref{Eq:T29m2} are of a similar form and can be estimated analogously.
Furthermore, we have 
\begin{align*}
| r^2 \partial^2_r K^{-}_1(r,s,\omega)| \lesssim   r^{-1} \omega \langle \omega^{-\frac12} r \rangle^{-\frac72 - \frac{b}{2}} 
e^{-20  \mathrm{Re}  \mu^{1/2}} e^{-\frac{s^2}{2} - \mathrm{Re}  \mu \xi(s,\mu)} \lesssim r^{-1 + \delta} s^{-\delta},
\end{align*} 
by applying similar arguments as above. An integration by parts shows that   
\begin{align*}
& \int_0^r r^2 \partial_r [ \tilde k^{-}_0(r,s,\omega) f(s) +   \tilde k^{-}_1(r,s,\omega) s f'(s)] ds \\
& \phantom{!!!!!!!!!!!!!!!!!!!!!}    = 
\sum_{k\in \{0,1\}} \tilde q^{-}_k(r,\omega)  r^{k} f^{(k)}(r)   + \sum_{j=0}^{2} \int_0^{r} \hat k^{-}_j(r,s,\omega)  s^{j} f^{(j)} (s) ds,
\end{align*}
where $|\hat q^{-}_k(r,\omega)| \lesssim 1$ and  
\[ \hat k^{-}_j(r,s,\omega) = \Oc(r^{-3} s^0 \omega^0)    \Oc( \langle \omega^{-\frac12} r \rangle^{-\frac72 })
 \Oc( \langle \omega^{-\frac12} s \rangle^{-\frac52})
e^{-10 \mu^{1/2}}e^{\frac{r^2}{2}- \mu \xi(r,\mu)}  e^{-\frac{s^2}{2} -\mu \xi(s,\mu)}  W_{+}(\cdot,\omega). 
\]
By the same arguments as above, we obtain that $|\hat k^{-}_j(r,s,\omega)| \lesssim r^{-1 + \delta} s^{-\delta}$.
The case $m=3$ is treated similarly. 
Finally, we investigate $T^{+}_{2,7}(\omega)$. 
We integrate by parts using that $[r -  \mu \partial_r \xi(r,\mu)]^{-1} = \mu^{-\frac12} \Oc(\langle  \omega^{-\frac12} r \rangle)$ to obtain
\begin{align}\label{Eq:T29p1}
\begin{split}
r^2[T^{+}_{2,7}(\omega)f](r) & = \int_0^{r} r^2 G^{+}_{2,7}(r,s,\omega) f(s)ds \\
 \phantom{!!!!!!} & =  r q_0^{+}(r,\omega) r f(r)  + \int_0^r  r^2 [ k_0^{+}(r,s,\omega)   s f(s) + k^{+}_1(r,s,\omega) s^2 f'(s)] ds, 
\end{split}
\end{align}
where the kernels are supported for $r \geq 20$ and $s \geq 4$. 
The boundary term is of the form
\[  q_0^{+}(r,\omega) =   \chi_{W_2}(r) \Oc(r^{0} \omega^{-\frac12})
\Oc(\langle r \omega^{-\frac12} \rangle^{-1}) = \Oc(r^{-1} \omega^{0}), \]
such that $| r q_0^{+}(r,\omega)| \lesssim 1$. 
The integral kernels are given by 
\[ k_j^{+}(r,s,\omega) = \Oc(r^{-3} s \omega^{-\frac12}) e^{\frac{r^2}{2}- \mu \xi(r,\mu)} 
\Oc( \langle \omega^{-\frac12} r \rangle^{-\frac32 })  e^{- \frac{s^2}{2} + \mu \xi(s,\mu)} 
\Oc( \langle \omega^{-\frac12} s \rangle^{\frac12}),\]
for $j=0,1$. By Lemma \ref{Le:RepXiQ}, for $4 \leq s \leq 10$, 
\begin{align*}
| r^2 k_j^{+}(r,s,\omega) |  & \lesssim \omega^{-\frac12} r^{-1}  
\langle \omega^{-\frac12} r \rangle^{-\frac32 - \frac{b}{2}} e^{-\varphi(10,\omega)} 
e^{-\mathrm{Re}  \mu^{1/2}(10 - s )} e^{\tilde \varphi(10,\omega)}  \\
& \lesssim \omega^{-\frac12} r^{-1} 
\langle \omega^{-\frac12} r \rangle^{-\frac32 - \frac{b}{2}} \lesssim r^{-1 + \delta} s^{-\delta}.
\end{align*}
For $s > 10$ we use that $r^{-1} s \leq \langle \omega^{-\frac12} r \rangle^{-1} \langle \omega^{-\frac12} s \rangle$ for $s \leq  r$ to estimate
\begin{align*}
| r^2 k_j^{+}(r,s,\omega) | &  \lesssim \omega^{-\frac12}  r^{-1}  s 
\langle \omega^{-\frac12} r \rangle^{-\frac32 - \frac{b}{2}} \langle \omega^{-\frac12} s \rangle^{\frac12 
+ \frac{b}{2}}  e^{-[\varphi(r,\omega)-\varphi(s,\omega)]} \\
&   \lesssim 
\omega^{-\frac12}  \langle \omega^{-\frac12} r \rangle^{-\frac52 - \frac{b}{2}} \langle \omega^{-\frac12} s \rangle^{\frac32 + \frac{b}{2}} 
\lesssim r^{-1 + \delta} s^{-\delta},
\end{align*}
where the last step follows from the fact that 
\[\omega^{-\frac12}  \langle \omega^{-\frac12} r \rangle^{-\frac52 - \frac{b}{2}} 
\langle \omega^{-\frac12} s \rangle^{\frac32 + \frac{b}{2}}  \leq \omega^{-\frac12}  \langle \omega^{-\frac12} r \rangle^{-1} \lesssim r^{-1}, 
\]
for  $b \geq -3$, whereas 
\begin{align*}
\omega^{-\frac12}  \langle \omega^{-\frac12} r \rangle^{-\frac52 - \frac{b}{2}} 
\langle \omega^{-\frac12} s \rangle^{\frac32 + \frac{b}{2}}  & \leq r^{-\frac52 - \frac{b}{2}} s^{\frac32 + \frac{b}{2}}
(\omega^{-\frac12} r )^{\frac52 + \frac{b}{2}}   \langle \omega^{-\frac12} r \rangle^{-\frac52 - \frac{b}{2}}
(\omega^{-\frac12} s )^{-\frac32 - \frac{b}{2}}\langle \omega^{-\frac12} s \rangle^{\frac32 + \frac{b}{2}} \\
&   \leq r^{-\frac52 - \frac{b}{2}} s^{\frac32 + \frac{b}{2}},
\end{align*}
for $- 4 < b < -3$. 
This proves that
\[ | r^2 k_j^{+}(r,s,\omega) | \lesssim r^{-1 + \delta} s^{-\delta},\]
for all $r \leq s > 0$ and $\omega \gg c^2$. 

Starting from Eq.~\eqref{Eq:T29p1}, we show that  
\begin{align}\label{Eq:T29p2}
r^3 [T^{+}_{2,7}(\omega)f]'(r) =  \sum_{j\in \{0,1\}} \tilde q_j^{+}(r,\omega) r^{j+1} f^{(j)}(r)   
+  \sum_{n=0}^{2} \int_0^{r} r^3\tilde k^{+}_n(r,s,\omega) s^{n+1} f^{(n)} (s) ds  
\end{align}
with $|\tilde q_j^{+}(r,\omega)  | \lesssim 1$ and $|\tilde k^{+}_n(r,s,\omega)| \lesssim r^{-1+\delta} s^{-\delta}$.
In fact, $\partial_r q_0^{+}(r,\omega)  = \Oc(r^{-2} \omega^0)$ and
\[ k^{+}_j(r,r,\omega) = \Oc(r^{-2} \omega^{-\frac12}) \Oc( \langle  \omega^{-\frac12} r \rangle^{-1} ) = \Oc(r^{-3} \omega^0), \]
for $j=0,1$. Moreover,
\[  \partial_r k^{+}_j(r,s,\omega) =  \Oc(r^{-3} s \omega^{0}) e^{\frac{r^2}{2}- \mu \xi(r,\mu)} 
\Oc( \langle \omega^{-\frac12} r \rangle^{-\frac52 })  e^{- \frac{s^2}{2} + \mu \xi(s,\mu)} 
\Oc( \langle \omega^{-\frac12} s \rangle^{\frac12}),\]
and an integration by parts yields
\begin{align*}
\sum_{j\in \{0,1\}} & \int_0^{r} \partial_r k^{+}_j(r,s,\omega)  s^{j+1} f^{(j)}(s) ds  \\
&   = 
\sum_{j\in \{0,1\}} \Oc(r^{-3} \omega^0) r^{j+1} f^{(j)}(r)
+  \sum_{n=0}^{2} \int_0^{r} \tilde k^{+}_n(r,s,\omega) s^{n+1} f^{(n)} (s) ds 
\end{align*}
where 
\[ \tilde k^{+}_n(r,s,\omega) = \Oc(r^{-3} s^{0} \omega^{-\frac12}) e^{\frac{r^2}{2}- \mu \xi(r,\mu)} 
\Oc( \langle \omega^{-\frac12} r \rangle^{-\frac52 })  e^{- \frac{s^2}{2} + \mu \xi(s,\mu)} 
\Oc( \langle \omega^{-\frac12} s \rangle^{\frac32}).\]
As above, for $4 \leq s \leq 10$, 
\begin{align*}
|r^3 \tilde k^{+}_n(r,s,\omega)| \lesssim  \omega^{-\frac12} \langle \omega^{-\frac12} r \rangle^{-\frac52-\frac{b}{2}} 
e^{-\mathrm{Re}  \mu^{1/2}(10 - s )} \lesssim r^{-1+\delta} s^{-\delta},
\end{align*}
since for $r > \omega^{\frac12}$, 
$r^{-1}  (r \omega^{-\frac12})\langle \omega^{-\frac12} r \rangle^{-\frac52-\frac{b}{2}}  \leq r^{-1}$ if $b  \geq -3$ and 
$ \omega^{-\frac12} \langle \omega^{-\frac12} r \rangle^{-\frac52-\frac{b}{2}} \leq r^{-\frac52-\frac{b}{2}}$ if $-4 < b < -3$. 
For $s > 10$, 
\begin{align*}
|r^3 \tilde k^{+}_n(r,s,\omega)| \lesssim  \omega^{-\frac12} \langle \omega^{-\frac12} r \rangle^{-\frac52-\frac{b}{2}} 
\langle \omega^{-\frac12} s \rangle^{\frac32+\frac{b}{2}}  \lesssim r^{-1+\delta} s^{-\delta}.
\end{align*}
This implies Eq.~\eqref{Eq:T29p2}. With this, it is not hard to see that for higher derivatives we have
\begin{align*}
r^2 [T^{+}_{2,7}(\omega)f]''(r) =  \sum_{j= 0}^2 c_j(r,\omega)  r^{j} f^{(j)}(r)   +  \sum_{k=0}^{3} \int_0^{r} i_k(r,s,\omega) s^{k} f^{(k)} (s) ds  
\end{align*}
with $|c_j(r,\omega)   | \lesssim 1$ and 
\[i_k(r,s,\omega) = \Oc(r^{-1} s^{0} \omega^{-\frac12}) e^{\frac{r^2}{2}- \mu \xi(r,\mu)} 
\Oc( \langle \omega^{-\frac12} r \rangle^{-\frac72 })  e^{- \frac{s^2}{2} + \mu \xi(s,\mu)} 
\Oc( \langle \omega^{-\frac12} s \rangle^{\frac52}).  \]
such that $|i_k(r,s,\omega)| \lesssim r^{-1+\delta} s^{-\delta}$. Integrating by parts once more shows that 
\begin{align*}
\begin{split}
r^3 [T^{+}_{2,7}(\omega)f]'''(r)  & =   \tilde c_0(r,\omega)  f(r) + \sum_{j= 1}^3 \tilde  c_j(r,\omega) r^{j-1} f^{(j)}(r) \\
& +
\int_0^{r} \tilde i_0(r,s,\omega) f(s) ds +  \sum_{k=1}^{4}  \int_0^{r} \tilde i_k(r,s,\omega) s^{k-1} f^{(k)} (s) ds,
\end{split}
\end{align*}
where $|\tilde  c_j(r,\omega)| \lesssim 1$ and $|\tilde i_k(r,s,\omega)| \lesssim r^{-1+\delta} s^{-\delta}$ for $k = 0, \dots,4$, all $0 <s \leq r$ and
all $\omega \gg c^2$. Summing up the individual contributions shows that $\mc T_2(\omega)f$ as well as its weighted derivatives are of the claimed form.
\end{proof}

\pagestyle{plain}
\bibliography{references}
\bibliographystyle{plain}

\end{document}